\documentclass[12pt,vatola]{article}

\usepackage{graphicx}

\def\c{\centerline}

\def\no{\noindent}

\begin{document}

\c{\Large\bf Smarandache Multi-Space Theory(II) }\vskip 5mm

\hskip 70mm {\it -Multi-spaces on graphs}\vskip 10mm

\c{Linfan Mao}\vskip 3mm

\c{\small Academy of Mathematics and System Sciences}

\c{\small Chinese Academy of Sciences, Beijing 100080}

\c{\small maolinfan@163.com}

\vskip 10mm

\begin{minipage}{130mm}

\no{\bf Abstract.} {\small A Smarandache multi-space is a union of
$n$ different spaces equipped with some different structures for
an integer $n\geq 2$, which can be both used for discrete or
connected spaces, particularly for geometries and spacetimes in
theoretical physics. This monograph concentrates on characterizing
various multi-spaces including three parts altogether. The first
part is on {\it algebraic multi-spaces with structures}, such as
those of multi-groups, multi-rings, multi-vector spaces,
multi-metric spaces, multi-operation systems and multi-manifolds,
also multi-voltage graphs, multi-embedding of a graph in an
$n$-manifold,$\cdots$, etc.. The second discusses {\it Smarandache
geometries}, including those of map geometries, planar map
geometries and pseudo-plane geometries, in which the {\it Finsler
geometry}, particularly the {\it Riemann geometry} appears as a
special case of these Smarandache geometries. The third part of
this book considers the {\it applications of multi-spaces to
theoretical physics}, including the relativity theory, the
M-theory and the cosmology. Multi-space models for $p$-branes and
cosmos are constructed and some questions in cosmology are
clarified by multi-spaces. The first two parts are relative
independence for reading and in each part open problems are
included for further research of interested readers.}

\vskip 3mm \no{\bf Key words:} {\small  graph, multi-voltage
graph, Cayley graph of a multi-group, multi-embedding of a graph,
map, graph model of a multi-space, graph phase.}

 \vskip 3mm \no{{\bf
Classification:} AMS(2000) 03C05,05C15,51D20,51H20,51P05,83C05,
83E50}
\end{minipage}

\newpage

\no{\large\bf Contents}\vskip 8mm

\no $2.$ Multi-Spaces on graphs\dotfill 3\vskip 3mm

\no $\S 2.1$ \ Graphs\dotfill 3\vskip 2mm

\no$2.1.1$ What is a graph\dotfill 3

\no $2.1.2$ Subgraphs in a graph\dotfill 7

\no $2.1.3$ Classes of graphs with decomposition\dotfill 18

\no $2.1.4$ Operations on graphs\dotfill 26\vskip 2mm

\no $\S 2.2$ \ Multi-Voltage Graphs\dotfill 28\vskip 2mm

\no $2.2.1$ Type $1$\dotfill 28

\no $2.2.2$ Type $2$\dotfill 34\vskip 2mm

\no $\S 2.3$ \ Graphs in a Space\dotfill 40\vskip 2mm

\no $2.3.1$ Graphs in an $n$-manifold\dotfill 40

\no $2.3.2$ Graphs on a surface\dotfill 45

\no $2.3.3$ Multi-Embeddings in an $n$-manifold\dotfill 56

\no $2.3.4$ Classification of graphs in an $n$-manifold\dotfill
59\vskip 2mm

\no $\S 2.4$ \ Multi-Spaces on Graphs\dotfill 63\vskip 2mm

\no $2.4.1$ A graph model for an operation system\dotfill 63

\no $2.4.2$ Multi-Spaces on graphs\dotfill 66

\no $2.4.3$ Cayley graphs of a multi-group\dotfill 68\vskip 2mm

\no $\S 2.5$ \ Graph Phase Spaces\dotfill 70\vskip 2mm

\no $2.5.1$ Graph phase in a multi-space\dotfill 70

\no $2.5.2$ Transformation of a graph phase\dotfill 73\vskip 2mm

\no $\S 2.6$ \ Remarks and Open Problems\dotfill 76\vskip 3mm

\newpage

\no{\large\bf $2.$ Multi-spaces on graphs}\vskip 8mm

\no As a useful tool for dealing with relations of events, graph
theory has rapidly grown in theoretical results as well as its
applications to real-world problems, for example see $[9],[11]$
and $[80]$ for graph theory, $[42]-[44]$ for topological graphs
and combinatorial map theory, $[7],[12]$ and $[104]$ for its
applications to probability, electrical network and real-life
problems. By applying the Smarandache's notion, graphs are models
of multi-spaces and matters in the natural world. For the later,
graphs are a generalization of $p$-branes and seems to be useful
for mechanics and quantum physics.

\vskip 8mm

\no{\bf \S $2.1$ \ Graphs}

\vskip 5mm

\no{\bf $2.1.1.$ What is a graph?}

\vskip 3mm

\no A {\it graph} $G$ is an ordered $3$-tuple $(V,E;I)$, where $V,
E$ are finite sets, $V\not=\emptyset$ and $I: E\rightarrow V\times
V$. Call $V$ the {\it vertex set} and $E$ the {\it edge set} of
$G$, denoted by $V(G)$ and $E(G)$, respectively. Two elements
$v\in V(G)$ and $e\in E(G)$ are said to be {\it incident} if
$I(e)= (v,x)$ or $(x,v)$, where $x\in V(G)$. If $(u,v)=(v,u)$ for
$\forall u,v\in V$, the graph $G$ is called a graph, otherwise, a
directed graph with an orientation $u\rightarrow v$ on each edge
$(u,v)$. Unless Section $2.4$, graphs considered in this chapter
are non-directed.

The cardinal numbers of $|V(G)|$ and $|E(G)|$ are called the {\it
order} and the {\it size} of a graph $G$, denoted by $|G|$ and
$\varepsilon (G)$, respectively.

We can draw a graph $G$ on a plane $\sum$ by representing each
vertex $u$ of $G$ by a point $p(u)$, $p(u)\not=p(v)$ if $u\not= v$
and an edge $(u,v)$ by a plane curve connecting points $p(u)$ and
$p(v)$ on $\sum$, where $p: G\rightarrow P$ is a mapping from the
graph $G$ to $P$.

For example, a graph $G = (V, E; I)$ with $V = \{v_1, v_2, v_3,
v_4 \},$ $E = \{e_1,e_2,e_3,e_4,e_5,$ $e_6,e_7,e_8,e_9,e_{10}\}$
and $I(e_i)=(v_i,v_i), 1\leq i\leq 4; I(e_5)=(v_1,v_2)=(v_2,v_1),
I(e_8)=(v_3,v_4)=(v_4,v_3), I(e_6)=I(e_7)=(v_2,v_3)=(v_3,v_2),
I(e_8)=I(e_9)=(v_4,v_1)=(v_1,v_4)$ can be drawn on a plane as
shown in Fig.$2.1$

\includegraphics[bb=5 10 200 140]{sgm4.eps}\vskip 2mm

\c{\bf Fig $2.1$} \vskip 3mm

In a graph $G = (V,E;I)$, for $\forall e\in E$, if $I(e) = (u,u),
u\in V$, then $e$ is called a {\it loop}. For $\forall e_1, e_2\in
E$, if $I(e_1)=I(e_2)$ and they are not loops, then $e_1$ and
$e_2$ are called {\it multiple edges} of $G$. A graph is {\it
simple} if it is loopless and without multiple edges, i.e.,
$\forall e_1,e_2\in E(\Gamma )$, $I(e_1) \not= I(e_2)$ if
$e_1\not=e_2$ and for $\forall e\in E$, if $I(e) = (u,v)$, then
$u\not=v$. In a simple graph, an edge $(u,v)$ can be abbreviated
to $uv$.

An edge $e\in E(G)$ can be divided into two semi-arcs $e_u, e_v$
if $I(e)=(u,v)$. Call $u$ the {\it root vertex} of the semi-arc
$e_u$. Two semi-arc $e_u,f_v$ are said to be {\it $v$-incident} or
{\it $e-$incident} if $u=v$ or $e=f$. The set of all semi-arcs of
a graph $G$ is denoted by $X_{\frac{1}{2}}(G)$.

A {\it walk} of a graph $\Gamma$ is an alternating sequence of
vertices and edges $u_1, e_1,u_2,e_2,$ $\cdots,e_n, u_{n_1}$ with
$e_i=(u_i,u_{i+1})$ for $1\leq i\leq n$. The number $n$ is the
{\it length of the walk}. If $u_1=u_{n+1}$, the walk is said to be
{\it closed}, and {\it open} otherwise. For example,
$v_1e_1v_1e_5v_2e_6v_3e_3v_3e_7v_2e_2v_2$ is a walk in Fig.$2.1$.
A walk is called a {\it trail} if all its edges are distinct and a
{\it path} if all the vertices are distinct. A closed path is said
to be a {\it circuit}.

A graph $G = (V,E;I)$ is {\it connected} if there is a path
connecting any two vertices in this graph. In a graph, a maximal
connected subgraph is called a {\it component}. A graph $G$ is
{\it $k$-connected} if removing vertices less than $k$ from $G$
still remains a connected graph. Let $G$ be a graph. For $\forall
u\in V(G)$, the neighborhood $N_{G}(u)$ of the vertex $u$ in $G$
is defined by $N_{G}(u)=\{v|\forall(u,v)\in E(G)\}$. The cardinal
number $|N_{G}(u)|$ is called the {\it valency of the vertex $u$}
in the graph $G$ and denoted by $\rho_{G}(u)$. A vertex $v$ with
$\rho_{G}(v)=0$ is called an {\it isolated vertex} and
$\rho_{G}(v)=1$ a {\it pendent vertex}. Now we arrange all
vertices valency of $G$ as a sequence $\rho_G (u)\geq\rho_G
(v)\geq\cdots\geq\rho_G (w)$. Call this sequence the {\it valency
sequence} of $G$. By enumerating edges in $E(G)$, the following
result holds.

$$\sum\limits_{u\in V(G)}\rho_{G}(u)= 2|E(G)|.$$

Give a sequence $\rho_1, \rho_2, \cdots ,\rho_p$ of non-negative
integers. If there exists a graph whose valency sequence is
$\rho_1\geq\rho_2\geq\cdots\geq\rho_p$, then we say that $\rho_1,
\rho_2, \cdots ,\rho_p$ is a {\it graphical sequence}. We have
known the following results (see $[11]$ for details).

\vskip 4mm

\no{\bf Theorem $2.1.1$}(Havel,1955 and Hakimi,1962) \ {\it A
sequence $\rho_1, \rho_2, \cdots ,\rho_p$ of non-negative integers
with $\rho_1\geq\rho_2\geq\cdots\geq\rho_p$, $p\geq 2, \rho_1\geq
1$ is graphical if and only if the sequence $\rho_2-1,
\rho_3-1,\cdots ,\rho_{\rho_1+1}-1, \rho_{\rho_1+2},\cdots
,\rho_p$ is graphical.}

\vskip 4mm

\no{\bf Theorem $2.1.2$}(Erd$\ddot{o}$s and Gallai,1960) \ {\it A
sequence $\rho_1, \rho_2, \cdots ,\rho_p$ of non-negative integers
with $\rho_1\geq\rho_2\geq\cdots\geq\rho_p$ is graphical if and
only if $\sum\limits_{i=1}^p\rho_i$ is even and for each integer
$n, 1\leq n\leq p-1$,}

$$\sum\limits_{i=1}^n\rho_i\leq n(n-1)+\sum\limits_{i=n+1}^pmin\{n,\rho_i\}.$$

A graph $G$ with a vertex set $V(G) = \{v_1,v_2,\cdots , v_p\}$
and an edge set $E(G) = \{e_1,e_2,\cdots , e_q\}$ can be also
described by means of matrix. One such matrix is a $p\times q$
{\it adjacency matrix} $A(G) = [a_{ij}]_{p\times q}$, where
$a_{ij}= |I^{-1}(v_i,v_j)|$. Thus, the adjacency matrix of a graph
$G$ is symmetric and is a $0,1$-matrix having $0$ entries on its
main diagonal if $G$ is simple. For example, the adjacency matrix
$A(G)$ of the graph in Fig.$2.1$ is

\[
A(G)=\left[\begin{array}{cccc}
1 & 1 & 0 & 2\\
1 & 1 & 2 & 0\\
0 & 2 & 1 & 1\\
2 & 0 & 1 & 1
\end{array}
\right]
\]

Let $G_1 = (V_1,E_1;I_1)$ and $G_2 = (V_2,E_2;I_2)$ be two graphs.
They are {\it identical}, denoted by $G_1=G_2$ if $V_1=V_2,
E_1=E_2$ and $I_1=I_2$. If there exists a $1-1$ mapping $\phi:
E_1\rightarrow E_2$ and $\phi: V_1\rightarrow V_2$ such that $\phi
I_1(e) = I_2\phi (e)$ for $\forall e\in E_1$ with the convention
that $\phi (u,v) = (\phi (u),\phi (v))$, then we say that $G_1$ is
{\it isomorphic} to $G_2$, denoted by $G_1\cong G_2$ and $\phi$ an
{\it isomorphism} between $G_1$ and $G_2$. For simple graphs $H_1,
H_2$, this definition can be simplified by $(u,v)\in I_1(E_1)$ if
and only if $(\phi (u),\phi (v))\in I_2(E_2)$ for $\forall u,v\in
V_1$.

For example, let $G_1 = (V_1,E_1;I_1)$ and $G_2 = (V_2,E_2;I_2)$
be two graphs with

$$V_1 = \{v_1,v_2,v_3\},$$

$$E_1 = \{e_1,e_2,e_3,e_4\},$$

$$I_1(e_1)=(v_1,v_2), I_1(e_2)=(v_2,v_3), I_1(e_3)=(v_3,v_1), I_1(e_4)=(v_1,v_1)$$

\no and

$$V_2 = \{u_1,u_2,u_3\},$$

$$E_2 = \{f_1,f_2,f_3,f_4\},$$

$$I_2(f_1)=(u_1,u_2), I_2(f_2)=(u_2,u_3), I_2(f_3)=(u_3,u_1), I_2(f_4)=(u_2,u_2),$$

\no i.e., the graphs shown in Fig.$2.2$.

\includegraphics[bb=5 10 200 150]{sgm5.eps}\vskip 2mm

\c{\bf Fig $2.2$} \vskip 3mm

\no Then they are isomorphic since we can define a mapping $\phi:
E_1\bigcup V_1\rightarrow E_2\bigcup V_2$ by

$$\phi (e_1)=f_2, \phi (e_2)=f_3, \phi (e_3)=f_1, \phi (e_4)=f_4$$

\no and $\phi (v_i)=u_i$ for $1\leq i\leq 3$. It can be verified
immediately that $\phi I_1(e) = I_2\phi (e)$ for $\forall e\in
E_1$. Therefore, $\phi$ is an isomorphism between $G_1$ and $G_2$.

If $G_1=G_2=G$, an isomorphism between $G_1$ and $G_2$ is said to
be an {\it automorphism} of $G$. All automorphisms of a graph $G$
form a group under the composition operation, i.e., $\phi\theta
(x)= \phi(\theta (x))$, where $x\in E(G)\bigcup V(G)$. We denote
the automorphism group of a graph $G$ by ${\rm Aut}G$.

For a simple graph $G$ of $n$ vertices, it is easy to verify that
${\rm Aut}G\leq S_n$, the symmetry group action on these $n$
vertices of $G$. But for non-simple graph, the situation is more
complex. The automorphism groups of graphs $K_m, m = |V(K_m)|$ and
$B_n, n = |E(B_n)|$ in Fig.$2.3$ are ${\rm Aut}K_m = S_m$ and
${\rm Aut}B_n = S_n$.

\includegraphics[bb=5 10 200 140]{sgm6.eps}\vskip 2mm

\c{\bf Fig $2.3$}

\vskip 2mm

For generalizing the conception of automorphisms, the semi-arc
automorphisms of a graph were introduced in $[53]$, which is
defined in the following definition.

\vskip 4mm

\no{\bf Definition $2.1.1$}  {\it A one-to-one mapping $\xi$ on
$X_{\frac{1}{2}}(G)$ is called a semi-arc automorphism of a graph
$G$ if $\xi (e_u)$ and $\xi (f_v)$ are $v-$incident or
$e-$incident if $e_u$ and $f_v$ are $v-$incident or $e-$incident
for $\forall e_u,f_v\in X_{\frac{1}{2}}(G)$.}

\vskip 3mm

All semi-arc automorphisms of a graph also form a group, denoted
by ${\rm Aut}_{\frac{1}{2}}G$. For example, ${\rm
Aut}_{\frac{1}{2}}B_n = S_n[S_2]$.

For $\forall g\in {\rm Aut}G$, there is an induced action
$g|^{\frac{1}{2}}: X_{\frac{1}{2}}(G)\rightarrow
X_{\frac{1}{2}}(G)$ on $X_{\frac{1}{2}}(G)$ defined by

$$\forall e_u\in X_{\frac{1}{2}}(G), g(e_u)=g(e)_{g(u)}. $$

\no{All induced action of elements in ${\rm Aut}G$ is denoted by
${\rm Aut}G|^{\frac{1}{2}}.$}

The graph $B_n$ shows that ${\rm Aut}_{\frac{1}{2}}G$ may be not
the same as ${\rm Aut}G|^{\frac{1}{2}}$. However, we get a result
in the following.

\vskip 4mm \no{\bf Theorem $2.1.3$}([56]) \ {\it For a graph
$\Gamma$ without loops,}

$${\rm Aut}_{\frac{1}{2}}\Gamma = {\rm Aut}\Gamma|^{\frac{1}{2}}.$$

\vskip 3mm

Various applications of this theorem to graphs, especially, to
combinatorial maps can be found in references $[55]-[56]$ and
$[66]-[67]$.

\vskip 5mm

\no{\bf $2.1.2.$ Subgraphs in a graph}

\vskip 4mm

\no A graph $H = (V_1,E_1;I_1)$ is a {\it subgraph} of a graph $G
= (V,E;I)$ if $V_1\subseteq V$, $E_1\subseteq E$ and $I_1:
E_1\rightarrow V_1\times V_1$. We denote that $H$ is a subgraph of
$G$ by $H\subset G$. For example, graphs $G_1,G_2,G_3$ are
subgraphs of the graph $G$ in Fig.$2.4$.

\includegraphics[bb=5 10 200 150]{sgm7.eps}\vskip 3mm

\c{\bf Fig $2.4$}

\vskip 3mm

For a nonempty subset $U$ of the vertex set $V(G)$ of a graph $G$,
the subgraph $\left<U\right>$ of $G$ {\it induced} by $U$ is a
graph having vertex set $U$ and whose edge set consists of these
edges of $G$ incident with elements of $U$. A subgraph $H$ of $G$
is called {\it vertex-induced} if $H\cong\left<U\right>$ for some
subset $U$ of $V(G)$.  Similarly, for a nonempty subset $F$ of
$E(G)$, the subgraph $\left<F\right>$ induced by $F$ in $G$ is a
graph having edge set $F$ and whose vertex set consists of
vertices of $G$ incident with at least one edge of $F$. A subgraph
$H$ of $G$ is {\it edge-induced} if $H\cong\left<F\right>$ for
some subset $F$ of $E(G)$. In Fig.$2.4$, subgraphs $G_1$ and $G_2$
are both vertex-induced subgraphs $\left<\{u_1,u_4\}\right>$,
$\left<\{u_2,u_3\}\right>$ and edge-induced subgraphs
$\left<\{(u_1,u_4)\}\right>$, $\left<\{(u_2,u_3)\}\right>$.

For a subgraph $H$ of $G$, if $|V(H)| = |V(G)|$, then $H$ is
called a {\it spanning subgraph} of $G$. In Fig.$2.4$, the
subgraph $G_3$ is a spanning subgraph of the graph $G$. Spanning
subgraphs are useful for constructing multi-spaces on graphs, see
also Section $2.4$.

A spanning subgraph without circuits is called a {\it spanning
forest}. It is called a {\it spanning tree} if it is connected.
The following characteristic for spanning trees of a connected
graph is well-known.

\vskip 4mm

\no{\bf Theorem $2.1.4$} \ {\it A subgraph $T$ of a connected
graph $G$ is a spanning tree if and only if $T$ is connected and
$E(T)= |V(G)|-1$.}

\vskip 3mm

{\it Proof} \ The necessity is obvious. For its sufficiency, since
$T$ is connected and $E(T)= |V(G)|-1$, there are no circuits in
$T$. Whence, $T$ is a spanning tree. \quad\quad $\natural$

A path is also a tree in which each vertex has valency $2$ unless
the two pendent vertices valency $1$. We denote a path with $n$
vertices by $P_n$ and define the {\it length} of $P_n$ to be
$n-1$. For a connected graph $G$, $x,y\in V(G)$, the distance
$d(x,y)$ of $x$ to $y$ in $G$ is defined by

$$d_G(x,y)= min\{ \ |V(P(x,y))|-1 \ | \ P(x,y) \  is \ a \ path \ connecting \ x \
 and \ y \ \}.$$

For $\forall u\in V(G)$, the {\it eccentricity} $e_G(u)$ of $u$ is
defined by

$$e_G(u) = max\{ \ d_G(u,x) \ | x\in V(G)\}.$$

\no A vertex $u^+$ is called an {\it ultimate vertex} of a vertex
$u$ if $d(u,u^+)=e_G(u)$. Not loss of generality, we arrange these
eccentricities of vertices in $G$ in an order
$e_G(v_1),e_G(v_2),\cdots,$ $e_G(v_n)$ with $e_G(v_1)\leq
e_G(v_2)\leq\cdots\leq e_G(v_n)$, where $\{v_1,v_2,\cdots ,v_n\}$
$= V(G)$. The sequence $\{e_G(v_i)\}_{1\leq i\leq s}$ is called an
{\it eccentricity sequence} of $G$. If $\{e_1,e_2,\cdots ,$
$e_s\}=\{e_G(v_1),e_G(v_2),\cdots ,e_G(v_n)\}$ and $e_1 < e_2<
\cdots < e_s$, the sequence $\{e_i\}_{1\leq i\leq s}$ is called an
{\it eccentricity value sequence} of $G$. For convenience, we
abbreviate an integer sequence $\{r-1+i\}_{1\leq i\leq s+1}$ to
$[r,r+s]$.

The {\it radius} $r(G)$ and the {\it diameter} $D(G)$ of $G$ are
defined by

$$r(G) = min \{  e_G(u)  |  u\in V(G) \} \ \ {\rm and} \ \
D(G)= max \{  e_G(u)  |  u\in V(G)  \},$$

\no respectively. For a given graph $G$, if $r(G) = D(G)$, then
$G$ is called a {\it self-centered graph}, i.e., the eccentricity
value sequence of $G$ is $[r(G),r(G)]$. Some characteristics of
self-centered graphs can be found in $[47],[64]$ and $[108]$.

For $\forall x\in V(G)$, we define a {\it distance decomposition
$\{V_i(x)\}_{1\leq i\leq e_G(x)}$ of $G$ with root} $x$ by

$$G = V_1(x)\bigoplus V_2(x)\bigoplus\cdots\bigoplus
V_{e_G(x)}(x)$$

\no where $V_i(x)=\{ \ u \ | d(x,u)=i,u\in V(G)\}$ for any integer
$i, 0\leq i\leq e_G(x)$. We get a necessary and sufficient
condition for the eccentricity value sequence of a simple graph in
the following.

\vskip 4mm

\no{\bf Theorem $2.1.5$}\ {\it A non-decreasing integer sequence
$\{r_i\}_{1\leq i\leq s}$ is a graphical eccentricity value
sequence if and only if}

($i$) \ $r_1\leq r_s\leq 2r_1$;

($ii$) \ {\it $\triangle (r_{i+1},r_i) = | r_{i+1}-r_i |= 1$ for
any integer $i, 1\leq i\leq s-1$.}

\vskip 3mm

{\it Proof} \ If there is a graph $G$ whose eccentricity value
sequence is $\{r_i\}_{1\leq i\leq s}$, then $r_1\leq r_s$ is
trivial. Now we choose three different vertices $u_1,u_2,u_3$ in
$G$ such that $e_G(u_1)=r_1$ and $d_G(u_2,u_3)=r_s$. By
definition, we know that $d(u_1,u_2)\leq r_1$ and $d(u_1,u_3)\leq
r_1$. According to the triangle inequality for distances, we get
that $r_s= d(u_2,u_3) \leq d_G(u_2,u_1)+ d_G(u_1,u_3) =
d_G(u_1,u_2)+ d_G(u_1,u_3)\leq 2r_1$. So $r_1\leq r_s\leq 2r_1$.

Assume $\{e_i\}_{1\leq i\leq s}$ is the eccentricity value
sequence of a graph $G$. Define $\triangle (i)= e_{i+1}-e_i$,
$1\leq i\leq n-1$. We assert that $0\leq\triangle (i)\leq 1$. If
this assertion is not true, then there must exists a positive
integer $I,1\leq I\leq n-1$ such that $\triangle (I)=
e_{I+1}-e_I\geq 2$. Choose a vertex $x\in V(G)$ such that $e_G(x)
= e_I$ and consider the distance decomposition $\{V_i(x)\}_{0\leq
i\leq e_G(x)}$ of $G$ with root $x$.

Notice that it is obvious that $e_G(x)-1\leq e_G(u_1)\leq
e_G(x)+1$ for any vertex $u_1\in V_1(G)$. Since $\triangle (I)\geq
2$, there does not exist a vertex with the eccentricity
$e_G(x)+1$. Whence, we get $e_G(u_1)\leq e_G(x)$ for $\forall
u_1\in V_1(x)$. If we have proved that $e_G(u_j)\leq e_G(x)$ for
$\forall u_j\in V_j(x)$, $1\leq j<e_G(x)$, we consider these
eccentricity values of vertices in $V_{j+1}(x)$. Let $u_{j+1}\in
V_{j+1}(x)$. According to the definition of $\{V_i(x)\}_{0\leq
i\leq e_G(x)}$, there must exists a vertex $u_j\in V_j(x)$ such
that $(u_j,u_{j+1})\in E(G)$. Now consider the distance
decomposition $\{V_i(u_j)\}_{0\leq j\leq e_G(u)}$ of $G$ with root
$u_j$. Notice that $u_{j+1}\in V_1(u_j)$. Thereby we get that

$$e_G(u_{j+1}) \leq e_G(u_j)+1\leq e_G(x)+1.$$

Because we have assumed that there are no vertices with the
eccentricity $e_G(x)+1$, so $e_G(u_{j+1})\leq e_G(x)$ for any
vertex $u_{j+1}\in V_{j+1}(x)$. Continuing this process, we know
that $e_G(y)\leq e_G(x)=e_I$ for any vertex $y\in V(G)$. But then
there are no vertices with the eccentricity $e_I+1$, which
contradicts the assumption that $\triangle (I)\geq 2$. Therefore
$0\leq \triangle (i)\leq 1$ and $\triangle (r_{i+1},r_i)=1, 1\leq
i\leq s-1$.

For any integer sequence $\{r_i\}_{1\leq i\leq s}$ with conditions
($i$) and ($ii$) hold, it can be simply written as $\{r,r+1,\cdots
,r+s-1\}=[r,r+s-1]$, where $s\leq r$. We construct a graph with
the eccentricity value sequence $[r,r+s-1]$ in the following.

\vskip 3mm

\no{\bf Case $1$} \ \ \ $s = 1$ \vskip 2mm

In this case, $\{r_i\}_{1\leq i\leq s}=[r,r]$. We can choose any
self-centered graph with $r(G)=r$, especially, the circuit
$C_{2r}$ of order $2r$. Then its eccentricity value sequence is
$[r,r]$.

\vskip 3mm

\no{\bf Case $2$} \ \ \  $s \geq 2$

\vskip 2mm

Choose a self-centered graph $H$ with $r(H)=r, x\in V(H)$ and a
path $P_s=u_0u_1\cdots u_{s-1}$. Define a new graph $G
=P_s\bigodot H$ as follows:

$V(G)=V(P_s)\bigcup V(H)\setminus\{u_0\}$,

$E(G)=(E(P_s)\bigcup \{(x,u_1)\}\bigcup E(H)\setminus
\{(u_1,u_0)\}$

\no such as the graph $G$ shown in Fig.$2.5$.

\includegraphics[bb=5 10 200 130]{sgm8.eps}\vskip 3mm

\c{\bf Fig $2.5$}

\vskip 3mm

\no Then we know that $e_G(x)=r$,  $e_G(u_{s-1})=r+s-1$  and
$r\leq e_G(x)\leq r+s-1$ for all other vertices $x\in V(G)$.
Therefore, the eccentricity value sequence of G is $[r,r+s-1]$.
This completes the proof. \quad\quad $\natural$

For a given eccentricity value $l$, the {\it multiplicity set}
$N_G(l)$ is defined by $N_G(l)= \{ \ x \ | \ x\in V(G),e(x)=l \
\}$. Jordan proved that the $\left<N_G(r(G))\right>$ in a tree is
a vertex or two adjacent vertices in 1869([$11$]). For a graph
must not being a tree, we get the following result which
generalizes Jordan's result for trees.

\vskip 4mm

\no{\bf Theorem $2.1.6$} \ {\it Let $\{r_i\}_{1\leq i\leq s}$ be a
graphical eccentricity value sequence. If $|N_G(r_I)|$ $=1$, then
there must be $I=1$, i.e., $|N_G(r_i)|\geq 2$ for any integer $i,
2\leq i\leq s$.}

\vskip 3mm

{\it Proof} \ Let $G$ be a graph with the eccentricity value
sequence $\{r_i\}_{1\leq i\leq s}$ and $N_G(r_I)=\{x_0\},
e_G(x_0)=r_I$. We prove that $e_G(x)>e_G(x_0)$ for any vertex
$x\in V(G)\setminus\{x_0\}$. Consider the distance decomposition
$\{V_i(x_0)\}_{0\leq i\leq e_G(x_0)}$ of $G$ with root $x_0$.
First, we prove that $e_G(v_1)=e_G(x_0)+1$ for any vertex $v_1\in
V_1(x_0)$. Since $e_G(x_0)-1 \leq e_G(v_1)\leq e_G(x_0)+1$ for any
vertex $v_1\in V_1(x_0)$, we only need to prove that $e_G(v_1) >
e_G(x_0)$ for any vertex $v_1\in V_1(x_0)$. In fact, since for any
ultimate vertex $x_0^+$ of $x_0$, we have that $d_G(x_0,
x_0^+)=e_G(x_0)$. So $e_G(x_0^+)\geq e_G(x_0)$. Since
$N_G(e_G(x_0))=\{x_0\}, x_0^+\not\in N_G(e_G(x_0))$. Therefore,
$e_G(x_0^+)
> e_G(x_0)$. Choose $v_1\in V_1(x_0)$. Assume the shortest path from $v_1$ to
$x_0^+$ is $P_1=v_1v_2\cdots v_sx_0^+$ and $x_0\not\in V(P_1)$.
Otherwise, we already have $e_G(v_1)
> e_G(x_0)$. Now consider the distance decomposition
$\{V_i(x_0^+)\}_{0\leq i\leq e_G(x_0^+)}$ of $G$ with root
$x_0^+$. We know that $v_s\in V_1(x_0^+)$. So we get that

$$e_G(x_0^+)-1 \leq e_G(v_s)\leq e_G(x_0^+)+1.$$

\no Thereafter we get that $e_G(v_s)\geq e_G(x_0^+)-1\geq
e_G(x_0)$. Because $N_G(e_G(x_0))=\{x_0\}$, so $v_s\not\in
N_G(e_G(x_0))$. We finally get that $e_G(v_s)
> e_G(x_0)$.

Similarly, choose $v_s,v_{s-1},\cdots ,v_2$ to be root vertices
respectively and consider these distance decompositions of $G$
with roots $v_s,v_{s-1},\cdots ,v_2$, we get that

$e_G(v_s) > e_G(x_0),$

$e_G(v_{s-1}) > e_G(x_0),$

$\cdots\cdots\cdots\cdots\cdots ,$

\no and

$e_G(v_1) > e_G(x_0).$

\no Therefore, $e_G(v_1)=e_G(x_0)+1$ for any vertex $v_1\in
V_1(x_0)$.

Now consider these vertices in $V_2(x_0)$. For $\forall v_2\in
V_2(x_0)$, assume that $v_2$ is adjacent to $u_1, u_1\in
V_1(x_0)$. We know that $e_G(v_2)\geq e_G(u_1)-1\geq e_G(x_0)$.
Since $| N_G(e_G(x_0))|=| N_G(r_I)|=1$, we get $e_G(v_2)\geq
e_G(x_0)+1$.

Now assume that we have proved $e_G(v_k)\geq e_G(x_0)+1$ for any
vertex $v_k\in V_1(x_0)$ $\bigcup V_2(x_0)\bigcup\cdots\bigcup
V_k(x_0)$ for $1\leq k <e_G(x_0)$. Let $v_{k+1}\in V_{k+1}(x_0)$
and assume that $v_{k+1}$ is adjacent to $u_k$ in $V_k(x_0)$. Then
we know that $e_G(v_{k+1})\geq e_G(u_k)-1\geq e_G(x_0)$. Since
$|N_G(e_G(x_0))|=1$, we get that $e_G(v_{k+1})\geq e_G(x_0)+1$.
Therefore, $e_G(x)>e_G(x_0)$ for any vertex $x, x\in
V(G)\setminus\{x_0\}$. That is, if $|N_G(r_I)|=1$, then there must
be $I = 1$. \quad\quad $\natural$\vskip 2mm

Theorem $2.1.6$ is the best possible in some cases of trees. For
example, the eccentricity value sequence of a path $P_{2r+1}$ is
$[r,2r]$ and we have that $| N_G(r)|=1$ and $|N_G(k)|=2$  for
$r+1\leq k\leq 2r$. But for graphs not being trees, we only found
some examples satisfying $|N_G(r_1)|=1$ and $| N_G(r_i)|>2$. A
non-tree graph with the eccentricity value sequence $[2,3]$ and
$|NG(2)|=1$ can be found in Fig.$2$ in the reference $[64]$.

For a given graph $G$ and $V_1, V_2\in V(G)$, define an {\it edge
cut} $E_G(V_1,V_2)$ by

$$E_G(V_1,V_2) = \{ \ (u,v)\in E(G) \ | \ u\in V_1, v\in V_2\}.$$

\no A graph $G$ is {\it hamiltonian} if it has a circuit
containing all vertices of $G$. This circuit is called a {\it
hamiltonian circuit}. A path containing all vertices of a graph
$G$ is called a {\it hamiltonian path}. For hamiltonian circuits,
we have the following characteristic.

\vskip 4mm

\no{\bf Theorem $2.1.7$} \ {\it A circuit $C$ of a graph $G$
without isolated vertices is a hamiltonian circuit if and only if
for any edge cut ${\mathcal C}$, $|E(C)\bigcap E({\mathcal
C})|\equiv 0(mod2)$ and $|E(C)\bigcap E({\mathcal C})|\geq 2$.}

\vskip 3mm

{\it Proof} \ For any circuit $C$ and an edge cut ${\mathcal C}$,
the times crossing ${\mathcal C}$ as we travel along $C$ must be
even. Otherwise, we can not come back to the initial vertex. if
$C$ is a hamiltonian circuit, then $|E(C)\bigcap E({\mathcal
C})|\not = 0$. Whence, $|E(C)\bigcap E({\mathcal C})|\geq 2$ and
$|E(C)\bigcap E({\mathcal C})|\equiv 0(mod2)$ for any edge cut
${\mathcal C}$.

Now if a circuit $C$ satisfies $|E(C)\bigcap E({\mathcal C})|\geq
2$ and $|E(C)\bigcap E({\mathcal C})|\equiv 0(mod2)$ for any edge
cut ${\mathcal C}$, we prove that $C$ is a hamiltonian circuit of
$G$. In fact, if $V(G)\setminus V(C)\not= \emptyset$, choose $x\in
V(G)\setminus V(C)$. Consider an edge cut
$E_G(\{x\},V(G)\setminus\{x\})$. Since $\rho_G(x)\not= 0$, we know
that $|E_G(\{x\},V(G)\setminus\{x\})|\geq 1$. But since
$V(C)\bigcap (V(G)\setminus V(C))=\emptyset$, we know that
$|E_G(\{x\},V(G)\setminus\{x\})\bigcap E(C)| = 0$. Contradicts the
fact that $|E(C)\bigcap E({\mathcal C})|\geq 2$ for any edge cut
${\mathcal C}$. Therefore $V(C)=V(G)$ and $C$ is a hamiltonian
circuit of $G$. \quad\quad $\natural$

Let $G$ be a simple graph. The {\it closure} of $G$, denoted by
$C(G)$, is a graph obtained from $G$ by recursively joining pairs
of non-adjacent vertices whose valency sum is at least $|G|$ until
no such pair remains. In $1976$, Bondy and Chv$\acute{a}$tal
proved a very useful theorem for hamiltonian graphs.

\vskip 4mm

\no{\bf Theorem $2.1.8$}([5][8]) \ {\it A simple graph is
hamiltonian if and only if its closure is hamiltonian.}

\vskip 3mm

This theorem generalizes Dirac's and Ore's theorems simultaneously
stated as follows:\vskip 3mm

{Dirac (1952): \it Every connected simple graph $G$ of order
$n\geq 3$ with the minimum valency$\geq\frac{n}{2}$ is
hamiltonian.}\vskip 2mm

{Ore (1960): \it If $G$ is a simple graph of order $n\geq 3$ such
that $\rho_G(u)+\rho_G(v)\geq n$ for all distinct non-adjacent
vertices $u$ and $v$, then $G$ is hamiltonian.}\vskip 2mm

In $1984$, Fan generalized Dirac's theorem to a localized form
([$41$]). He proved that \vskip 3mm

{\it Let $G $ be a $2$-connected simple graph of order $n$. If
Fan's condition:\vskip 2mm

$max\{\rho_G(u),\rho_G(v)\}\geq\frac{n}{2}$\vskip 2mm

\no holds for $\forall u,v\in V(G)$ provided $d_G(u,v)=2$, then
$G$ is hamiltonian.}

After Fan's paper [$17$], many researches concentrated on
weakening Fan's condition and found new localized conditions for
hamiltonian graphs. For example, those results in references
$[4],[48]-[50],[52],[63]$ and $[65]$ are this type. The next
result on hamiltonian graphs is obtained by Shi in 1992 ([$84$]).

\vskip 4mm

\no{\bf Theorem $2.1.9$}(Shi, 1992) \ {\it Let $G$ be a
$2$-connected simple graph of order $n$. Then $G$ contains a
circuit passing through all vertices of valency$\geq\frac{n}{2}$.}

\vskip 3mm

{\it Proof} \ Assume the assertion is false. Let $C = v_1v_2\cdots
v_kv_1$ be a circuit containing as many vertices of
valency$\geq\frac{n}{2}$ as possible and with an orientation on
it. For $\forall v\in V(C)$, $v^+$ denotes the successor and $v^-$
the predecessor of $v$ on $C$. Set $R = V(G)\setminus V(C)$. Since
$G$ is $2$-connected, there exists a path length than $2$
connecting two vertices of $C$ that is internally disjoint from
$C$ and containing one internal vertex $x$ of
valency$\geq\frac{n}{2}$ at least. Assume $C$ and $P$ are chosen
in such a way that the length of $P$ as small as possible. Let
$N_R(x)=N_G(x)\bigcap R$, $N_C(x)=N_G(x)\bigcap C$, $N_C^+(x) = \{
v | v^-\in N_C(x)\}$ and $N_C^-(x) = \{ v | v^+\in N_C(x)\}$.

Not loss of generality, we may assume $v_1\in V(P)\bigcap V(C)$.
Let $v_t$ be the other vertex in $V(P)\bigcap V(C)$. By the way
$C$ was chosen, there exists a vertex $v_s$ with $1 < s < t$ such
that $\rho_G(v_s)\geq\frac{n}{2}$ and $\rho (v_i) < \frac{n}{2}$
for $1 < i < s$.

If $s\geq 3$, by the choice of $C$ and $P$ the sets

$$N_C^-(v_s)\setminus \{v_1\}, N_C(x), N_R(v_s), N_R(x), \{x, v_{s-1}\}$$

\no are pairwise disjoint, implying that

\begin{eqnarray*}
n &\geq& |N_C^-(v_s)\setminus
\{v_1\}|+|N_C(x)|+|N_R(v_s)|+|N_R(x)|+|\{x, v_{s-1}\}|\\
&=& \rho_G(v_s)+\rho_G(x)+1\geq n+1,
\end{eqnarray*}

\no a contradiction. If $s=2$, then the sets

$$N_C^-(v_s), N_C(x), N_R(v_s), N_R(x), \{x\}$$

\no are pairwise disjoint, which yields a similar contradiction.
\quad\quad $\natural$

Three induced subgraphs used in the next result for hamiltonian
graphs are shown in Fig.$2.6$.

\includegraphics[bb=10 10 200 140]{sgm9.eps}\vskip 3mm

\c{\bf Fig $2.6$}

\vskip 3mm

\no For an induced subgraph $L$ of a simple graph $G$, a condition
is called a {\it localized condition} $D_L(l)$ if $D_L(x,y)=l$
implies that $max\{ \rho_G(x), \rho_G(y) \}\geq\frac{|G|}{2}$ for
$\forall x,y\in V(L)$. Then we get the following result.

\vskip 4mm

\no{\bf Theorem $2.1.10$} \ {\it Let $G$ be a $2$-connected simple
graph. If the localized condition $D_L(2)$ holds for induced
subgraphs $L\cong K_{1.3}$ or $Z_2$ in $G$, then $G$ is
hamiltonian.}

\vskip 3mm

{\it Proof} \ By Theorem $2.1.9$, we denote by
$c_{\frac{n}{2}}(G)$ the maximum length of circuits passing
through all vertices$\geq\frac{n}{2}$. Similar to Theorem $2.1.8$,
we know that for $x,y\in V(G)$, if $\rho_G(x)\geq\frac{n}{2}$,
$\rho_G(y)\geq\frac{n}{2}$ and $xy\not\in E(G)$, then
$c_{\frac{n}{2}}(G\bigcup\{xy\})= c_{\frac{n}{2}}(G)$. Otherwise,
if $c_{\frac{n}{2}}(G\bigcup\{xy\}) > c_{\frac{n}{2}}(G)$, there
exists a circuit of length $c_{\frac{n}{2}}(G\bigcup\{xy\})$ and
passing through all vertices$\geq\frac{n}{2}$. Let
$C_{\frac{n}{2}}$ be such a circuit and
$C_{\frac{n}{2}}=xx_1x_2\cdots x_syx$ with
$s=c_{\frac{n}{2}}(G\bigcup\{xy\})-2$. Notice that

$$N_G(x)\bigcap (V(G)\setminus
V(C_{\frac{n}{2}}(G\bigcup\{xy\})))=\emptyset$$

\no and

$$N_G(y)\bigcap (V(G)\setminus
V(C_{\frac{n}{2}}(G\bigcup\{xy\})))=\emptyset.$$

\no If there exists an integer $i, 1\leq i\leq s$, $xx_i\in E(G)$,
then $x_{i-1}y\not\in E(G)$. Otherwise, there is a circuit
$C'=xx_ix_{i+1}\cdots x_syx_{i-1}x_{i-2}\cdots x$ in $G$ passing
through all vertices$\geq\frac{n}{2}$ with length
$c_{\frac{n}{2}}(G\bigcup\{xy\})$. Contradicts the assumption that
$c_{\frac{n}{2}}(G\bigcup\{xy\}) > c_{\frac{n}{2}}(G)$. Whence,

$$\rho_G(x)+\rho_G(y)\leq |V(G)\setminus V(C(C_{\frac{n}{2}}))|+
|V(C(C_{\frac{n}{2}}))|-1 = n-1,$$

\no also contradicts that $\rho_G(x)\geq\frac{n}{2}$ and
$\rho_G(y)\geq\frac{n}{2}$. Therefore,
$c_{\frac{n}{2}}(G\bigcup\{xy\})= c_{\frac{n}{2}}(G)$ and
generally, $c_{\frac{n}{2}}(C(G))= c_{\frac{n}{2}}(G)$.

Now let $C$ be a maximal circuit passing through all
vertices$\geq\frac{n}{2}$ in the closure $C(G)$ of $G$ with an
orientation $\overrightarrow{C}$. According to Theorem $2.1.8$, if
$C(G)$ is non-hamiltonian, we can choose $H$ be a component in
$C(G)\setminus C$. Define $N_C(H)=(\bigcup\limits_{x\in
H}N_{C(G)}(x))\bigcap V(C)$. Since $C(G)$ is $2$-connected, we get
that $|N_C(H)|\geq 2$. This enables us choose vertices $x_1,x_2\in
N_C(H)$, $x_1\not=x_2$ and $x_1$ can arrive at $x_2$ along
$\overrightarrow{C}$. Denote by $x_1\overrightarrow{C} x_2$ the
path from $x_1$ to $x_2$ on $\overrightarrow{C}$ and
$x_2\overleftarrow{C}x_1$ the reverse. Let $P$ be a shortest path
connecting $x_1, x_2$ in $C(G)$ and

$$u_1\in N_{C(G)}(x_1)\bigcap V(H)\bigcap V(P), \ \ \
u_2\in N_{C(G)}(x_2)\bigcap V(H)\bigcap V(P).$$

\no Then

$$E(C(G))\bigcap (\{x_1^-x_2^-,x_1^+x_2^+\}\bigcup
E_{C(G)}(\{u_1,u_2\},\{x_1^-,x_1^+,x_2^-,x_2^+\}))=\emptyset$$

\no and

$$\left<\{x_1^-,x_1,x_1^+,u_1\}\right>\not\cong K_{1.3} \ {\rm or} \
\left<\{x_2^-,x_2,x_2^+,u_2\}\right>\not\cong K_{1.3}.$$

\no Otherwise, there exists a circuit longer than $C$, a
contradiction. To prove this theorem, we consider two cases.

\vskip 3mm

\no{\bf Case $1$} \ $\left<\{x_1^-,x_1,x_1^+,u_1\}\right>\not\cong
K_{1.3}$ and $\left<\{x_2^-,x_2,x_2^+,u_2\}\right>\not\cong
K_{1.3}$

\vskip 2mm

In this case, $x_1^-x_1^+\in E(C(G))$ and $x_2^-x_2^+\in E(C(G))$.
By the maximality of $C$ in $C(G)$, we have two claims.

\vskip 3mm

\no{\bf Claim $1.1$} \ \ $u_1=u_2=u$

\vskip 2mm

Otherwise, let $P=x_1u_1y_1\cdots y_lu_2$. By the choice of $P$,
there must be

$$\left<\{x_1^-,x_1,x_1^+,u_1,y_1\}\right>\cong Z_2 \ {\rm and} \
\left<\{x_2^-,x_2,x_2^+,u_2,y_l\}\right>\cong Z_2$$

Since $C(G)$ also has the $D_L(2)$ property, we get that

$$max\{\rho_{C(G)}(x_1^-),\rho_{C(G)}(u_1)\}\geq\frac{n}{2}, \ \
max\{\rho_{C(G)}(x_12^-),\rho_{C(G)}(u_2)\}\geq\frac{n}{2}.$$

\no Whence, $\rho_{C(G)}(x_1^-)\geq\frac{n}{2}$,
$\rho_{C(G)}(x_2^-)\geq\frac{n}{2}$ and $x_1^-x_2^-\in E(C(G))$, a
contradiction.

\vskip 3mm

\no{\bf Claim $1.2$} \ \ $x_1x_2\in E(C(G))$

\vskip 2mm

If $x_1x_2\not\in E(C(G))$, then
$\left<\{x_1^-,x_1,x_1^+,u,x_2\}\right>\cong Z_2$. Otherwise, we
get $x_2x_1^-\in E(C(G))$ or $x_2x_1^+\in E(C(G))$. But then there
is a circuit

$$C_1=x_2^+\overrightarrow{C}x_1^-x_2ux_1\overrightarrow{C}x_2^-x_2^+
\ {\rm or} \
C_2=x_2^+\overrightarrow{C}x_1ux_2x_1^+\overrightarrow{C}x_2^-x_2^+.
$$

\no Contradicts the maximality of $C$. Therefore, we know that

$$\left<\{x_1^-,x_1,x_1^+,u,x_2\}\right>\cong Z_2.$$

\no By the property $D_L(2)$, we get that
$\rho_{C(G)}(x_1^-)\geq\frac{n}{2}$

Similarly, consider the induced subgraph
$\left<\{x_2^-,x_2,x_2^+,u,x_2\}\right>$, we get that
$\rho_{C(G)}(x_2^-)$ $\geq\frac{n}{2}$. Whence, $x_1^-x_2^-\in
E(C(G))$, also a contradiction. Thereby we know the structure of
$G$ as shown in Fig.$2.7$.

\includegraphics[bb=0 10 200 150]{sgm10.eps}\vskip 3mm

\c{\bf Fig $2.7$}

By the maximality of $C$ in $C(G)$, it is obvious that
$x_1^{--}\not=x_2^+$. We construct an induced subgraph sequence
$\{G_i\}_{1\leq i\leq m}$, $m=|V(x_1^-\overleftarrow{C}x_2^+)|-2$
and prove there exists an integer $r, 1\leq r\leq m$ such that
$G_r\cong Z_2$.

First, we consider the induced subgraph $G_1 =
\left<\{x_1,u,x_2,x_1^-,x_1^{--}\}\right>$. If $G_1\cong Z_2$,
take $r=1$. Otherwise, there must be

$$\{x_1^-x_2, x_1^{--}x_2, x_1^{--}u, x_1^{--}x_1\}\bigcap E(C(G))\not=\emptyset.$$

If $x_1^-x_2\in E(C(G))$, or $x_1^{--}x_2\in E(C(G))$, or
$x_1^{--}u\in E(C(G))$, there is a circuit
$C_3=x_1^-\overleftarrow{C}x_2^+x_2^-\overleftarrow{C}x_1ux_2x_1^-$,
or
$C_4=x_1^{--}\overleftarrow{C}x_2^+x_2^-\overleftarrow{C}x_1^+x_1^-x_1ux_2x_1^{--}$,
or $C_5 = x_1^{--}\overleftarrow{C}x_1^+x_1^-x_1ux_1^{--}$. Each
of these circuits contradicts the maximality of $C$. Therefore,
$x_1^{--}x_1\in E(C(G))$.

Now let $x_1^-\overleftarrow{C}x_2^+ = x_1^-y_1y_2\cdots
y_mx_2^+$, where $y_0=x_1^-, y_1=x_1^{--}$ and $y_m=x_2^{++}$. If
we have defined an induced subgraph $G_k$ for any integer $k$ and
have gotten $y_ix_1\in E(C(G))$ for any integer $i, 1\leq i\leq k$
and $y_{k+1}\not=x_2^{++}$, then we define

$$G_{k+1} = \left<\{y_{k+1},y_k,x_1,x_2,u\}\right>.$$

\no If $G_{k+1}\cong Z_2$, then $r=k+1$. Otherwise, there must be

$$\{y_ku,y_kx_2,y_{k+1}u,y_{k+1}x_2,y_{k+1}x_1\}\bigcap E(C(G))\not=\emptyset.$$

If $y_ku\in E(C(G))$, or $y_kx_2\in E(C(G))$, or $y_{k+1}u\in
E(C(G))$, or $y_{k+1}x_2\in E(C(G))$, there is a circuit $C_6 =
y_k\overleftarrow{C}x_1^+x_1^-\overleftarrow{C}y_{k-1}x_1uy_k$, or
$C_7 =
y_k\overleftarrow{C}x_2^+x_2^-\overleftarrow{C}x_1^+x_1^-\overleftarrow{C}$
$y_{k-1}x_1ux_2y_k$, or $C_8 =
y_{k+1}\overleftarrow{C}x_1^+x_1^-\overleftarrow{C}y_kx_1uy_{k+1}$,
or $C_9 =
y_{k+1}\overleftarrow{C}x_2^+x_2^-\overleftarrow{C}x_1^+x_1^-\overleftarrow{C}y_k
x_1u$ $x_2y_{k+1}$. Each of these circuits contradicts the
maximality of $C$. Thereby, $y_{k+1}x_1\in E(C(G))$.

Continue this process. If there are no subgraphs in
$\{G_i\}_{1\leq i\leq m}$ isomorphic to $Z_2$, we finally get
$x_1x_2^{++}\in E(C(G))$. But then there is a circuit $C_{10} =
x_1^-\overleftarrow{C}x_2^{++}x_1ux_2x_2^+$
$\overleftarrow{C}x_1^+x_1^-$ in $C(G)$. Also contradicts the
maximality of $C$ in $C(G)$. Therefore, there must be an integer
$r, 1\leq r\leq m$ such that $G_r\cong Z_2$.

Similarly, let $x_2^-\overleftarrow{C}x_1^+ = x_2^-z_1z_2\cdots
z_tx_1^-$, where $t=|V(x_2^-\overleftarrow{C}x_1^+)|-2$,
$z_0=x_2^-, z_1^{++}=x_2, z_t=x_1^{++}$. We can also construct an
induced subgraph sequence $\{G^i\}_{1\leq i\leq t}$ and know that
there exists an integer $h, 1\leq h\leq t$ such that $G^h\cong
Z_2$ and $x_2z_i\in E(C(G))$ for $0\leq i\leq h-1$.

Since the localized condition $D_L(2)$ holds for an induced
subgraph $Z_2$ in $C(G)$, we get that
$max\{\rho_{C(G)}(u),\rho_{C(G)}(y_{r-1})\}\geq\frac{n}{2}$ and
$max\{\rho_{C(G)}(u),\rho_{C(G)}(z_{h-1})\}\geq\frac{n}{2}$.
Whence $\rho_{C(G)}(y_{r-1})\geq\frac{n}{2}$,
$\rho_{C(G)}(z_{h-1})\geq\frac{n}{2}$ and  $y_{r-1}z_{h-1}\in
E(C(G))$. But then there is a circuit

$$C_{11}=y_{r-1}\overleftarrow{C}x_2^+x_2^-\overleftarrow{C}z_{h-2}x_2ux_1y_{r-2}
\overrightarrow{C}x_1^-x_1^+\overrightarrow{C}z_{h-1}y_{r-1}$$

\no in $C(G)$, where if $h=1$, or $r=1$,
$x_2^-\overleftarrow{C}z_{h-2}=\emptyset$, or
$y_{r-2}\overrightarrow{C}x_1^-=\emptyset$. Also contradicts the
maximality of $C$ in $C(G)$.

\vskip 3mm

\no{\bf Case $2$} \ $\left<\{x_1^-,x_1,x_1^+,u_1\}\right>\not\cong
K_{1.3}$, $\left<\{x_2^-,x_2,x_2^+,u_2\}\right>\cong K_{1.3}$ or
$\left<\{x_1^-,x_1,x_1^+,u_1\}\right>\cong K_{1.3}$,
$\left<\{x_2^-,x_2,x_2^+,u_2\}\right>\not\cong K_{1.3}$

\vskip 2mm

Not loss of generality, we assume that
$\left<\{x_1^-,x_1,x_1^+,u_1\}\right>\not\cong K_{1.3}$,
$\left<\{x_2^-,x_2,x_2^+,u_2\}\right>$ $\cong K_{1.3}$. Since each
induced subgraph $K_{1.3}$ in $C(G)$ possesses $D_L(2)$, we get
that $max\{\rho_{C(G)}(u),$ $\rho_{C(G)}(x_2^-)\}\geq\frac{n}{2}$
and $max\{\rho_{C(G)}(u),\rho_{C(G)}(x_2^+)\}\geq\frac{n}{2}$.
Whence $\rho_{C(G)}(x_2^-)$ $\geq\frac{n}{2}$,
$\rho_{C(G)}(x_2^+)\geq\frac{n}{2}$ and $x_2^-x_2^+\in E(C(G))$.
Therefore, the discussion of Case $1$ also holds in this case and
yields similar contradictions.

Combining Case $1$ with Case $2$, the proof is complete.\quad\quad
$\natural$

Let $G$, $F_1,F_2,\cdots , F_k$ be $k+1$ graphs. If there are no
induced subgraphs of $G$ isomorphic to $F_i, 1\leq i\leq k$, then
$G$ is called {\it $\{F_1,F_2,\cdots ,F_k\}$-free}. we get a
immediately consequence by Theorem $2.1.10$.

\vskip 4mm

\no{\bf Corollary $2.1.1$} \ {\it Every $2$-connected $\{K_{1.3},
Z_2\}$-free graph is hamiltonian.}

\vskip 3mm

Let $G$ be a graph. For $\forall u\in V(G)$, $\rho_G(u)= d$, let
$H$ be a graph with $d$ pendent vertices $v_1,v_2,\cdots ,v_d$.
Define a splitting operator $\vartheta : G\rightarrow G^{\vartheta
(u)}$ on $u$ by

$$V(G^{\vartheta (u)})= (V(G)\setminus\{u\})\bigcup (V(H)\setminus\{v_1,v_2,\cdots
,v_d\}),$$

\begin{eqnarray*}
E(G^{\vartheta (u)})&=& (E(G)\setminus\{ux_i\in E(G), 1\leq i\leq
d\})\\
&\bigcup& (E(H)\setminus\{v_iy_i\in E(H), 1\leq i\leq
d\})\bigcup\{x_iy_i, 1\leq i\leq d\}.
\end{eqnarray*}

\no We call $d$ the {\it degree of the splitting operator
$\vartheta$} and $N(\vartheta(u))= H\setminus\{x_iy_i, 1\leq i\leq
d\}$ the {\it nucleus of $\vartheta$ }. A splitting operator is
shown in Fig.$2.8$.

\includegraphics[bb=0 10 200 130]{sgm11.eps}\vskip 3mm

\c{\bf Fig $2.9$}

Erd$\ddot{o}$s and R$\acute{e}$nyi raised a question in 1961 (
[$7$]): {\it in what model of random graphs is it true that almost
every graph is hamiltonian?} P$\acute{o}$sa and Korshuuov proved
independently that for some constant $c$ almost every labelled
graph with $n$ vertices and at least $n\log n$ edges is
hamiltonian in 1974. Contrasting this probabilistic result, there
is another property for hamiltonian graphs, i.e., there is a
splitting operator $\vartheta$ such that $G^{\vartheta (u)}$ is
non-hamiltonian for $\forall u\in V(G)$ of a graph $G$.

\vskip 4mm

\no{\bf Theorem $2.1.11$} \ {\it Let $G$ be a graph. For $\forall
u\in V(G)$, $\rho_G(u)=d$, there exists a splitting operator
$\vartheta$ of degree $d$ on $u$ such that $G^{\vartheta (u)}$ is
non-hamiltonian.}

\vskip 3mm

{\it Proof} \ For any positive integer $i$, define a simple graph
$\Theta_i$ by $V(\Theta_i )=\{x_i,y_i,z_i,$ $u_i\}$ and
$E(\Theta_i )=\{x_iy_i,x_iz_i,y_iz_i,y_iu_i,z_iu_i\}$. For
integers $\forall i,j\geq 1$, the point product
$\Theta_i\odot\Theta_j$ of $\Theta_i$ and $\Theta_j$ is defined by

$$V(\Theta_i\odot\Theta_j)=V(\Theta_i)\bigcup V(\Theta_j)\setminus\{u_j\},$$

$$E(\Theta_i\odot\Theta_j)=E(\Theta_i)\bigcup E(\Theta_j)
\bigcup\{x_iy_j,x_iz_j\}\setminus\{x_jy_j,x_jz_j\}.$$

Now let $H_d$ be a simple graph with

$$V(H_d)= V(\Theta_1\odot\Theta_2\odot\cdots\Theta_{d+1})\bigcup\{v_1,v_2,\cdots ,v_d\},$$

$$E(H_d)= E(\Theta_1\odot\Theta_2\odot\cdots\Theta_{d+1})\bigcup
\{v_1u_1,v_2u_2,\cdots ,v_du_d\}.$$

\no Then $\vartheta :G\rightarrow G^{\vartheta (w)}$ is a
splitting operator of degree $d$ as shown in Fig.$2.10$.

\includegraphics[bb=0 10 200 130]{sgm12.eps}\vskip 3mm

\c{\bf Fig $2.10$} \vskip 2mm

For any graph $G$ and $w\in V(G), \rho_G(w)=d$, we prove that
$G^{\vartheta (w)}$ is non-hamiltonian. In fact, If $G^{\vartheta
(w)}$ is a hamiltonian graph, then there must be a hamiltonian
path $P(u_i,u_j)$ connecting two vertices $u_i,u_j$ for some
integers $i,j, 1\leq i,j\leq d$ in the graph
$H_d\setminus\{v_1,v_2,\cdots , v_d\}$. However, there are no
hamiltonian path connecting vertices $u_i,u_j$ in the graph
$H_d\setminus\{v_1,v_2,\cdots , v_d\}$ for any integer $i,j, 1\leq
i,j\leq d$. Therefore, $G^{\vartheta (w)}$ is
non-hamiltonian.\quad\quad $\natural$

\vskip 5mm

\no{\bf $2.1.3.$ Classes of graphs with decomposition}

\vskip 4mm

\no{\bf $(1)$ Typical classes of graphs}

\vskip 3mm

\no{\bf C1. Bouquets and Dipoles}

\vskip 3mm

\no In graphs, two simple cases is these graphs with one or two
vertices, which are just bouquets or dipoles. A graph $B_n =
(V_b,E_b;I_b)$ with $V_b =\{ \ O \ \}$, $E_b =\{e_1,e_2,\cdots
,e_n\}$ and $I_b(e_i) = (O,O)$ for any integer $i, 1\leq i\leq n$
is called a {\it bouquet} of $n$ edges. Similarly, a graph
$D_{s.l.t} = (V_d,E_d;I_d)$ is called a {\it dipole} if $V_d =
\{O_1,O_2\}$, $E_d = \{e_1,e_2,\cdots ,e_s,e_{s+1},\cdots ,
e_{s+l}, e_{s+l+1},\cdots , e_{s+l+t}\}$ and

\vskip 3mm

\[
I_d(e_i)=\left\{\begin{array}{ll}
(O_1,O_1), & \ {\rm if} \ 1\leq i\leq s,\\
(O_1,O_2), & \ {\rm if} \ s+1\leq i\leq s+l,\\
(O_2,O_2), & \ {\rm if} \ s+l+1\leq i\leq s+l+t.\\
\end{array}
\right.
\]

\vskip 2mm

For example, $B_3$ and $D_{2,3,2}$ are shown in Fig.$2.11$.

\includegraphics[bb=0 10 200 110]{sgm13.eps}\vskip 3mm

\c{\bf Fig $2.11$} \vskip 2mm

In the past two decades, the behavior of bouquets on surfaces
fascinated many mathematicians. A typical example for its
application to mathematics is the classification theorem of
surfaces. By a combinatorial view, these connected sums of tori,
or these connected sums of projective planes used in this theorem
are just bouquets on surfaces. In Section $2.4$, we will use them
to construct completed multi-spaces.

\vskip 3mm

\no{\bf C2. Complete graphs}

\vskip 3mm

\no A {\it complete graph} $K_n = (V_c,E_c;I_c)$ is a simple graph
with $V_c = \{v_1,v_2,\cdots ,v_n\}$, $E_c = \{e_{ij},1\leq
i,j\leq n, i\not=j\}$ and $I_c(e_{ij})= (v_i,v_j)$. Since $K_n$ is
simple, it can be also defined by a pair $(V,E)$ with $V =
\{v_1,v_2,\cdots ,v_n\}$ and $E = \{v_iv_j, 1\leq i,j\leq n,
i\not=j\}$. The one edge graph $K_2$ and the triangle graph $K_3$
are both complete graphs.

A complete subgraph in a graph is called a {\it clique}.
Obviously, every graph is a union of its cliques.

\vskip 3mm

\no{\bf C3. $r$-Partite graphs}

\vskip 3mm

\no A simple graph $G = (V, E; I)$ is {\it r-partite} for an
integer $r\geq 1$ if it is possible to partition $V$ into $r$
subsets $V_1,V_2,\cdots , V_r$ such that for $\forall e\in E$,
$I(e)=(v_i,v_j)$ for $v_i\in V_i$, $v_j\in V_j$ and $i\not=j,
1\leq i,j\leq r$. Notice that by definition, there are no edges
between vertices of $V_i$, $1\leq i\leq r$. A vertex subset of
this kind in a graph is called an {\it independent vertex subset}.

For $n=2$, a $2$-partite graph is also called a {\it bipartite
graph}. It can be shown that {\it a graph is bipartite if and only
if there are no odd circuits in this graph}. As a consequence, a
tree or a forest is a bipartite graph since they are circuit-free.

Let $G = (V, E; I)$ be an r-partite graph and let $V_1,V_2,\cdots
, V_r$ be its $r$-partite vertex subsets. If there is an edge
$e_{ij}\in E$ for $\forall v_i\in V_i$ and $\forall v_j\in V_j$,
where $1\leq i,j\leq r, i\not=j$ such that $I(e)=(v_i,v_j)$, then
we call $G$ a {\it complete $r$-partite graph}, denoted by $G=
K(|V_1|,|V_2|,\cdots , |V_r|)$. Whence, a complete graph is just a
complete $1$-partite graph. For an integer $n$, the complete
bipartite graph $K(n,1)$ is called a {\it star}. For a graph $G$,
we have an obvious formula shown in the following, which
corresponds to the neighborhood decomposition in topology.

$$E(G) = \bigcup\limits_{x\in V(G)}E_G(x, N_G(x)).$$

\vskip 3mm

\no{\bf C4. Regular graphs}

\vskip 3mm

\no A graph $G$ is {\it regular of valency $k$} if $\rho_G(u)=k$
for $\forall u\in V(G)$. These graphs are also called {\it
$k$-regular}. There $3$-regular graphs are referred to as {\it
cubic graphs}. A $k$-regular vertex-spanning subgraph of a graph
$G$ is also called a $k$-factor of $G$.

For a $k$-regular graph $G$, since $k|V(G)|=2|E(G)|$, thereby one
of $k$ and $|V(G)|$ must be an even number, i.e., there are no
$k$-regular graphs of odd order with $k\equiv 1(mod 2)$. A
complete graph $K_n$ is ($n-1$)-regular and a complete $s$-partite
graph $K(p_1,p_2,\cdots ,p_s)$ of order $n$ with $p_1=p_2=\cdots
=p_s=p$ is $(n-p)$-regular.

In regular graphs, those of simple graphs with high symmetry are
particularly important to mathematics. They are related
combinatorics with group theory and crystal geometry. We briefly
introduce them in the following.

Let $G$ be a simple graph and $H$ a subgroup of ${\rm Aut}G$. $G$
is said to be {\it $H$-vertex transitive}, {\it $H$-edge
transitive} or {\it $H$-symmetric} if $H$ acts transitively on the
vertex set $V(G)$, the edge set $E(G)$ or the set of ordered
adjacent pairs of vertex of $G$. If $H = {\rm Aut}G$, an
$H$-vertex transitive, an $H$-edge transitive or an $H$-symmetric
graph is abbreviated to a {\it vertex-transitive}, an {\it
edge-transitive} or a {\it symmetric} graph.

Now let $\Gamma$ be a finite generated group and
$S\subseteq\Gamma$ such that $1_{\Gamma}\not\in S$ and $S^{-1}=\{
x^{-1} | x\in S \}=S$. A {\it Cayley graph} $Cay(\Gamma :S)$ is a
simple graph with vertex set $V(G) = \Gamma$ and edge set $E(G)=\{
(g,h) | g^{-1}h\in S \}$. By the definition of Cayley graphs, we
know that {\it a Cayley graph Cay$(\Gamma :S)$ is complete if and
only if $S=\Gamma\setminus\{1_{\Gamma}\}$ and connected if and
only if $\Gamma = \left<S\right>$.}

\vskip 4mm

\no{\bf Theorem $2.1.12$} \ {\it A Cayley graph Cay$(\Gamma :S)$
is vertex-transitive.}

\vskip 3mm

{\it Proof} \ For $\forall g\in\Gamma$, define a permutation
$\zeta_g$ on $V({\rm Cay}(\Gamma :S))=\Gamma$ by $\zeta_g(h)=gh,
h\in\Gamma$. Then $\zeta_g$ is an automorphism of Cay$(\Gamma :S)$
for $(h,k)\in E({\rm Cay}(\Gamma :S))\Rightarrow h^{-1}k\in
S\Rightarrow (gh)^{-1}(gk)\in S\Rightarrow
(\zeta_g(h),\zeta_g(k))\in E({\rm Cay}(\Gamma :S))$.

Now we know that $\zeta_{kh^{-1}}(h)=(kh^{-1})h=k$ for $\forall
h,k\in\Gamma$. Whence, Cay$(\Gamma :S)$ is vertex-transitive.
\quad\quad $\natural$

Not every vertex-transitive graph is a Cayley graph of a finite
group. For example, the Petersen graph is vertex-transitive but
not a Cayley graph(see $[10],[21]]$ and $[110]$ for details).
However, every vertex-transitive graph can be constructed almost
like a Cayley graph. This result is due to Sabidussi in 1964. The
readers can see $[110]$ for a complete proof of this result.

\vskip 4mm

\no{\bf Theorem $2.1.13$} \ {\it Let $G$ be a vertex-transitive
graph whose automorphism group is $A$. Let $H=A_b$ be the
stabilizer of $b\in V(G)$. Then $G$ is isomorphic with the
group-coset graph $C(A/H,S)$, where $S$ is the set of all
automorphisms $x$ of $G$ such that $(b,x(b))\in E(G)$,
$V(C(A/H,S))=A/H$ and $E(C(A/H,S))=\{ (xH,yH) | x^{-1}y\in HSH
\}$.}

\vskip 5mm

\no{\bf C5. Planar graphs}

\vskip 4mm

\no Every graph is drawn on the plane. A graph is {\it planar} if
it can be drawn on the plane in such a way that edges are disjoint
expect possibly for endpoints. When we remove vertices and edges
of a planar graph $G$ from the plane, each remained connected
region is called a {\it face} of $G$. The length of the boundary
of a face is called its {\it valency}. Two planar graphs are shown
in Fig.$2.12$.

\includegraphics[bb=0 10 200 105]{sgm14.eps}\vskip 2mm

\c{\bf Fig $2.12$} \vskip 2mm

For a planar graph $G$, its order, size and number of faces are
related by a well-known formula discovered by Euler.

\vskip 4mm

\no{\bf Theorem $2.1.14$} \ {\it let $G$ be a planar graph with
$\phi (G)$ faces. Then}

$$|G|-\varepsilon (G)+\phi (G) = 2.$$

\vskip 3mm

{\it Proof} \ We can prove this result by employing induction on
$\varepsilon (G)$. See $[42]$ or $[23],[69]$ for a complete
proof.\quad\quad $\natural$

For an integer $s, s\geq 3$, an $s$-regular planar graph with the
same length $r$ for all faces is often called an {\it
$(s,r)$-polyhedron}, which are completely classified by the
ancient Greeks.

\vskip 4mm

\no{\bf Theorem $2.1.15$} \ {There are exactly five polyhedrons,
two of them are shown in Fig.$2.12$, the others are shown in
Fig.$2.13$.}

\includegraphics[bb=0 10 200 140]{sgm15.eps}\vskip 3mm

\c{\bf Fig $2.13$} \vskip 2mm

\vskip 3mm

{\it Proof} \ Let $G$ be a $k$-regular planar graph with $l$
faces. By definition, we know that $|G|k = \phi (G)l =
2\varepsilon (G)$. Whence, we get that $|G|=\frac{2\varepsilon
(G)}{k}$ and $\phi (G)=\frac{2\varepsilon (G)}{l}$. According to
Theorem $2.1.14$, we get that

$$\frac{2\varepsilon
(G)}{k}-\varepsilon (G)+\frac{2\varepsilon (G)}{l} = 2.$$

\no i.e.,

$$\varepsilon (G) = \frac{2}{\frac{2}{k}-1+\frac{2}{l}}.$$

\no Whence, $\frac{2}{k}+\frac{2}{l}-1 > 0$. Since $k,l$ are both
integers and $k\geq 3, l\geq 3$, if $k\geq 6$, we get

$$\frac{2}{k}+\frac{2}{l}-1\leq\frac{2}{3}+\frac{2}{6}-1= 0.$$

\no Contradicts that $\frac{2}{k}+\frac{2}{l}-1 > 0$. Therefore,
$k\leq 5$. Similarly, $l\leq 5$. So we have $3\leq k\leq 5$ and
$3\leq l\leq 5$. Calculation shows that all possibilities for
$(k,l)$ are $(k,l) = (3,3), (3,4), (3,5), (4,3)$ and $(5,3)$. The
$(3,3)$ and $(3,4)$ polyhedrons have be shown in Fig.$2.12$ and
the remainder $(3,5), (4,3)$ and $(5,3)$ polyhedrons are shown in
Fig.$2.13$. \quad\quad $\natural$

An {\it elementary subdivision} on a graph $G$ is a graph obtained
from $G$ replacing an edge $e=uv$ by a path $uwv$, where,
$w\not\in V(G)$. A {\it subdivision} of $G$ is a graph obtained
from $G$ by a succession of elementary subdivision. A graph $H$ is
defined to be a {\it homeomorphism of $G$} if either $H\cong G$ or
$H$ is isomorphic to a subdivision of $G$. Kuratowski found the
following characterization for planar graphs in $1930$. For its a
complete proof, see $[9],[11]$ for details.

\vskip 4mm

\no{\bf Theorem $2.1.16$} \ {\it A graph is planar if and only if
it contains no subgraph homeomorphic with $K_5$ or $K(3,3)$.}

\vskip 5mm

\no{\bf $(2)$ Decomposition of graphs}

\vskip 4mm

\no A complete graph $K_6$ with vertex set $\{1,2,3,4,5,6\}$ has
two families of subgraphs $\{C_6,C_3^1, C_3^2, P_2^1, P_2^2,
P_2^3\}$ and $\{S_{1.5}, S_{1.4}, S_{1.3}, S_{1.2}, S_{1.1}\}$,
such as those shown in Fig.$2.14$ and Fig.$2.15$.

\includegraphics[bb=10 10 450 100]{sgm16.eps}

\vskip 2mm

\c{\bf Fig $2.14$}

\includegraphics[bb=0 10 450 140]{sgm17.eps}

\vskip 3mm

\c{\bf Fig $2.15$}

 \vskip 3mm

\no We know that

$$E(K_6) = E(C_6)\bigcup E(C_3^1)\bigcup E(C_3^2)\bigcup E(P_2^1)\bigcup E(P_2^2)\bigcup E(P_2^3);$$

$$E(K_6) = E(S_{1.5})\bigcup E(S_{1.4})\bigcup E(S_{1.3})\bigcup E(S_{1.2})\bigcup E(S_{1.1}).$$

\no These formulae imply the conception of decomposition of
graphs. For a graph $G$, a {\it decomposition} of $G$ is a
collection $\{H_i\}_{1\leq i\leq s}$ of subgraphs of $G$ such that
for any integer $i, 1\leq i\leq s$, $H_i=\left<E_i\right>$ for
some subsets $E_i$ of $E(G)$ and $\{E_i\}_{1\leq i\leq s}$ is a
partition of $E(G)$, denoted by $G=H_1\bigoplus
H_2\bigoplus\cdots\bigoplus H_s$. The following result is obvious.

\vskip 4mm

\no{\bf Theorem $2.1.17$} \ {\it Any graph $G$ can be decomposed
to bouquets and dipoles, in where $K_2$ is seen as a dipole
$D_{0.1.0}$.}

\vskip 4mm

\no{\bf Theorem $2.1.18$} \ {\it For every positive integer $n$,
the complete graph $K_{2n+1}$ can be decomposed to $n$ hamiltonian
circuits.}

\vskip 3mm

{\it Proof} \ For $n=1$, $K_3$ is just a hamiltonian circuit. Now
let $n\geq 2$ and $V(K_{2n+1})=\{v_0,v_1,v_2,\cdots , v_{2n}\}$.
Arrange these vertices $v_1,v_2,\cdots , v_{2n}$ on vertices of a
regular $2n$-gon and place $v_0$ in a convenient position not in
the $2n$-gon. For $i=1,2,\cdots,n$, we define the edge set of
$H_i$ to be consisted of $v_0v_i,v_0v_{n+i}$ and edges parallel to
$v_iv_{i+1}$ or edges parallel to $v_{i-1}v_{i+1}$, where the
subscripts are expressed modulo $2n$. Then we get that

$$K_{2n+1}= H_1\bigoplus H_2\bigoplus\cdots\bigoplus H_n$$

\no with each $H_i, 1\leq i\leq n$ being a hamiltonian circuit

$$v_0v_iv_{i+1}v_{i-1}v_{i+1}v_{i-2}\cdots v_{n+i-1}v_{n+i+1}v_{n+i}v_0. \quad\quad \natural$$

Every Cayley graph of a finite group $\Gamma$ can be decomposed
into $1$-factors or $2$-factors in a natural way as stated in the
following theorems.

\vskip 4mm

\no{{\bf Theorem $2.1.19$} {\it Let $G$ be a vertex-transitive
graph and let $H$ be a regular subgroup of ${\rm Aut}G$. Then for
any chosen vertex $x, x\in V(G)$, there is a factorization}

$$
G = ( \bigoplus\limits_{y\in N_G
(x),|H_{(x,y)}|=1}(x,y)^H)\bigoplus ( \bigoplus\limits_{y\in N_G
(x),|H_{(x,y)}|=2}(x,y)^H) ,
$$

\no{\it for $G$ such that $(x,y)^H$ is a $2$-factor if
$|H_{(x,y)}|=1$ and a $1$-factor if $|H_{(x,y)}|=2$.}} \vskip 3mm

{\it Proof} \ First, We prove the following claims.

\vskip 3mm

\no{\bf Claim $1$} \ {\it $\forall x\in V(G),x^H=V(G)$ and
$H_x=1_H$.}\vskip 3mm

\no{\bf Claim $2$} \ {\it For $\forall (x,y),(u,w)\in E(G)$,
$(x,y)^H \bigcap (u,w)^H=\emptyset$  or $(x,y)^H=(u,w)^H$.} \vskip
2mm

Claims $1$ and $2$ are holden by definition. \vskip 3mm

\no{\bf Claim $3$} \ {\it For $\forall (x,y)\in E(G),|H_{(x,y)}|=1
\quad or\quad 2$.}\vskip 2mm

Assume that $|H_{(x,y)}| \not= 1$. Since we know that
$(x,y)^h=(x,y)$, i.e., $(x^h,y^h)=(x,y)$ for any element $h\in
H_{(x,y)}$. Thereby we get that $x^h=x$ and $y^h=y$ or $x^h=y$ and
$y^h=x$. For the first case we know $h=1_H$ by Claim $1$. For the
second, we get that $x^{h^2}=x$. Therefore, $h^2=1_H$ .

Now if there exists an element $g\in H_{(x,y)}\backslash \{1_H,h
\}$, then we get $x^g=y= x^h$ and $y^g=x= y^h$. Thereby we get
$g=h$ by Claim $1$, a contradiction. So we get that
$|H_{(x,y)}|=2$.\vskip 3mm

\no{\bf Claim $4$} \ {\it For any $(x,y)\in E(G)$, if
$|H_{(x,y)}|=1$, then $(x,y)^H$ is a $2$-factor.}\vskip 2mm

Because $x^H=V(G)\subset V(\left<(x,y)^H\right>)\subset V(G)$, so
$V(\left<(x,y)^H\right>) = V(G)$. Therefore, $(x,y)^H$ is a
spanning subgraph of $G$.

Since $H$ acting on $V(G)$ is transitive, there exists an element
$h\in H$ such that $x^h=y$. It is obvious that $o(h)$ is finite
and $o(h)\not= 2$. Otherwise, we have $|H_{(x,y)}|\geq 2$, a
contradiction. Now $(x,y)^{\left<h\right>} =xx^hx^{h^2}\cdots
x^{h^{o(h)-1}}x$ is a circuit in the graph $G$. Consider the right
coset decomposition of $H$ on $\left<h\right>$. Suppose
$H=\bigcup\limits_{i=1}^s \left<h\right> a_i$, $\left<h\right> a_i
\bigcap \left<h\right> a_j=\emptyset$, if $i\not= j$, and
$a_1=1_H$.

Now let $X = \{a_1,a_2,...,a_s \}$. We know that for any $a,b\in
X, (\left<h\right> a) \bigcap(\left<h\right> b)= \emptyset$ if
$a\not= b$. Since $(x,y)^{\left<h\right>
a}=((x,y)^{\left<h\right>})^a$ and $(x,y)^{\left<h\right> b}=
((x,y)^{\left<h\right>})^b$ are also circuits, if
$V(\left<(x,y)^{\left<h\right> a}\right>)\bigcap
V(\left<(x,y)^{\left<h\right> b}\right>) \not=\emptyset$ for some
$a,b\in X, a\not= b$, then there must be two elements
$f,g\in\left<h\right>$ such that $x^{fa}=x^{gb}$ . According to
Claim $1$, we get that $fa=gb$, that is
$ab^{-1}\in\left<h\right>$. So $\left<h\right> a=\left<h\right> b$
and $a=b$, contradicts to the assumption that $a\not= b$.

Thereafter we know that $(x,y)^H=\bigcup\limits_{a\in X}
(x,y)^{\left<h\right> a}$ is a disjoint union of circuits. So
$(x,y)^H$ is a $2$-factor of the graph $G$.\vskip 3mm

\no{\bf Claim $5$} \ {\it For any $(x,y)\in E(G)$, $(x,y)^H$ is an
$1$-factor if $|H_{(x,y)}|=2$.} \vskip 2mm

Similar to the proof of Claim $4$, we know that
$V(\left<(x,y)^H\right>) = V(G)$ and $(x,y)^H$ is a spanning
subgraph of the graph $G$.

Let $H_{(x,y)}=\{1_H,h \}$, where $x^h=y$ and $y^h=x$. Notice that
$(x,y)^a=(x,y)$ for $\forall a\in H_{(x,y)}$. Consider the coset
decomposition of $H$ on $H_{(x,y)}$, we know that
$H=\bigcup\limits_{i=1}^t H_{(x,y)}b_i$ , where
$H_{(x,y)}b_i\bigcap H_{(x,y)}b_j=\emptyset$ if $i\not= j, 1\leq
i,j\leq t$. Now let $L=\{ H_{(x,y)}b_i , 1\leq i\leq t \}$. We get
a decomposition

$$(x,y)^H= \bigcup\limits_{b\in L} (x,y)^b$$

\no for $(x,y)^H$. Notice that if $b=H_{(x,y)}b_i \in L$,
$(x,y)^b$ is an edge of $G$. Now if there exist two elements
$c,d\in L, c= H_{(x,y)}f$ and $d= H_{(x,y)}g$, $f\not= g$ such
that $V(\left<(x,y)^c\right>)\bigcap$
$V(\left<(x,y)^d\right>)\not=\emptyset$, there must be $x^f=x^g$
or $x^f=y^g$. If $x^f=x^g$, we get $f=g$ by Claim $1$, contradicts
to the assumption that $f\not= g$. If $x^f=y^g=x^{hg}$, where
$h\in H_{(x,y)}$, we get $f=hg$ and $fg^{-1}\in H_{(x,y)}$, so
$H_{(x,y)}f= H_{(x,y)}g$. According to the definition of $L$, we
get $f=g$, also contradicts to the assumption that $f\not= g$.
Therefore, $(x,y)^H$ is an $1$-factor of the graph $G$.

Now we can prove the assertion in this theorem.  According to
Claim $1$- Claim $4$, we get that

$$
G = ( \bigoplus\limits_{y\in N_G
(x),|H_{(x,y)}|=1}(x,y)^H)\bigoplus ( \bigoplus\limits_{y\in N_G
(x),|H_{(x,y)}|=2}(x,y)^H) .
$$

\no for any chosen vertex $x, x\in V(G)$. By Claims $5$ and $6$,
we know that $(x,y)^H$ is a $2$-factor if $|H_{(x,y)}|=1$ and is a
$1$-factor if $|H_{(x,y)}|=2$. Whence, the desired factorization
for $G$ is obtained. \quad\quad $\natural$

Now for a Cayley graph ${\rm Cay}(\Gamma:S)$, by Theorem $2.1.13$,
we can always choose the vertex $x=1_{\Gamma}$ and $H$ the right
regular transformation group on $\Gamma$. After then, Theorem
$2.1.19$ can be restated as follows.

\vskip 4mm

\no{\bf Theorem $2.1.20$} \ {\it Let $\Gamma$ be a finite group
with a subset $S, S^{-1}=S$, $1_{\Gamma}\not\in S$ and $H$ is the
right transformation group on $\Gamma$. Then there is a
factorization}

$$
G = ( \bigoplus\limits_{s\in S,s^2\not=
1_{\Gamma}}(1_{\Gamma},s)^H)\bigoplus (\bigoplus\limits_{s\in
S,s^2 = 1_{\Gamma}} (1_{\Gamma},s)^H)
$$

\no{\it for the Cayley graph ${\rm Cay}(\Gamma:S)$ such that
$(1_{\Gamma},s)^H$ is a $2$-factor if $s^2\not=1_{\Gamma}$ and
$1$-factor if $s^2=1_{\Gamma}$.}

\vskip 3mm

{\it Proof} \ For any $h\in H_{(1_{\Gamma},s)}$, if $h\not=
1_{\Gamma}$, then we get that $1_{\Gamma}h=s$ and $sh =
1_{\Gamma}$, that is $s^2=1_{\Gamma}$. According to Theorem
$2.1.19$, we get the factorization for the Cayley graph ${\rm
Cay}(\Gamma:S)$. \quad\quad $\natural$

More factorial properties for Cayley graphs of a finite group can
be found in the reference $[51]$.

\vskip 5mm

\no{\bf $2.1.4.$ Operations on graphs}

\vskip 4mm

\no For two given graphs $G_1= (V_1.E_1;I_1)$ and $G_2 =
(V_2,E_2;I_2)$, there are a number of ways to produce new graphs
from $G_1$ and $G_2$. Some of them are described in the following.

\vskip 4mm

\no{\bf Operation $1.$ \ Union}

\vskip 3mm

\no The {\it union} $G_1\bigcup G_2$ of graphs $G_1$ and $G_2$ is
defined by

$$V(G_1\bigcup G_2)= V_1\bigcup V_2, \ E(G_1\bigcup G_2)= E_1\bigcup E_2 \ {\rm and} \
I(E_1\bigcup E_2)= I_1(E_1)\bigcup I_2(E_2).$$

\no If a graph consists of $k$ disjoint copies of a graph $H$,
$k\geq 1$, then we write $G = kH$. Therefore, we get that
$K_6=C_6\bigcup 3K_2\bigcup 2K_3 = \bigcup\limits_{i=1}^5S_{1.i}$
for graphs in Fig.$2.14$ and Fig.$2.15$ and generally, $K_n =
\bigcup\limits_{i=1}^{n-1}S_{1.i}$. For an integer $k, k\geq 2$
and a simple graph $G$, $kG$ is a multigraph with edge multiple
$k$ by definition.

By the definition of a union of two graphs, we get decompositions
for some well-known graphs such as

$$B_n = \bigcup_{i=1}^nB_1(O), \ \ \
D_{k,m,n} = (\bigcup\limits_{i=1}^kB_1(O_1))\bigcup
(\bigcup\limits_{i=1}^mK_2)\bigcup
(\bigcup\limits_{i=1}^nB_1(O_2)),$$

\no where $V(B_1)(O_1)=\{O_1\}, V(B_1)(O_2)=\{O_2\}$ and
$V(K_2)=\{O_1,O_2\}$. By Theorem $1.18$, we get that

$$K_{2n+1} = \bigcup\limits_{i=1}^n H_i$$

\no with $H_i=v_0v_iv_{i+1}v_{i-1}v_{i+1}v_{i-2}\cdots
v_{n+i-1}v_{n+i+1}v_{n+i}v_0$.

In Fig.$2.16$, we show two graphs $C_6$ and $K_4$ with a nonempty
intersection and their union $C_6\bigcup K_4$.

\includegraphics[bb=60 10 450 120]{sgm18.eps}

\vskip 3mm

\c{\bf Fig $2.16$}

\vskip 4mm

\no{\bf Operation $2.$ \ Join}

\vskip 3mm

\no The {\it complement} $\overline{G}$ of a graph $G$ is a graph
with the vertex set $V(G)$ such that two vertices are adjacent in
$\overline{G}$ if and only if these vertices are not adjacent in
$G$. The {\it join} $G_1+G_2$ of $G_1$ and $G_2$ is defined by

$$V(G_1+G_2)=V(G_1)\bigcup V(G_2),$$

$$E(G_1+G_2)= E(G_1)\bigcup E(G_2)\bigcup\{(u,v) |u\in V(G_1), v\in V(G_2)\}$$

\no and

$$I(G_1+G_2)= I(G_1)\bigcup I(G_2)\bigcup\{I(u,v)= (u,v) |u\in V(G_1), v\in V(G_2)\}.$$

\no Using this operation, we can represent
$K(m,n)\cong\overline{K_m}+\overline{K_n}$. The join graph of
circuits $C_3$ and $C_4$ is given in Fig.$2.17$.

\includegraphics[bb=60 10 450 120]{sgm19.eps}

\vskip 3mm

\c{\bf Fig $2.17$}

\vskip 4mm

\no{\bf Operation $3.$ \ Cartesian product}

\vskip 3mm

\no The {\it cartesian product} $G_1\times G_2$ of graphs $G_1$
and $G_2$ is defined by $V(G_1\times G_2)= V(G_1)\times V(G_2)$
and two vertices $(u_1,u_2)$ and $(v_1,v_2)$ of $G_1\times G_2$
are adjacent if and only if either $u_1=v_1$ and $(u_2,v_2)\in
E(G_2)$ or $u_2=v_2$ and $(u_1,v_1)\in E(G_1)$.

For example, the cartesian product $C_3\times C_3$ of circuits
$C_3$ and $C_3$ is shown in Fig.$2.18$.

\includegraphics[bb=60 10 450 150]{sgm20.eps}

\vskip 3mm

\c{\bf Fig $2.18$}

\vskip 8mm

\no{\bf \S $2.2$ \ Multi-Voltage Graphs}

\vskip 5mm

\no There is a convenient way for constructing a covering space of
a graph $G$ in topological graph theory, i.e., by a voltage graph
$(G, \alpha)$ of $G$ which was firstly introduced by Gustin in
1963 and then generalized by Gross in 1974. Youngs extensively
used voltage graphs in proving Heawood map coloring
theorem([$23$]). Today, it has become a convenient way for finding
regular maps on surface. In this section, we generalize voltage
graphs to two types of multi-voltage graphs by using finite
multi-groups.

\vskip 4mm

\no{\bf $2.2.1.$ Type $1$}

\vskip 4mm

\no{\bf Definition $2.2.1$} \ {\it Let $\widetilde{\Gamma} =
\bigcup\limits_{i=1}^n\Gamma_i$ be a finite multi-group with an
operation set $O(\widetilde{\Gamma})=\{\circ_i | 1\leq i\leq n\}$
and $G$ a graph. If there is a mapping $\psi:
X_{\frac{1}{2}}(G)\rightarrow\widetilde{\Gamma}$ such that
$\psi(e^{-1}) = (\psi (e^+))^{-1}$ for $\forall e^+\in
X_{\frac{1}{2}}(G)$, then $(G, \psi )$ is called a multi-voltage
graph of type $1$.}

\vskip 3mm

Geometrically, a multi-voltage graph is nothing but a weighted
graph with weights in a multi-group. Similar to voltage graphs,
the importance of a multi-voltage graph is in its {\it lifting}
defined in the next definition.

\vskip 4mm

\no{\bf Definition $2.2.2$} \ {\it For a multi-voltage graph
$(G,\psi )$ of type $1$, the lifting graph $G^{\psi}=
(V(G^{\psi}),$ $E(G^{\psi});I(G^{\psi}))$ of $(G,\psi )$ is
defined by}

$$V(G^{\psi}) = V(G)\times\widetilde{\Gamma},$$

$$E(G^{\psi}) = \{(u_a, v_{a\circ b}) | e^+ =(u,v)\in X_{\frac{1}{2}}(G),
\psi (e^+) =b, a\circ b\in\widetilde{\Gamma} \}$$

\no{\it and}

$$I(G^{\psi})= \{(u_a, v_{a\circ b}) | I(e)=(u_a, v_{a\circ b}) \ if \
e=(u_a, v_{a\circ b})\in E(G^{\psi})\}.$$

\vskip 3mm

For abbreviation, a vertex $(x,g)$ in $G^{\psi}$ is denoted by
$x_g$. Now for $\forall v\in V(G)$, $v\times\widetilde{\Gamma} =
\{v_g | g\in\widetilde{\Gamma}\}$ is called a {\it fiber over
$v$}, denoted by $F_v$. Similarly, for $\forall e^+ = (u,v)\in
X_{\frac{1}{2}}(G)$ with $\psi (e^+)= b$, all edges
$\{(u_g,v_{g\circ b}) | g, g\circ b\in\widetilde{\Gamma} \}$ is
called the {\it fiber over $e$}, denoted by $F_e$.

For a multi-voltage graph $(G,\psi )$ and its lifting $G^{\psi}$,
there is a {\it natural projection} $p: G^{\psi}\rightarrow G$
defined by $p(F_v)=v$ for $\forall v\in V(G)$. It can be verfied
that $p(F_e)=e$ for $\forall e\in E(G)$.

Choose $\widetilde{\Gamma}=\Gamma_1\bigcup \Gamma_2$ with
$\Gamma_1 = \{1,a,a^2\}$, $\Gamma_2=\{1,b,b^2\}$ and $a\not=b$. A
multi-voltage graph and its lifting are shown in Fig.$2.19$.

\includegraphics[bb=60 10 450 140]{sgm21.eps}

\vskip 3mm

\c{\bf Fig $2.19$}\vskip 2mm

Let $\widetilde{\Gamma} = \bigcup\limits_{i=1}^n\Gamma_i$ be a
finite multi-group with  groups $(\Gamma_i;\circ_i), 1\leq i\leq
n$. Similar to the unique walk lifting theorem for voltage graphs,
we know the following {\it walk multi-lifting theorem} for
multi-voltage graphs of type $1$.

\vskip 4mm

\no{\bf Theorem $2.2.1$} \ {\it Let $W = e^1e^2\cdots e^k$ be a
walk in a multi-voltage graph $(G,\psi)$ with initial vertex $u$.
Then there exists a lifting $W^{\psi}$ start at $u_a$ in
$G^{\psi}$ if and only if there are integers $i_1,i_2,\cdots ,i_k$
such that}

$$a\circ_{i_1}\psi (e_1^+)\circ_{i_2}\cdots\circ_{i_{j-1}}\psi (e_j^+)\in \Gamma_{i_{j+1}}
\ and \ \psi (e_{j+1}^+)\in \Gamma_{i_{j+1}}$$

\no{\it for any integer $j, 1\leq j\leq k$}

\vskip 3mm

{\it Proof} \ Consider the first semi-arc in the walk $W$, i.e.,
$e_1^+$. Each lifting of $e_1$ must be $(u_a,u_{a\circ \psi
(e_1^+)})$. Whence, there is a lifting of $e_1$ in $G^{\psi}$ if
and only if there exists an integer $i_1$ such that $\circ =
\circ_{i_1}$ and $a, a\circ_{i_1}\psi (e_1^+)\in \Gamma_{i_1}$.

Now if we have proved there is a lifting of a sub-walk
$W_l=e_1e_2\cdots e_l$ in $G^{\psi}$ if and only if there are
integers $i_1,i_2,\cdots , i_l$, $1\leq l < k$ such that

$$a\circ_{i_1}\psi (e_1^+)\circ_{i_2}\cdots\circ_{i_{j-1}}\psi (e_j^+)\in \Gamma_{i_{j+1}},
\ \ \ \psi (e_{j+1}^+)\in \Gamma_{i_{j+1}}$$

\no for any integer $j, 1\leq j\leq l$, we consider the semi-arc
$e_{l+1}^+$. By definition, there is a lifting of $e_{l+1}^+$ in
$G^{\psi}$ with initial vertex $u_{a\circ_{i_1}\psi
(e_1^+)\circ_{i_2}\cdots\circ_{i_{j-1}}\psi (e_l^+)}$ if and only
if there exists an integer $i_{l+1}$ such that

$$a\circ_{i_1}\psi
(e_1^+)\circ_{i_2}\cdots\circ_{i_{j-1}}\psi (e_l^+)\in
\Gamma_{l+1} \ {\rm and} \ \psi (e_{l+1}^+)\in \Gamma_{l+1}.$$

According to the induction principle, we know that there exists a
lifting $W^{\psi}$ start at $u_a$ in $G^{\psi}$ if and only if
there are integers $i_1,i_2,\cdots ,i_k$ such that

$$a\circ_{i_1}\psi (e_1^+)\circ_{i_2}\cdots\circ_{i_{j-1}}\psi (e_j^+)\in \Gamma_{i_{j+1}},
\ {\rm and} \ \psi (e_{j+1}^+)\in \Gamma_{i_{j+1}}$$

\no for any integer $j, 1\leq j\leq k$.\quad\quad $\natural$

For two elements $g,h\in\widetilde{\Gamma}$, if there exist
integers $i_1,i_2,\cdots , i_k$ such that
$g,h\in\bigcap\limits_{j=1}^k\Gamma_{i_j}$ but for $\forall
i_{k+1}\in\{1,2,\cdots ,n\}\setminus\{i_1,i_2,\cdots ,i_k\}$, $g,
h\not\in\bigcap\limits_{j=1}^{k+1}\Gamma_{i_j}$, we call $k= \Pi
[g,h]$ the {\it joint number of $g$ and $h$}. Denote $O
(g,h)=\{\circ_{i_j}; 1\leq j\leq k\}$. Define
$\widetilde{\Pi}[g,h]=\sum\limits_{\circ\in
O(\widetilde{\Gamma})}\Pi [g,g\circ h]$, where $\Pi [g,g\circ
h]=\Pi [g\circ h,h]=0$ if $g\circ h$ does not exist in
$\widetilde{\Gamma}$. According to Theorem $2.2.1$, we get an
upper bound for the number of liftings in $G^{\psi}$ for a walk
$W$ in $(G,\psi)$.

\vskip 4mm

\no{\bf Corollary $2.2.1$} \ {\it If those conditions in Theorem
$2.2.1$ hold, the number of liftings of $W$ with initial vertex
$u_a$ in $G^{\psi}$ is not in excess of}

\begin{eqnarray*}
& \ &\widetilde{\Pi}[a,\psi (e_1^+)]\times\\
& \ &\prod\limits_{i=1}^k\sum\limits_{\circ_1\in O(a,\psi
(e_1^+))}\cdots\sum\limits_{\circ_i\in O(a;\circ_j,\psi
(e_j^+),1\leq j\leq i-1)}\widetilde{\Pi}[a\circ_1\psi
(e_1^+)\circ_2\cdots\circ_i\psi (e_i^+),\psi (e_{i+1}^+)],
\end{eqnarray*}

\no{\it where $O(a;\circ_j,\psi (e_j^+),1\leq j\leq i-1)=
O(a\circ_1\psi (e_1^+)\circ_2\cdots\circ_{i-1}\psi
(e_{i-1}^+),\psi (e_i^+))$.}\vskip 3mm

The natural projection of a multi-voltage graph is not regular in
general. For finding a regular covering of a graph, a typical
class of multi-voltage graphs is the case of $\Gamma_i= \Gamma$
for any integer $i, 1\leq i\leq n$ in these multi-groups
$\widetilde{\Gamma}= \bigcup\limits_{i=1}^n\Gamma_i$. In this
case, we can find the exact number of liftings in $G^{\psi}$ for a
walk in $(G,\psi)$.

\vskip 4mm

\no{\bf Theorem $2.2.2$} \ {\it Let $\widetilde{\Gamma} =
\bigcup\limits_{i=1}^n\Gamma$ be a finite multi-group with groups
$(\Gamma ;\circ_i), 1\leq i\leq n$ and let $W = e^1e^2\cdots e^k$
be a walk in a multi-voltage graph $(G,\psi)$ , $\psi :
X_{\frac{1}{2}}(G)\rightarrow\widetilde{\Gamma}$ of type $1$ with
initial vertex $u$. Then there are $n^k$ liftings of $W$ in
$G^{\psi}$ with initial vertex $u_a$ for $\forall
a\in\widetilde{\Gamma}$.}

\vskip 3mm

{\it Proof} \ The existence of lifting of $W$ in $G^{\psi}$ is
obvious by Theorem $2.2.1$. Consider the semi-arc $e_1^+$. Since
$\Gamma_i=\Gamma$ for $1\leq i\leq n$, we know that there are $n$
liftings of $e_1$ in $G^{\psi}$ with initial vertex $u_a$ for any
$a\in\widetilde{\Gamma}$, each with a form $(u_a,u_{a\circ\psi
(e_1^+)}), \circ\in O(\widetilde{\Gamma})$.

Now if we have gotten $n^s, 1\leq s\leq k-1$ liftings in
$G^{\psi}$ for a sub-walk $W_s= e^1e^2\cdots e^s$. Consider the
semi-arc $e_{s+1}^+$. By definition we know that there are also
$n$ liftings of $e_{s+1}$ in $G^{\psi}$ with initial vertex
$u_{a\circ_{i_1}\psi (e_1^+)\circ_{i_2}\cdots\circ_{s}\psi
(e_{s}^+)}$, where $\circ_i\in O(\widetilde{\Gamma}), 1\leq i\leq
s$. Whence, there are $n^{s+1}$ liftings in $G^{\psi}$ for a
sub-walk $W_s= e^1e^2\cdots e^{s+1}$ in $(G;\psi )$.

By the induction principle, we know the assertion is true.
\quad\quad $\natural$

\vskip 4mm

\no{\bf Corollary $2.2.2$}([$23$]) \ {\it Let $W$ be a walk in a
voltage graph $(G,\psi), \psi
:X_{\frac{1}{2}}(G)\rightarrow\Gamma$ with initial vertex $u$.
Then there is an unique lifting of $W$ in $G^{\psi}$ with initial
vertex $u_a$ for $\forall a\in \Gamma$.}

\vskip 3mm

If a lifting $W^{\psi}$ of a multi-voltage graph $(G,\psi)$ is the
same as the lifting of a voltage graph $(G,\alpha),
\alpha:X_{\frac{1}{2}}(G)\rightarrow\Gamma_i$, then this lifting
is called a {\it homogeneous lifting of $\Gamma_i$}. For lifting a
circuit in a multi-voltage graph, we get the following result.

\vskip 4mm

\no{\bf Theorem $2.2.3$} \ {\it Let $\widetilde{\Gamma} =
\bigcup\limits_{i=1}^n\Gamma$ be a finite multi-group with groups
$(\Gamma ;\circ_i), 1\leq i\leq n$, $C = u_1u_2\cdots u_mu_1$ a
circuit in a multi-voltage graph $(G,\psi)$ and  $\psi :
X_{\frac{1}{2}}(G)\rightarrow\widetilde{\Gamma}$. Then there are
$\frac{|\Gamma |}{o(\psi (C, \circ_i))}$ homogenous liftings of
length $o(\psi (C, \circ_i))m$ in $G^{\psi}$ of $C$ for any
integer $i, 1\leq i\leq n$, where $\psi (C, \circ_i)= \psi
(u_1,u_2)\circ_i\psi (u_2,u_3)\circ_i\cdots\circ_i\psi
(u_{m-1},u_m)\circ_i\psi (u_m,u_1)$ and there are }

$$\sum\limits_{i=1}^n\frac{|\Gamma |}{o(\psi (C, \circ_i))}$$

\no{\it homogenous liftings of $C$ in $G^{\psi}$ altogether.}

\vskip 3mm

{\it Proof} \ According to Theorem $2.2.2$, there are liftings
with initial vertex $(u_1)_a$ of $C$ in $G^{\psi}$ for $\forall
a\in\widetilde{\Gamma}$. Whence, for any integer $i, 1\leq i\leq
n$, walks

$$W_a=(u_1)_a(u_2)_{a\circ_i\psi (u_1,u_2)}\cdots
(u_m)_{a\circ_i\psi (u_1,u_2) \circ_i\cdots\circ_i\psi
(u_{m-1},u_m)}(u_1)_{a\circ_i\psi (C, \circ_i)},$$

\begin{eqnarray*}W_{a\circ_i\psi (C, \circ_i)}&=& (u_1)_{a\circ_i\psi (C,
\circ_i)}(u_2)_{a\circ_i\psi (C, \circ_i)\circ_i\psi (u_1,u_2)}\\
&\cdots& (u_m)_{a\circ_i\psi (C, \circ_i)\circ_i\psi (u_1,u_2)
\circ_i\cdots\circ_i\psi (u_{m-1},u_m)}(u_1)_{a\circ_i\psi^2 (C,
\circ_i)},\end{eqnarray*}

$$\cdots\cdots\cdots\cdots\cdots\cdots\cdots,$$

\no and

\begin{eqnarray*}W_{a\circ_i\psi^{o(\psi (C, \circ_i))-1} (C,
\circ_i)}&=&(u_1)_{a\circ_i\psi^{o(\psi (C, \circ_i))-1} (C,
\circ_i)}(u_2)_{a\circ_i\psi^{o(\psi (C, \circ_i))-1} (C,
\circ_i)\circ_i\psi (u_1,u_2)}\\
&\cdots& (u_m)_{a\circ_i\psi^{o(\psi (C, \circ_i))-1} (C,
\circ_i)\circ_i\psi (u_1,u_2) \circ_i\cdots\circ_i\psi
(u_{m-1},u_m)}(u_1)_a \end{eqnarray*}

\no are attached end-to-end to form a circuit of length $o(\psi
(C, \circ_i))m$. Notice that there are $\frac{|\Gamma |}{o(\psi
(C, \circ_i))}$ left cosets of the cyclic group generated by $\psi
(C, \circ_i)$ in the group $(\Gamma ,\circ_i)$ and each is
correspondent with a homogenous lifting of $C$ in $G^{\psi}$.
Therefore, we get

$$\sum\limits_{i=1}^n\frac{|\Gamma |}{o(\psi (C, \circ_i))}$$

\no homogenous liftings of $C$ in $G^{\psi}$. \quad\quad
$\natural$

\vskip 4mm

\no{\bf Corollary $2.2.3$}([$23$]) \ {\it Let $C$ be a $k$-circuit
in a voltage graph $(G,\psi )$ such that the order of $\psi (C,
\circ)$ is $m$ in the voltage group $(\Gamma ;\circ )$. Then each
component of the preimage $p^{-1}(C)$ is a $km$-circuit, and there
are $\frac{|\Gamma |}{m}$ such components.}

\vskip 3mm

The lifting $G^{\zeta}$ of a multi-voltage graph $(G,\zeta)$ of
type $1$ has a natural decomposition described in the next result.

\vskip 4mm

\no{\bf Theorem $2.2.4$} \ {\it  Let $(G, \zeta), \zeta
:X_{\frac{1}{2}}(G)\rightarrow\widetilde{\Gamma}=\bigcup\limits_{i=1}^n\Gamma_i$,
be a multi-voltage graph of type $1$. Then}

$$G^{\zeta} =\bigoplus\limits_{i=1}^nH_i,$$

\no{\it where $H_i$ is an induced subgraph $\left<E_i\right>$ of
$G^{\zeta}$ for an integer $i, 1\leq i\leq n$ with}

$$E_i=\{ (u_a, v_{a\circ_ib}) | a,b\in\Gamma_i \ and \ (u,v)\in E(G), \zeta (u,v)=b \}.$$

For a finite multi-group $\widetilde{\Gamma}
=\bigcup\limits_{i=1}^n\Gamma_i$ with an operation set
$O(\widetilde{\Gamma}) = \{\circ_i, 1\leq i\leq n\}$ and a graph
$G$, if there exists a decomposition $G
=\bigoplus\limits_{j=1}^nH_i$ and we can associate each element
$g_i\in\Gamma_i$ a homeomorphism $\varphi_{g_i}$ on the vertex set
$V(H_i)$ for any integer $i, 1\leq i\leq n$ such that

($i$) \ $\varphi_{g_i\circ_ih_i} =
\varphi_{g_i}\times\varphi_{h_i}$ for all $g_i, h_i\in\Gamma_i$,
where ¡°$\times$ ¡±is an operation between homeomorphisms;

($ii$) \ $\varphi_{g_i}$ is the identity homeomorphism if and only
if $g_i$ is the identity element of the group $(\Gamma_i;
\circ_i)$,

\no then we say this association to be a {\it subaction of a
multi-group $\widetilde{\Gamma}$ on the graph $G$}. If there
exists a subaction of $\widetilde{\Gamma}$ on $G$ such that
$\varphi_{g_i}(u) = u$ only if $g_i= {\bf 1}_{\Gamma_i}$ for any
integer $i, 1\leq i\leq n$, $g_i\in\Gamma_i$ and $u\in V_i$, then
we call it a {\it fixed-free subaction}.\vskip 3mm

A {\it left subaction $lA$} of $\widetilde{\Gamma}$ on $G^{\psi}$
is defined as follows:

{\it For any integer $i, 1\leq i\leq n$, let $V_i=\{u_a | u\in
V(G), a\in\widetilde{\Gamma}\}$ and $g_i\in\Gamma_i$. Define
$lA(g_i)(u_a) = u_{g_i\circ_ia}$ if $a\in V_i$. Otherwise,
$g_i(u_a) = u_a$.}\vskip 2mm

\no Then the following result holds.

\vskip 4mm

\no{\bf Theorem $2.2.5$} \ {\it Let $(G, \psi )$ be a
multi-voltage graph with $\psi:
X_{\frac{1}{2}}(G)\rightarrow\widetilde{\Gamma}=\bigcup\limits_{i=1}^n\Gamma_i$
and $G =\bigoplus\limits_{j=1}^nH_i$ with $H_i=\left<E_i\right>$,
$1\leq i\leq n$, where $E_i=\{ (u_a, v_{a\circ_ib}) |
a,b\in\Gamma_i \ and \ (u,v)\in E(G), \zeta (u,v)=b \}$. Then for
any integer $i, 1\leq i\leq n$,}

($i$) \ {\it for $\forall g_i\in\Gamma_i$, the left subaction
$lA(g_i)$ is a fixed-free subaction of an automorphism of $H_i$;}

($ii$) \ {\it $\Gamma_i$ is an automorphism group of $H_i$}.

\vskip 3mm

{\it Proof} \ Notice that $lA(g_i)$ is a one-to-one mapping on
$V(H_i)$ for any integer $i, 1\leq i\leq n$, $\forall
g_i\in\Gamma_i$. By the definition of a lifting, an edge in $H_i$
has the form $(u_a,v_{a\circ_ib})$ if $a,b\in\Gamma_i$. Whence,

$$(lA(g_i)(u_a), lA(g_i)(v_{a\circ_ib}))= (u_{g_i\circ_ia}, v_{g_i\circ_ia\circ_ib})\in E(H_i).$$

\no As a result, $lA(g_i)$ is an automorphism of the graph $H_i$.

Notice that $lA: \Gamma_i\rightarrow {\rm Aut}H_i$ is an injection
from $\Gamma_i$ to ${\rm Aut}G^{\psi}$. Since
$lA(g_i)\not=lA(h_i)$ for $\forall g_i, h_i\in\Gamma_i,
g_i\not=h_i, 1\leq i\leq n$. Otherwise, if $lA(g_i)=lA(h_i)$ for
$\forall a\in \Gamma_i$, then $g_i\circ_ia=h_i\circ_ia$. Whence,
$g_i=h_i$, a contradiction. Therefore, $\Gamma_i$ is an
automorphism group of $H_i$.

For any integer $i, 1\leq i\leq n$, $g_i\in\Gamma_i$,  it is
implied by definition that $lA(g_i)$ is a fixed-free subaction on
$G^{\psi}$. This completes the proof. \quad\quad $\natural$

\vskip 4mm

\no{\bf Corollary $2.2.4$}([$23$]) \ {\it Let $(G,\alpha )$ be a
voltage graph with $\alpha : X_{\frac{1}{2}}(G)\rightarrow\Gamma$.
Then $\Gamma$ is an automorphism group of $G^{\alpha}$.}

For a finite multi-group
$\widetilde{\Gamma}=\bigcup\limits_{i=1}^n\Gamma_i$ action on a
graph $\widetilde{G}$, the vertex orbit $orb(v)$ of a vertex $v\in
V(\widetilde{G})$ and the edge orbit $orb(e)$ of an edge $e\in
E(\widetilde{G})$ are defined as follows:

$$orb(v) = \{ g(v) | g\in\widetilde{\Gamma} \} \ {\rm and} \ orb(e) = \{ g(e) | g\in\widetilde{\Gamma} \}.$$

The {\it quotient graph $\widetilde{G}/\widetilde{\Gamma}$} of
$\widetilde{G}$ under the action of $\widetilde{\Gamma}$ is
defined by

$$V(\widetilde{G}/\widetilde{\Gamma})=\{ \ orb(v) \ | \ v\in
V(\widetilde{G})\}, \ \
E(\widetilde{G}/\widetilde{\Gamma})=\{orb(e) | e\in
E(\widetilde{G})\}$$

\no and

$$I(orb(e))= (orb(u),orb(v)) \ if \ there \ exists \ (u,v)\in
E(\widetilde{G})$$.

For example, a quotient graph is shown in Fig.$2.20$, where,
$\widetilde{\Gamma}= Z_5$.

\includegraphics[bb=10 10 400 90]{sgm22.eps}

\vskip 3mm

\c{\bf Fig $2.20$}\vskip 3mm

\no Then we get a necessary and sufficient condition for the
lifting of a multi-voltage graph in next result.

\vskip 4mm

\no{\bf Theorem $2.2.6$} \ {\it If the subaction ${\mathcal A}$ of
a finite multi-group
$\widetilde{\Gamma}=\bigcup\limits_{i=1}^n\Gamma_i$ on a graph
$\widetilde{G}=\bigoplus\limits_{i=1}^nH_i$ is fixed-free, then
there is a multi-voltage graph
$(\widetilde{G}/\widetilde{\Gamma},\zeta )$, $\zeta:
X_{\frac{1}{2}}(\widetilde{G}/\widetilde{\Gamma})
\rightarrow\widetilde{\Gamma}$ of type $1$ such that}

$$\widetilde{G} \ \cong \ (\widetilde{G}/\widetilde{\Gamma})^{\zeta}.$$

\vskip 3mm

{\it Proof} \  First, we choose positive directions for edges of
$\widetilde{G}/\widetilde{\Gamma}$ and $\widetilde{G}$ so that the
quotient map
$q_{\widetilde{\Gamma}}:\widetilde{G}\rightarrow\widetilde{G}/\widetilde{\Gamma}$
is direction-preserving and that the action ${\mathcal A}$ of
$\widetilde{\Gamma}$ on $\widetilde{G}$ preserves directions.
Next, for any integer $i, 1\leq i\leq n$ and $\forall v\in
V(\widetilde{G}/\widetilde{\Gamma})$, label one vertex of the
orbit $q_{\widetilde{\Gamma}}^{-1}(v)$ in $\widetilde{G}$ as
$v_{1_{\Gamma_i}}$ and for every group element $g_i\in\Gamma_i,
g_i\not=1_{\Gamma_i}$, label the vertex ${\mathcal
A}(g_i)(v_{1_{\Gamma_i}})$ as $v_{g_i}$. Now if the edge $e$ of
$\widetilde{G}/\widetilde{\Gamma}$ runs from $u$ to $w$, we
assigns the label $e_{g_i}$ to the edge of the orbit
$q_{\Gamma_i}^{-1}(e)$ that originates at the vertex $u_{g_i}$.
Since $\Gamma_i$ acts freely on $H_i$, there are just $|\Gamma_i|$
edges in the orbit $q_{\Gamma_i}^{-1}(e)$ for each integer $i,
1\leq i\leq n$, one originating at each of the vertices in the
vertex orbit $q_{\Gamma_i}^{-1}(v)$. Thus the choice of an edge to
be labelled $e_{g_i}$ is unique for any integer $i, 1\leq i\leq
n$. Finally, if the terminal vertex of the edge $e_{1_{\Gamma_i}}$
is $w_{h_i}$, one assigns a voltage $h_i$ to the edge $e$ in the
quotient $\widetilde{G}/\widetilde{\Gamma}$, which enables us to
get a multi-voltage graph $(\widetilde{G}/\widetilde{\Gamma},\zeta
)$. To show that this labelling of edges in $q_{\Gamma_i}^{-1}(e)$
and the choice of voltages $h_i, 1\leq i\leq n$ for the edge $e$
really yields an isomorphism $\vartheta :\widetilde{G}\rightarrow
(\widetilde{G}/\widetilde{\Gamma})^{\zeta}$, one needs to show
that for $\forall g_i\in\Gamma_i, 1\leq i\leq n$ that the edge
$e_{g_i}$ terminates at the vertex $w_{g_i\circ_ih_i}$. However,
since $e_{g_i}={\mathcal A}(g_i)(e_{1_{\Gamma_i}})$, the terminal
vertex of the edge $e_{g_i}$ must be the terminal vertex of the
edge ${\mathcal A}(g_i)(e_{1_{\Gamma_i}})$, which is

$${\mathcal A}(g_i)(w_{h_i})={\mathcal A}(g_i){\mathcal A}(h_i)(w_{1_{\Gamma_i}})=
{\mathcal A}(g_i\circ_ih_i)(w_{1_{\Gamma_i}})=w_{g_i\circ_ih_i}.$$

\no Under this labelling process, the isomorphism $\vartheta
:\widetilde{G}\rightarrow
(\widetilde{G}/\widetilde{\Gamma})^{\zeta}$ identifies orbits in
$\widetilde{G}$ with fibers of $G^{\zeta}$. Moreover, it is
defined precisely so that the action of $\widetilde{\Gamma}$ on
$\widetilde{G}$ is consistent with the left subaction $lA$ on the
lifting graph $G^{\zeta}$. This completes the proof.\quad\quad
$\natural$

\vskip 4mm

\no{\bf Corollary $2.2.5$}([$23$]) \ {\it  Let $\Gamma$ be a group
acting freely on a graph $\widetilde{G}$ and let $G$ be the
resulting quotient graph. Then there is an assignment $\alpha$ of
voltages in $\Gamma$ to the quotient graph $G$ and a labelling of
the vertices $\widetilde{G}$ by the elements of $V(G)\times\Gamma$
such that $\widetilde{G} = G^{\alpha}$ and that the given action
of $\Gamma$ on $\widetilde{G}$ is the natural left action of
$\Gamma$ on $G^{\alpha}$.}

\vskip 4mm

\no{\bf $2.2.2.$ Type $2$}

\vskip 4mm

\no{\bf Definition $2.2.3$} \ {\it Let $\widetilde{\Gamma} =
\bigcup\limits_{i=1}^n\Gamma_i$ be a finite multi-group and let
$G$ be a graph with vertices partition
$V(G)=\bigcup\limits_{i=1}^nV_i$. For any integers $i,j, 1\leq
i,j\leq n$, if there is a mapping $\tau:
X_{\frac{1}{2}}(\left<E_G(V_i,V_j)\right>)\rightarrow\Gamma_i\bigcap\Gamma_j$
and $\varsigma :V_i\rightarrow\Gamma_i$ such that $\tau (e^{-1}) =
(\tau (e^+))^{-1}$ for $\forall e^+\in X_{\frac{1}{2}}(G)$ and the
vertex subset $V_i$ is associated with the group $(\Gamma_i,
\circ_i)$ for any integer $i, 1\leq i\leq n$, then $(G, \tau
,\varsigma )$ is called a multi-voltage graph of type $2$.}

\vskip 3mm

Similar to multi-voltage graphs of type $1$, we construct a
lifting from a multi-voltage graph of type $2$.

\vskip 4mm

\no{\bf Definition $2.2.4$} \ {\it For a multi-voltage graph
$(G,\tau ,\varsigma )$ of type $2$, the lifting graph $G^{(\tau ,
\varsigma )}= (V(G^{(\tau , \varsigma )}),$ $E(G^{(\tau ,
\varsigma )});I(G^{(\tau , \varsigma )}))$ of $(G,\tau , \varsigma
)$ is defined by}

$$V(G^{(\tau ,
\varsigma )}) = \bigcup\limits_{i=1}^n\{V_i\times\Gamma_i\},$$

$$E(G^{(\tau ,
\varsigma )}) = \{(u_a, v_{a\circ b}) | e^+ =(u,v)\in
X_{\frac{1}{2}}(G), \psi (e^+) =b, a\circ b\in\widetilde{\Gamma}
\}$$

\no{\it and}

$$I(G^{(\tau ,
\varsigma )})= \{(u_a, v_{a\circ b}) | I(e)=(u_a, v_{a\circ b}) \
if \ e=(u_a, v_{a\circ b})\in E(G^{(\tau , \varsigma )})\}.$$

\vskip 3mm

Two multi-voltage graphs of type $2$ are shown on the left and
their lifting on the right in (a) and (b) of Fig.$21$. In where,
$\widetilde{\Gamma}= Z_2\bigcup Z_3$, $V_1=\{u\}$, $V_2=\{v\}$ and
$\varsigma :V_1\rightarrow Z_2$, $\varsigma :V_2\rightarrow Z_3$.

\includegraphics[bb=5 10 400 200]{sgm23.eps}

\vskip 3mm

\c{\bf Fig $2.21$}\vskip 4mm

\no{\bf Theorem $2.2.7$} \ {\it Let $(G,\tau ,\varsigma )$ be a
multi-voltage graph of type $2$ and let $W_k = u_1u_2\cdots u_k$
be a walk in $G$. Then there exists a lifting of $W^{(\tau
,\varsigma )}$ with an initial vertex $(u_1)_a,
a\in\varsigma^{-1}(u_1)$ in $G^{(\tau ,\varsigma )}$ if and only
if $a\in\varsigma^{-1}(u_1)\bigcap\varsigma^{-1}(u_2)$ and for any
integer $s, 1\leq s < k$, $a\circ_{i_1}\tau
(u_1u_2)\circ_{i_2}\tau
(u_2u_3)\circ_{i_3}\cdots\circ_{i_{s-1}}\tau
(u_{s-2}u_{s-1})\in\varsigma^{-1}(u_{s-1})\bigcap\varsigma^{-1}(u_s)$,
where ¡°$\circ_{i_j}$¡± is an operation in the group
$\varsigma^{-1}(u_{j+1})$ for any integer $j, 1\leq j\leq s$.}

\vskip 3mm

{\it Proof} \ By the definition of the lifting of a multi-voltage
graph of type $2$, there exists a lifting of the edge $u_1u_2$ in
$G^{(\tau ,\varsigma )}$ if and only if $a\circ_{i_1}\tau
(u_1u_2)\in\varsigma^{-1}(u_2)$, where ¡°$\circ_{i_j}$¡± is an
operation in the group $\varsigma^{-1}(u_2)$. Since $\tau
(u_1u_2)\in\varsigma^{-1}(u_1)\bigcap\varsigma^{-1}(u_2)$, we get
that $a\in\varsigma^{-1}(u_1)\bigcap\varsigma^{-1}(u_2)$.
Similarly, there exists a lifting of the subwalk $W_2= u_1u_2u_3$
in $G^{(\tau ,\varsigma )}$ if and only if
$a\in\varsigma^{-1}(u_1)\bigcap\varsigma^{-1}(u_2)$ and
$a\circ_{i_1}\tau
(u_1u_2)\in\varsigma^{-1}(u_2)\bigcap\varsigma^{-1}(u_3)$.

Now assume there exists a lifting of the subwalk $W_l=
u_1u_2u_3\cdots u_l$ in $G^{(\tau ,\varsigma )}$ if and only if
 $a\circ_{i_1}\tau
(u_1u_2)\circ_{i_2}\cdots\circ_{i_{t-2}}\tau
(u_{t-2}u_{t-1})\in\varsigma^{-1}
(u_{t-1})\bigcap\varsigma^{-1}(u_t)$ for any integer $t, 1\leq
t\leq l$, where¡°$\circ_{i_j}$¡±is an operation in the group
$\varsigma^{-1}(u_{j+1})$ for any integer $j, 1\leq j\leq l$. We
consider the lifting of the subwalk $W_{l+1}=u_1u_2u_3\cdots
u_{l+1}$. Notice that if there exists a lifting of the subwalk
$W_l$ in $G^{(\tau ,\varsigma )}$, then the terminal vertex of
$W_l$ in $G^{(\tau ,\varsigma )}$ is $(u_l)_{a\circ_{i_1}\tau
(u_1u_2)\circ_{i_2}\cdots\circ_{i_{l-1}}\tau (u_{l-1}u_l)}$. We
only need to find a necessary and sufficient condition for
existing a lifting of $u_lu_{l+1}$ with an initial vertex
$(u_l)_{a\circ_{i_1}\tau
(u_1u_2)\circ_{i_2}\cdots\circ_{i_{l-1}}}\tau (u_{l-1}u_l)$. By
definition, there exists such a lifting of the edge $u_lu_{l+1}$
if and only if $(a\circ_{i_1}\tau
(u_1u_2)\circ_{i_2}\cdots\circ_{i_{l-1}})\tau
(u_{l-1}u_l))\circ_l\tau (u_lu_{l+1})\in\varsigma^{-1}(u_{l+1})$.
Since $\tau (u_lu_{l+1})\in\varsigma^{-1}(u_{l+1})$ by the
definition of multi-voltage graphs of type $2$, we know that
$a\circ_{i_1}\tau (u_1u_2)\circ_{i_2}\cdots\circ_{i_{l-1}}\tau
(u_{l-1}u_l)\in\varsigma^{-1}(u_{l+1})$.

Continuing this process, we get the assertion of this theorem by
the induction principle. \quad\quad $\natural$

\vskip 4mm

\no{\bf Corollary $2.2.7$} \ {\it Let $G$ a graph with vertices
partition $V(G)=\bigcup\limits_{i=1}^nV_i$ and let $(\Gamma ;
\circ)$ be a finite group, $\Gamma_i\prec\Gamma$ for any integer
$i, 1\leq i\leq n$. If $(G,\tau ,\varsigma )$ is a multi-voltage
graph with $\tau: X_{\frac{1}{2}}(G)\rightarrow\Gamma$ and
$\varsigma :V_i\rightarrow\Gamma_i$ for any integer $i, 1\leq
i\leq n$, then for a walk $W$ in $G$ with an initial vertex $u$,
there exists a lifting $W^{(\tau ,\varsigma )}$ in $G^{(\tau
,\varsigma )}$ with the initial vertex $u_a,
a\in\varsigma^{-1}(u)$ if and only if $a\in\bigcap_{v\in
V(W)}\varsigma^{-1}(v)$.}

\vskip 3mm

Similar to multi-voltage graphs of type $1$, we can get the exact
number of liftings of a walk in the case of $\Gamma_i=\Gamma$ and
$V_i=V(G)$ for any integer $i, 1\leq i\leq n$.

\vskip 4mm

\no{\bf Theorem $2.2.8$} \ {\it Let $\widetilde{\Gamma} =
\bigcup\limits_{i=1}^n\Gamma$ be a finite multi-group with groups
$(\Gamma ;\circ_i), 1\leq i\leq n$ and let $W = e^1e^2\cdots e^k$
be a walk with an initial vertex $u$ in a multi-voltage graph
$(G,\tau ,\varsigma )$ , $\tau:
X_{\frac{1}{2}}(G)\rightarrow\bigcap\limits_{i=1}^n\Gamma$ and
$\varsigma :V(G)\rightarrow\Gamma$, of type $2$. Then there are
$n^k$ liftings of $W$ in $G^{(\tau ,\varsigma )}$ with an initial
vertex $u_a$ for $\forall a\in\widetilde{\Gamma}$.}

\vskip 3mm

{\it Proof} \ The proof is similar to the proof of Theorem
$2.2.2$.
 \quad\quad $\natural$

\vskip 4mm

\no{\bf Theorem $2.2.9$} \ {\it Let $\widetilde{\Gamma} =
\bigcup\limits_{i=1}^n\Gamma$ be a finite multi-group with groups
$(\Gamma ;\circ_i), 1\leq i\leq n$, $C = u_1u_2\cdots u_mu_1$ a
circuit in a multi-voltage graph $(G,\tau ,\varsigma )$ of type
$2$ where $\tau:
X_{\frac{1}{2}}(G)\rightarrow\bigcap\limits_{i=1}^n\Gamma$ and
$\varsigma :V(G)\rightarrow\Gamma$. Then  there are $\frac{|\Gamma
|}{o(\tau (C, \circ_i))}$ liftings of length $o(tau (C,
\circ_i))m$ in $G^{(\tau ,\varsigma )}$ of $C$ for any integer $i,
1\leq i\leq n$, where $\tau (C, \circ_i)= \tau
(u_1,u_2)\circ_i\tau (u_2,u_3)\circ_i\cdots\circ_i\tau
(u_{m-1},u_m)\circ_i\tau (u_m,u_1)$ and there are }

$$\sum\limits_{i=1}^n\frac{|\Gamma |}{o(\tau (C, \circ_i))}$$

\no{\it liftings of $C$ in $G^{(\tau ,\varsigma )}$ altogether.}

\vskip 3mm

{\it Proof} \ The proof is similar to the proof of Theorem
$2.2.3$. \quad\quad $\natural$

\vskip 4mm

\no{\bf Definition $2.2.5$} \ {\it Let $G_1, G_2$ be two graphs
and $H$ a subgraph of $G_1$ and $G_2$. A one-to-one mapping $\xi$
between $G_1$ and $G_2$ is called an $H$-isomorphism if for any
subgraph $J$ isomorphic to $H$ in $G_1$, $\xi (J)$ is also a
subgraph isomorphic to $H$ in $G_2$.

If $G_1=G_2=G$, then an $H$-isomorphism between $G_1$ and $G_2$ is
called an $H$-automorphism of $G$. Certainly, all
$H$-automorphisms form a group under the composition operation,
denoted by ${\rm Aut}_HG$ and ${\rm Aut}_HG = {\rm Aut}G$ if we
take $H=K_2$.}

\vskip 3mm

For example, let $H=\left<E(x, N_{G}(x))\right>$ for $\forall x\in
V(G)$. Then the $H$-automorphism group of a complete bipartite
graph $K(n,m)$ is ${\rm Aut}_HK(n,m)=S_n[S_m]={\rm Aut}K(n,m)$.
There $H$-automorphisms are called {\it star-automorphisms}.

\vskip 4mm

\no{\bf Theorem $2.2.10$} \ {\it Let $G$ be a graph. If there is a
decomposition $G = \bigoplus\limits_{i=1}^nH_i$ with $H_i\cong H$
for $1\leq i\leq n$ and $H = \bigoplus\limits_{j=1}^mJ_j$ with
$J_j\cong J$ for $1\leq j\leq m$, then

($i$) \ $\left<\iota_i, \iota_i:H_1\rightarrow H_i, \ an  \
isomorphism, \ 1\leq i\leq n\right> = S_n\preceq {\rm Aut}_HG$,
and particularly, $S_n\preceq {\rm Aut}_HK_{2n+1}$ if $H = C$, a
hamiltonian circuit in $K_{2n+1}$.

($ii$) \ ${\rm Aut}_JG\preceq {\rm Aut}_HG$, and particularly,
${\rm Aut}G \preceq {\rm Aut}_HG$ for a simple graph $G$.}

\vskip 3mm

{\it Proof} \ ($i$) \ For any integer $i, 1\leq i\leq n$, we prove
there is a such $H$-automorphism $\iota$ on $G$ that $\iota_i:
H_1\rightarrow H_i$. In fact, since $H_i\cong H$, $1\leq i\leq n$,
there is an isomorphism $\theta :H_1\rightarrow H_i$. We define
$\iota_i$ as follows:

\vskip 3mm

\[
\iota_i(e)=\left\{\begin{array}{ll} \theta
(e), & \ {\rm if} \ e\in V(H_1)\bigcup E(H_1),\\
e, & \ {\rm if} e\in (V(G)\setminus V(H_1))\bigcup (E(G)\setminus
E(H_1)).\\
\end{array}
\right.
\]

\vskip 2mm

\no Then $\iota_i$ is a one-to-one mapping on the graph $G$ and is
also an $H$-isomorphism by definition. Whence,

$$\left<\iota_i, \iota_i:H_1\rightarrow H_i, \ {\rm an  \
isomorphism}, \ 1\leq i\leq n\right>\preceq {\rm Aut}_HG.$$

Since $\left<\iota_i, 1\leq i\leq n\right>\cong\left<(1,i), 1\leq
i\leq n\right> = S_n$, thereby we get that $S_n\preceq {\rm
Aut}_HG.$

For a complete graph $K_{2n+1}$, we know a decomposition $K_{2n+1}
= \bigoplus\limits_{i=1}^nC_i$ with

$$C_i= v_0v_iv_{i+1}v_{i-1}v_{i-2}\cdots v_{n+i-1}v_{n+i+1}v_{n+i}v_0$$

\no for any integer $i, 1\leq i\leq n$ by Theorem $2.1.18$.
Therefore, we get that

$$S_n\preceq {\rm Aut}_HK_{2n+1}$$

\no if we choose a hamiltonian circuit $H$ in $K_{2n+1}$.

($ii$) \ Choose $\sigma\in {\rm Aut}_JG$. By definition, for any
subgraph $A$ of $G$, if $A\cong J$, then $\sigma (A)\cong J$.
Notice that $H = \bigoplus\limits_{j=1}^mJ_j$ with $J_j\cong J$
for $1\leq j\leq m$. Therefore, for any subgraph $B, B\cong H$ of
$G$, $\sigma (B)\cong\bigoplus\limits_{j=1}^m \sigma (J_j)\cong
H$. This fact implies that $\sigma\in {\rm Aut}_HG$.

Notice that for a simple graph $G$, we have a decomposition $G =
\bigoplus\limits_{i=1}^{\varepsilon (G)}K_2$ and ${\rm Aut}_{K_2}G
= {\rm Aut}G$. Whence, ${\rm Aut}G \preceq {\rm Aut}_HG$.
\quad\quad $\natural$\vskip 2mm

The equality in Theorem $2.2.10 (ii)$ does not always hold. For
example, a one-to-one mapping $\sigma$ on the lifting graph of
Fig.$2.21(a)$: $\sigma (u_0)=u_1$, $\sigma (u_1)=u_0$, $\sigma
(v_0)=v_1$, $\sigma (v_1)=v_2$ and $\sigma (v_2)=v_0$ is not an
automorphism, but it is an $H$-automorphism with $H$ being a star
$S_{1.2}$.

For automorphisms of the lifting $G^{(\tau ,\varsigma )}$ of a
multi-voltage graph $(G,\tau ,\varsigma )$ of type $2$, we get a
result in the following.

\vskip 4mm

\no{\bf Theorem $2.2.11$} \ {\it Let $(G,\tau ,\varsigma )$ be a
multi-voltage graph of type $2$ with $\tau:
X_{\frac{1}{2}}(G)\rightarrow\bigcap\limits_{i=1}^n\Gamma_i$ and
$\varsigma :V_i\rightarrow\Gamma_i$. Then for any integers $i,j,
1\leq i,j\leq n$, }

($i$) \ {\it for $\forall g_i\in\Gamma_i$, the left action
$lA(g_i)$ on $\left<V_i\right>^{(\tau ,\varsigma )}$ is a
fixed-free action of an automorphism of $\left<V_i\right>^{(\tau
,\varsigma )}$;}

($ii$) \ {\it for $\forall g_{ij}\in\Gamma_i\bigcap\Gamma_j$, the
left action $lA(g_{ij})$ on $\left<E_G(V_i,V_j)\right>^{(\tau
,\varsigma )}$ is a star-automorphism of
$\left<E_G(V_i,V_j)\right>^{(\tau ,\varsigma )}$.}

\vskip 3mm

{\it Proof} \ The proof of ($i$) is similar to the proof of
Theorem $2.2.4$. We prove the assertion ($ii$). A star with a
central vertex $u_a$, $u\in V_i, a\in\Gamma_i\bigcap\Gamma_j$ is
the graph $S_{star} =\left<\{(u_a, v_{a\circ_jb}) \ {\rm if} \
(u,v)\in E_G(V_i,V_j), \tau (u,v)=b \}\right>$. By definition, the
left action $lA(g_{ij})$ is a one-to-one mapping on
$\left<E_G(V_i,V_j)\right>^{(\tau ,\varsigma )}$. Now for any
element $g_{ij}, g_{ij}\in\Gamma_i\bigcap\Gamma_j$, the left
action $lA(g_{ij})$ of $g_{ij}$ on a star $S_{star}$ is

$$lA(g_{ij})(S_{star})= \left<\{(u_{g_{ij}\circ_ia},
v_{(g_{ij}\circ_ia)\circ_jb}) \ {\rm if} \ (u,v)\in E_G(V_i,V_j),
\tau (u,v)=b \}\right> = S_{star}.$$

\no Whence, $lA(g_{ij})$ is a star-automorphism of
$\left<E_G(V_i,V_j)\right>^{(\tau ,\varsigma )}$. \quad\quad
$\natural$

Let $\widetilde{G}$ be a graph and let
$\widetilde{\Gamma}=\bigcup\limits_{i=1}^n\Gamma_i$ be a finite
multi-group. If there is a partition for the vertex set
$V(\widetilde{G})=\bigcup\limits_{i=1}^nV_i$ such that the action
of $\widetilde{\Gamma}$ on $\widetilde{G}$ consists of $\Gamma_i$
action on $\left<V_i\right>$ and $\Gamma_i\bigcap\Gamma_j$ on
$\left<E_G(V_i,v_j)\right>$ for $1\leq i, j\leq n$, then we say
this action to be a {\it partially-action}. A partially-action is
called {\it fixed-free} if $\Gamma_i$ is fixed-free on
$\left<V_i\right>$ and the action of each element in
$\Gamma_i\bigcap\Gamma_j$ is a star-automorphism and fixed-free on
$\left<E_G(V_i,V_j)\right>$ for any integers $i, j, 1\leq i, j\leq
n$. These orbits of a partially-action are defined to be

$$orb_i(v) = \{g(v) | g\in\Gamma_i, v\in V_i \}$$

\no for any integer $i, 1\leq i\leq n$ and

$$orb(e) = \{g(e) | e\in E(\widetilde{G}),g\in\widetilde{\Gamma} \}.$$

A {\it partially-quotient graph}
$\widetilde{G}/_p\widetilde{\Gamma}$ is defined by

$$V(\widetilde{G}/_p\widetilde{\Gamma})=\bigcup\limits_{i=1}^n\{ \ orb_i(v) \ | \ v\in
V_i\}, \ \ E(\widetilde{G}/_p\widetilde{\Gamma})=\{orb(e) | e\in
E(\widetilde{G})\}$$

\no and
$I(\widetilde{G}/_p\widetilde{\Gamma})=\{I(e)=(orb_i(u),orb_j(v))
\ {\rm if} \ u\in V_i, v\in V_j \ {\rm and} \ (u,v)\in
E(\widetilde{G}), 1\leq i,j\leq n\}$. An example for
partially-quotient graph is shown in Fig.$2.22$, where
$V_1=\{u_0,u_1,u_2,u_3\}$, $V_2=\{v_0,v_1,v_2\}$ and $\Gamma_1 =
Z_4$, $\Gamma_2 = Z_3$.

\includegraphics[bb=5 10 350 110]{sgm24.eps}

\vskip 3mm

\c{\bf Fig $2.22$}\vskip 4mm

Then we have a necessary and sufficient condition for the lifting
of a multi-voltage graph of type $2$.

\vskip 4mm

\no{\bf Theorem $2.2.12$} \ {\it If the partially-action
${\mathcal P_a}$ of a finite multi-group
$\widetilde{\Gamma}=\bigcup\limits_{i=1}^n\Gamma_i$ on a graph
$\widetilde{G}$, $V(\widetilde{G})=\bigcup\limits_{i=1}^nV_i$ is
fixed-free, then there is a multi-voltage graph
$(\widetilde{G}/_p\widetilde{\Gamma},\tau ,\varsigma )$, $\tau:
X_{\frac{1}{2}}(\widetilde{G}/\widetilde{\Gamma})
\rightarrow\widetilde{\Gamma}$, $\varsigma:
V_i\rightarrow\Gamma_i$ of type $2$ such that}

$$\widetilde{G} \ \cong \ (\widetilde{G}/_p\widetilde{\Gamma})^{(\tau ,\varsigma )}.$$

\vskip 3mm

{\it Proof} \ Similar to the proof of Theorem $2.2.6$, we also
choose positive directions on these edges of
$\widetilde{G}/_p\widetilde{\Gamma}$ and $\widetilde{G}$ so that
the partially-quotient map
$p_{\widetilde{\Gamma}}:\widetilde{G}\rightarrow\widetilde{G}/_p\widetilde{\Gamma}$
is direction-preserving and the partially-action of
$\widetilde{\Gamma}$ on $\widetilde{G}$ preserves directions.

For any integer $i, 1\leq i\leq n$ and $\forall v^i\in V_i$, we
can label $v^i$ as $v^i_{1_{\Gamma_i}}$ and for every group
element $g_i\in\Gamma_i, g_i\not=1_{\Gamma_i}$, label the vertex
${\mathcal P_a}(g_i)((v_i)_{1_{\Gamma_i}})$ as $v^i_{g_i}$. Now if
the edge $e$ of $\widetilde{G}/_p\widetilde{\Gamma}$ runs from $u$
to $w$, we assign the label $e_{g_i}$ to the edge of the orbit
$p^{-1}(e)$ that originates at the vertex $u^i_{g_i}$ and
terminates at $w^j_{h_j}$.

Since $\Gamma_i$ acts freely on $\left<V_i\right>$, there are just
$|\Gamma_i|$ edges in the orbit $p_{\Gamma_i}^{-1}(e)$ for each
integer $i, 1\leq i\leq n$, one originating at each of the
vertices in the vertex orbit $p_{\Gamma_i}^{-1}(v)$. Thus for any
integer $i, 1\leq i\leq n$, the choice of an edge in $p^{-1}(e)$
to be labelled $e_{g_i}$ is unique. Finally, if the terminal
vertex of the edge $e_{g_i}$ is $w^j_{h_j}$, one assigns voltage
$g_i^{-1}\circ_jh_j$ to the edge $e$ in the partially-quotient
graph $\widetilde{G}/_p\widetilde{\Gamma}$ if $g_i,
h_j\in\Gamma_i\bigcap\Gamma_j$ for $1\leq i,j\leq n$.

Under this labelling process, the isomorphism $\vartheta
:\widetilde{G}\rightarrow
(\widetilde{G}/_p\widetilde{\Gamma})^{(\tau ,\varsigma )}$
identifies orbits in $\widetilde{G}$ with fibers of $G^{(\tau
,\varsigma )}$. \quad\quad $\natural$

\vskip 2mm

The multi-voltage graphs defined in this section enables us to
enlarge the application field of voltage graphs. For example, a
complete bipartite graph $K(n,m)$ is a lifting of a multi-voltage
graph, but it is not a lifting of a voltage graph in general if
$n\not=m$.

\vskip 5mm

\no{\bf \S $2.3$ \ Graphs in a Space}

\vskip 4mm

\no For two topological spaces ${\mathcal E}_1$ and ${\mathcal
E}_2$, an {embedding} of ${\mathcal E}_1$ in ${\mathcal E}_2$ is a
one-to-one continuous mapping $f:{\mathcal E}_1\rightarrow
{\mathcal E}_2$ (see $[92]$ for details). Certainly, the same
problem can be also considered for ${\mathcal E}_2$ being a metric
space. By a topological view, a graph is nothing but a
$1$-complex, we consider the embedding problem for graphs in
spaces or on surfaces in this section. The same problem had been
considered by Gr\"{u}mbaum in [$25$]-[$26$] for graphs in spaces
and in these references $[6],[23]$,$[42]-[44]$,$[56],[69]$ and
$[106]$ for graphs on surfaces.

\vskip 4mm

\no{\bf $2.3.1.$ Graphs in an $n$-manifold}

\vskip 3mm

\no For a positive integer $n$, an $n$-manifold ${\bf M}^n$ is a
Hausdorff space such that each point has an open neighborhood
homeomorphic to an open $n$-dimensional ball
$B^n=\{(x_1,x_2,\cdots,x_n)|x_1^2+x_2^2+\cdots+x_n^2< 1\}$. For a
given graph $G$ and an $n$-manifold ${\bf M}^n$ with $n\geq 3$,
the embeddability of $G$ in ${\bf M}^n$ is trivial. We
characterize an embedding of a graph in an $n$-dimensional
manifold ${\bf M}^n$ for $n\geq 3$ similar to the rotation
embedding scheme of a graph on a surface (see
$[23],[42]-[44],[69]$ for details) in this section.

For $\forall v\in V(G)$, a {\it space permutation $P(v)$} of $v$
is a permutation on $N_G(v)=\{u_1,u_2,\cdots,u_{\rho_G(v)}\}$ and
all space permutation of a vertex $v$ is denoted by ${\mathcal
P}_s(v)$. We define a {\it space permutation $P_s(G)$ of a graph
$G$} to be

$$P_s(G)=\{P(v)| \forall v\in V(G), P(v)\in{\mathcal
P}_s(v)\}$$

\no and a {\it permutation system ${\mathcal P}_s(G)$ of $G$} to
be all space permutation of $G$. Then we have the following
characteristic for an embedded graph in an $n$-manifold $\bf M^n$
with $n\geq 3$.

\vskip 4mm

\no{\bf Theorem $2.3.1$} \ {\it For an integer $n\geq 3$, every
space permutation $P_s(G)$ of a graph $G$ defines a unique
embedding of $G\rightarrow {\bf M}^n$. Conversely, every embedding
of a graph $G\rightarrow {\bf M}^n$ defines a space permutation of
$G$.}

\vskip 3mm

{\it Proof} \ Assume $G$ is embedded in an $n$-manifold ${\bf
M}^n$. For $\forall v\in V(G)$, define an $(n-1)$-ball
$B^{n-1}(v)$ to be $x_1^2+x_2^2+\cdots+x_n^2= r^2$ with center at
$v$ and radius $r$ as small as needed. Notice that all
autohomeomorphisms ${\rm Aut}B^{n-1}(v)$ of $B^{n-1}(v)$ is a
group under the composition operation and two points
$A=(x_1,x_2,\cdots,x_n)$ and $B=(y_1,y_2,\cdots,y_n)$ in
$B^{n-1}(v)$ are said to be combinatorially equivalent if there
exists an autohomeomorphism $\varsigma\in{\rm Aut}B^{n-1}(v)$ such
that $\varsigma(A)=B$. Consider intersection points of edges in
$E_G(v,N_G(v))$ with $B^{n-1}(v)$. We get a permutation $P(v)$ on
these points, or equivalently on $N_G(v)$ by $(A,B,\cdots,C,D)$
being a cycle of $P(v)$ if and only if there exists
$\varsigma\in{\rm Aut}B^{n-1}(v)$ such that $\varsigma^i(A) = B$,
$\cdots$, $\varsigma^j(C)=D$ and $\varsigma^l(D)=A$, where
$i,\cdots,j,l$ are integers. Thereby we get a space permutation
$P_s(G)$ of $G$.

Conversely, for a space permutation $P_s(G)$, we can embed $G$ in
${\bf M}^n$ by embedding each vertex $v\in V(G)$ to a point $X$ of
${\bf M}^n$ and arranging vertices in one cycle of $P_s(G)$ of
$N_G(v)$ as the same orbit of $\left<\sigma\right>$ action on
points of $N_G(v)$ for $\sigma\in{\rm Aut}B^{n-1}(X)$. Whence we
get an embedding of $G$ in the manifold ${\bf M}^n$.\quad\quad
$\natural$

Theorem $2.3.1$ establishes a relation for an embedded graph in an
$n$-dimensional manifold with a permutation, which enables us to
give a combinatorial definition for graphs embedded in
$n$-dimensional manifolds, see Definition $2.3.6$ in the finial
part of this section.

\vskip 4mm

\no{\bf Corollary $2.3.1$} \ {\it For a graph $G$, the number of
embeddings of G in ${\bf M}^n, n\geq 3$ is}

$$\prod\limits_{v\in V(G)}\rho_G(v)!.$$

\vskip 3mm

For applying graphs in spaces to theoretical physics, we consider
an embedding of a graph in an manifold with some additional
conditions which enables us to find good behavior of a graph in
spaces. On the first, we consider rectilinear embeddings of a
graph in an Euclid space.

\vskip 4mm

\no{\bf Definition $2.3.1$} \ {\it For a given graph $G$ and an
Euclid space ${\bf E}$, a rectilinear embedding of $G$ in ${\bf
E}$ is a one-to-one continuous mapping $\pi : G\rightarrow {\bf
E}$ such that}

($i$) \ {\it for $\forall e\in E(G)$, $\pi (e)$ is a segment of a
straight line in ${\bf E}$;}

($ii$) \ {\it for any two edges $e_1=(u,v),e_2=(x,y)$ in $E(G)$,
$(\pi (e_1)\setminus\{\pi (u),\pi (v)\})\bigcap$ $(\pi
(e_2)\setminus\{\pi (x),\pi (y)\})=\emptyset$.}

\vskip 3mm

In ${\bf R}^3$, a rectilinear embedding of $K_4$ and a cube $Q_3$
are shown in Fig.$2.23$.

\includegraphics[bb=5 10 350 130]{sgm25.eps}

\vskip 3mm

\c{\bf Fig $2.23$}\vskip 3mm

In general, we know the following result for rectilinear embedding
of $G$ in an Euclid space ${\bf R}^n, n\geq 3$.

 \vskip 4mm

\no{\bf Theorem $2.3.2$} \ {\it For any simple graph $G$ of order
$n$, there is a rectilinear embedding of $G$ in ${\bf R}^n$ with
$n\geq 3$.}

\vskip 3mm

{\it Proof} \ We only need to prove this assertion for $n=3$. In
${\bf R}^3$, choose $n$ points $(t_1,t_1^2,t_1^3),
(t_2,t_2^2,t_2^3),\cdots , (t_n,t_n^2,t_n^3)$, where
$t_1,t_2,\cdots ,t_n$ are $n$ different real numbers. For integers
$i,j,k,l, 1\leq i,j,k,l\leq n$, if a straight line passing through
vertices $(t_i,t_i^2,t_i^3)$ and $(t_j,t_j^2,t_j^3)$ intersects
with a straight line passing through vertices $(t_k,t_k^2,t_k^3)$
and $(t_l,t_l^2,t_l^3)$, then there must be

\[
\left|
\begin{array}{ccc}
t_k-t_i & t_j-t_i & t_l-t_k\\
t_k^2-t_i^2 & t_j^2-t_i^2 & t_l^2-t_k^2\\
t_k^3-t_i^3 & t_j^3-t_i^3 & t_l^3-t_k^3\\
\end{array}
\right|=0,
\]

\no which implies that there exist integers $s,f\in\{k,l,i,j\}$,
$s\not= f$ such that $t_s= t_f$, a contradiction.

Let $V(G)=\{v_1,v_2,\cdots ,v_n\}$. We embed the graph $G$ in
${\bf R}^3$ by a mapping $\pi : G\rightarrow {\bf R}^3$ with $\pi
(v_i)= (t_i,t_i^2,t_i^3)$ for $1\leq i\leq n$ and if $v_iv_j\in
E(G)$, define $\pi (v_iv_j)$ being the segment between points
$(t_i,t_i^2,t_i^3)$ and $(t_j,t_j^2,t_j^3)$ of a straight line
passing through points $(t_i,t_i^2,t_i^3)$ and
$(t_j,t_j^2,t_j^3)$. Then $\pi$ is a rectilinear embedding of the
graph $G$ in ${\bf R}^3$. \quad\quad $\natural$

For a graph $G$ and a surface $S$, an {\it immersion} $\iota $ of
$G$ on $S$ is a one-to-one continuous mapping $\iota :G\rightarrow
S$ such that for $\forall e\in E(G)$, if $e=(u,v)$, then $\iota
(e)$ is a curve connecting $\iota (u)$ and $\iota (v)$ on $S$. The
following two definitions are generalization of embedding of a
graph on a surface.

\vskip 4mm

\no{\bf Definition $2.3.2$} \ {\it Let $G$ be a graph and $S$ a
surface in a metric space ${\mathcal E}$. A pseudo-embedding of
$G$ on $S$ is a one-to-one continuous mapping $\pi : G\rightarrow
{\mathcal E}$ such that there exists vertices $V_1\subset V(G)$,
$\pi | _{\left<V_1\right>}$ is an immersion on $S$ with each
component of $S\setminus\pi (\left<V_1\right>)$ isomorphic to an
open $2$-disk.}

\vskip 4mm

\no{\bf Definition $2.3.3$} \ {\it Let $G$ be a graph with a
vertex set partition $V(G)=\bigcup\limits_{j=1}^kV_i$, $V_i\bigcap
V_j=\emptyset$ for $1\leq i,j\leq k$ and let $S_1,S_2,\cdots ,S_k$
be surfaces in a metric space ${\mathcal E}$ with $k\geq 1$. A
multi-embedding of $G$ on $S_1,S_2,\cdots ,S_k$ is a one-to-one
continuous mapping $\pi : G\rightarrow {\mathcal E}$ such that for
any integer $i, 1\leq i\leq k$, $\pi |_{\left<V_i\right>}$ is an
immersion with each component of $S_i\setminus\pi
(\left<V_i\right>)$ isomorphic to an open $2$-disk.}

\vskip 3mm

Notice that if $\pi (G)\bigcap (S_1\bigcup S_2\cdots\bigcup
S_k)=\pi(V(G))$, then every $\pi : G\rightarrow {\bf R}^3$ is a
multi-embedding of $G$. We say it to be a {\it trivial
multi-embedding} of $G$ on $S_1,S_2,\cdots ,S_k$. If $k=1$, then
every trivial multi-embedding is a trivial pseudo-embedding  of
$G$ on $S_1$. The main object of this section is to find
nontrivial multi-embedding of $G$ on $S_1,S_2,\cdots ,S_k$ with
$k\geq 1$. The existence pseudo-embedding of a graph $G$ is
obvious by definition. We concentrate our attention on
characteristics of multi-embeddings of a graph.

For a graph $G$, let $G_1,G_2,\cdots ,G_k$ be $k$ vertex-induced
subgraphs of $G$. If $V(G_i)\bigcap V(G_j)=\emptyset$ for any
integers $i,j, 1\leq i,j\leq k$, it is called a {\it block
decomposition} of $G$ and denoted by

$$G = \biguplus\limits_{i=1}^k G_i.$$

\no The {\it planar block number $n_p(G)$} of $G$ is defined by

$$n_p(G) = min\{k | G = \biguplus\limits_{i=1}^k G_i, {\rm For \ any \ integer \ } i,
1\leq i\leq k, G_i \ {\rm is \ planar}\}.$$

Then we get a result for the planar black number of a graph $G$ in
the following.

\vskip 4mm

\no{\bf Theorem $2.3.3$} \ {\it A graph $G$ has a nontrivial
multi-embedding on $s$ spheres $P_1,P_2,\cdots ,$ $P_s$ with empty
overlapping if and only if $n_p(G)\leq s\leq |G|$.}

\vskip 3mm

{\it Proof} \ Assume $G$ has a nontrivial multi-embedding on
spheres $P_1,P_2,\cdots , P_s$. Since $|V(G)\bigcap P_i|\geq 1$
for any integer $i, 1\leq i\leq s$, we know that

$$|G|=\sum\limits_{i=1}^s|V(G)\bigcap P_i|\geq s.$$

By definition, if $\pi : G\rightarrow {\bf R}^3$ is a nontrivial
multi-embedding of $G$ on $P_1,P_2,\cdots , P_s$, then for any
integer $i, 1\leq i\leq s$, $\pi^{-1}(P_i)$ is a planar induced
graph. Therefore,

$$G = \biguplus\limits_{i=1}^s\pi^{-1}(P_i),$$

\no and we get that $s\geq n_p(G)$.

Now if $n_p(G)\leq s\leq |G|$, there is a block decomposition $G =
\biguplus\limits_{i=1}^s G_s$ of $G$ such that $G_i$ is a planar
graph for any integer $i, 1\leq i\leq s$. Whence we can take $s$
spheres $P_1,P_2,\cdots , P_s$ and define an embedding $\pi_i:
G_i\rightarrow P_i$ of $G_i$ on sphere $P_i$ for any integer $i,
1\leq i\leq s$.

Now define an immersion $\pi : G\rightarrow {\bf R}^3$ of $G$ on
${\bf R}^3$ by

$$\pi (G) = (\bigcup\limits_{i=1}^s\pi (G_i))\bigcup\{(v_i,v_j)|v_i\in V(G_i),v_j\in V(G_j),
(v_i,v_j)\in E(G), 1\leq i,j\leq s\}.$$

\no Then $\pi : G\rightarrow {\bf R}^3$ is a multi-embedding of
$G$ on spheres $P_1,P_2,\cdots , P_s$.\quad\quad $\natural$

For example, a multi-embedding of $K_6$ on two spheres is shown in
Fig.$2.24$, in where, $\left<\{x,y,z\}\right>$ is on one sphere
and $\left<\{u,v,w\}\right>$ on another.

\includegraphics[bb=5 10 350 120]{sgm26.eps}

\vskip 2mm

\c{\bf Fig $2.24$}\vskip 3mm

For a complete or a complete bipartite graph, we get the number
$n_p(G)$ as follows.

\vskip 4mm

\no{\bf Theorem $2.3.4$} \ {\it For any integers $n,m, n,m\geq 1$,
the numbers $n_p(K_n)$ and $n_p(K(m,n))$ are}

$$n_p(K_n) = \lceil\frac{n}{4}\rceil \ \ { and} \ \
n_p(K(m,n)) = 2 ,$$

\no{\it if $m\geq 3, n\geq 3$, otherwise $1$, respectively.}

\vskip 3mm

{\it Proof} \ Notice that every vertex-induced subgraph of a
complete graph $K_n$ is also a complete graph. By Theorem
$2.1.16$, we know that $K_5$ is non-planar. Thereby  we get that

$$n_p(K_n) = \lceil\frac{n}{4}\rceil$$

\no by definition of $n_p(K_n)$. Now for a complete bipartite
graph K(m,n), any vertex-induced subgraph by choosing $s$ and $l$
vertices from its two partite vertex sets is still a complete
bipartite graph. According to Theorem $2.1.16$, $K(3,3)$ is
non-planar and $K(2,k)$ is planar. If $m\leq 2$ or $n\leq 2$, we
get that $n_p(K(m,n))=1$. Otherwise, $K(m,n)$ is non-planar.
Thereby we know that $n_p(K(m,n))\geq 2$.

Let $V(K(m,n))= V_1\bigcup V_2$, where $V_1, V_2$ are its partite
vertex sets. If $m\geq 3$ and $n\geq 3$, we choose vertices
$u,v\in V_1$ and $x,y\in V_2$. Then the vertex-induced subgraphs
$\left<\{u,v\}\bigcup V_2\setminus\{x,y\}\right>$ and
$\left<\{x,y\}\bigcup V_2\setminus\{u,v\}\right>$ in $K(m,n)$ are
planar graphs. Whence, $n_p(K(m,n))=2$ by definition.\quad\quad
$\natural$

The position of surfaces $S_1,S_2,\cdots, S_k$ in a metric space
${\mathcal E}$ also influences the existence of multi-embeddings
of a graph. Among these cases an interesting case is there exists
an arrangement $S_{i_1}, S_{i_2},\cdots, S_{i_k}$ for
$S_1,S_2,\cdots, S_k$ such that in ${\mathcal E}$, $S_{i_j}$ is a
subspace of $S_{i_{j+1}}$ for any integer $j,1\leq j\leq k$. In
this case, the multi-embedding is called an {\it including
multi-embedding} of $G$ on surfaces $S_1,S_2,\cdots, S_k$.

\vskip 4mm

\no{\bf Theorem $2.3.5$} \ {\it A graph $G$ has a nontrivial
including multi-embedding on spheres $P_1\supset
P_2\supset\cdots\supset P_s$ if and only if there is a block
decomposition $G=\biguplus\limits_{i=1}^sG_i$ of $G$ such that for
any integer $i, 1 < i < s$,}

($i$) \ {\it $G_i$ is planar;}

($ii$) \ {\it for $\forall v\in V(G_i)$, $N_G(x)\subseteq
(\bigcup\limits_{j=i-1}^{i+1}V(G_j))$.}

\vskip 3mm

{\it Proof} \ Notice that in the case of spheres, if the radius of
a sphere is tending to infinite, an embedding of a graph on this
sphere is tending to a planar embedding. From this observation, we
get the necessity of these conditions.

Now if there is a block decomposition
$G=\biguplus\limits_{i=1}^sG_i$ of $G$ such that $G_i$ is planar
for any integer $i, 1 < i < s$ and $N_G(x)\subseteq
(\bigcup\limits_{j=i-1}^{i+1}V(G_j))$ for $\forall v\in V(G_i)$,
we can so place $s$ spheres $P_1,P_2,\cdots , P_s$ in ${\bf R}^3$
that $P_1\supset P_2\supset\cdots\supset P_s$. For any integer $i,
1\leq i\leq s$, we define an embedding $\pi_i: G_i\rightarrow P_i$
of $G_i$ on sphere $P_i$.

Since $N_G(x)\subseteq (\bigcup\limits_{j=i-1}^{i+1}V(G_j))$ for
$\forall v\in V(G_i)$, define an immersion $\pi : G\rightarrow
{\bf R}^3$ of $G$ on ${\bf R}^3$ by

$$\pi (G) = (\bigcup\limits_{i=1}^s\pi (G_i))\bigcup\{(v_i,v_j) | j=i-1,i,
i+1 \ for \ 1< i< s \ and \ (v_i,v_j)\in E(G)\}.$$

\no Then $\pi : G\rightarrow {\bf R}^3$ is a multi-embedding of
$G$ on spheres $P_1,P_2,\cdots , P_s$.\quad\quad $\natural$

\vskip 4mm

\no{\bf Corollary $2.3.2$} \ {\it If a graph $G$ has a nontrivial
including multi-embedding on spheres $P_1\supset
P_2\supset\cdots\supset P_s$, then the diameter $D(G)\geq s-1$.}

\vskip 4mm

\no{\bf $2.3.2.$ Graphs on a surface}

\vskip 3mm

\no In recent years, many books concern the embedding problem of
graphs on surfaces, such as Biggs and White's [$6$], Gross and
Tucker's [$23$], Mohar and Thomassen's [$69$] and White's [$106$]
on embeddings of graphs on surfaces and Liu's [$42$]-[$44$], Mao's
[$56$] and Tutte's [$100$] for combinatorial maps. Two disguises
of graphs on surfaces, i.e., {\it graph embedding} and {\it
combinatorial map} consist of two main streams in the development
of topological graph theory in the past decades. For relations of
these disguises with Klein surfaces, differential geometry and
Riemman geometry, one can see in Mao's [$55$]-[$56$] for details.

\vskip 4mm

\no{\bf $(1)$ The embedding of a graph}

\vskip 3mm

\no For a graph $G = (V(G),E(G),I(G))$ and a surface $S$, an
embedding of $G$ on $S$ is the case of $k=1$ in Definition
$2.3.3$, which is also an embedding of a graph in a $2$-manifold.
It can be shown immediately that if there exists an embedding of
$G$ on $S$, then $G$ is connected. Otherwise, we can get a
component in $S\setminus\pi (G)$ not isomorphic to an open
$2$-disk. Thereafter all graphs considered in this subsection are
connected.

Let $G$ be a graph. For $v\in V(G)$, denote all of edges incident
with the vertex $v$ by $N_{G}^e(v)=\{e_1,e_2,\cdots, e_{\rho_G
(v)}\}$. A permutation $C(v)$ on $e_1,e_2,\cdots, e_{\rho_G (v)}$
is said a {\it pure rotation} of $v$. All pure rotations incident
with a vertex $v$ is denoted by $\varrho (v)$. A {\it pure
rotation system} of $G$ is defined by

$$\rho (G) = \{C(v)| C(v)\in\varrho (v) \ {\rm for} \ \forall v\in V(G)\}$$

\no{and all pure rotation systems of $G$ is denoted by $\varrho
(G)$.}

Notice that in the case of embedded graphs on surfaces, a
$1$-dimensional ball is just a circle. By Theorem $2.3.1$, we get
a useful characteristic for embedding of graphs on orientable
surfaces first found by Heffter in 1891 and then formulated by
Edmonds in 1962. It can be restated as follows.

\vskip 4mm

\no{\bf Theorem $2.3.6$}([23]) {\it Every pure rotation system for
a graph $G$ induces a unique embedding of $G$ into an orientable
surface. Conversely, every embedding of a graph $G$ into an
orientable surface induces a unique pure rotation system of $G$. }

\vskip 3mm

According to this theorem, we know that the number of all
embeddings of a graph $G$ on orientable surfaces is $\prod_{v\in
V(G)}(\rho_G (v)-1)!$.

By a topological view, an embedded vertex or face can be viewed as
a disk, and an embedded edge can be viewed as a $1$-band which is
defined as a topological space $B$ together with a homeomorphism
$h:I\times I \rightarrow B$, where $I=[0,1]$, the unit interval.
Whence, an edge in an embedded graph has two sides. One side is
$h((0,x)), x\in I$. Another is $h((1,x)), x\in I$.

For an embedded graph $G$ on a surface, the two sides of an edge
$e\in E(G)$ may lie in two different faces $f_1$ and $f_2$, or in
one face $f$ without a twist ,or in one face $f$ with a twist such
as those cases (a), or (b), or (c) shown in Fig.$25$.

\includegraphics[bb=5 10 350 160]{sgm27.eps}

\vskip 3mm

\c{\bf Fig $2.25$}\vskip 3mm

Now we define a rotation system $\rho^L (G)$ to be a pair $({\cal
J},\lambda)$, where ${\cal J}$ is a pure rotation system of $G$,
and $\lambda : E(G) \rightarrow Z_2$. The edge with $\lambda (e)=
0$ or $\lambda (e)=1$ is called {\it type $0$} or {\it type $1$}
edge, respectively.  The {\it rotation system} $\varrho^L (G)$ of
a graph $G$ are defined by

$$\varrho^L (G) = \{({\cal
J},\lambda)| {\cal J}\in \varrho (G),\lambda : E(G) \rightarrow
Z_2 \}.$$

By Theorem $2.3.1$ we know the following characteristic for
embedding graphs on locally orientable surfaces.

\vskip 4mm

\no{\bf Theorem $2.3.7$}([23],[91]) {\it Every rotation system on
a graph $G$ defines a unique locally orientable embedding of
$G\rightarrow S$. Conversely, every embedding of a graph
$G\rightarrow S$ defines a rotation system for $G$. }

\vskip 3mm

Notice that in any embedding of a graph $G$, there exists a
spanning tree $T$ such that every edge on this tree is type $0$
(see also [$23$],[$91$] for details). Whence, the number of all
embeddings of a graph $G$ on locally orientable surfaces is

$$2^{\beta (G)}\prod_{v\in V(G)}(\rho_G
(v)-1)!$$

\no{and the number of all embedding of $G$ on non-orientable
surfaces is}

$$(2^{\beta (G)}-1)\prod_{v\in V(G)}(\rho
(v)-1)!.$$

The following result is the famous {\it Euler-Poincar\'{e}}
formula for embedding a graph on a surface.

\vskip 4mm

\no{\bf Theorem $2.3.8$} {\it If a graph $G$ can be embedded into
a surface $S$, then

$$\nu(G)-\varepsilon(G)+\phi (G) = \chi (S),$$

\no{where $\nu(G), \varepsilon(G)$ and $\phi (G) $ are the order,
size and the number of faces of $G$ on $S$, and $\chi (S)$ is the
Euler characteristic of $S$, i.e.,}}

\[
\chi (S)=\left\{\begin{array}{lr}
2-2p, & if \ S \ is \  orientable,\\
2-q, &  if \ S \ is \  non-orientable.
\end{array}
\right.
\]

For a given graph $G$ and a surface $S$, whether $G$ embeddable on
$S$ is uncertain. We use the notation $G\rightarrow S$ denoting
that $G$ can be embeddable on $S$. Define the {\it orientable
genus range} $GR^O(G)$ and the {\it non-orientable genus range}
$GR^N(G)$ of a graph $G$ by

$$GR^O(G) = \{\frac{2-\chi (S)}{2} | G\rightarrow S,  S \ is \ an \ orientable \ surface\},$$

$$GR^N(G) = \{2-\chi (S) | G\rightarrow S, S \ is \ a \ non-orientable \ surface\},$$

\no respectively and the orientable or non-orientable genus
$\gamma(G)$, $\overline{\gamma}(G)$ by

$$\gamma (G)= min\{p | p\in GR^O(G)\}, \ \ \gamma_M(G)= max\{p | p\in GR^O(G)\},$$

$$\widetilde{\gamma} (G)= min\{q | q\in GR^N(G)\}, \ \ \widetilde{\gamma}_M(G)= max\{q | q\in GR^O(G)\}.$$

\vskip 4mm

\no{\bf Theorem $2.3.9$}(Duke 1966) \ {\it Let $G$ be a connected
graph. Then}

$$GR^O(G) = [\gamma (G), \gamma_M(G)].$$

\vskip 3mm

{\it Proof} \ Notice that if we delete an edge $e$ and its
adjacent faces from an embedded graph $G$ on a surface $S$, we get
two holes at most, see Fig.$25$ also. This implies that $|\phi
(G)-\phi (G-e)|\leq 1$.

Now assume $G$ has been embedded on a surface of genus $\gamma
(G)$ and $V(G)=\{u,v,\cdots,w\}$. Consider those of edges adjacent
with $u$. Not loss of generality, we assume the rotation of $G$ at
vertex $v$ is $(e_1,e_2,\cdots, e_{\rho_G(u)})$. Construct an
embedded graph sequence $G_1,G_2,\cdots, G_{\rho_G(u)!}$ by

\vskip 2mm

$\varrho (G_1)=\varrho (G)$;

$\varrho (G_2)=(\varrho (G)\setminus\{\varrho
(u)\})\bigcup\{(e_2,e_1,e_3,\cdots,e_{\rho_G(u)})\}$;

$\cdots\cdots\cdots\cdots\cdots\cdots\cdots\cdots$;

$\varrho (G_{\rho_G(u)-1})=(\varrho (G)\setminus\{\varrho
(u)\})\bigcup\{(e_2,e_3,\cdots,e_{\rho_G(u)},e_1)\}$;

$\varrho (G_{\rho_G(u)})=(\varrho (G)\setminus\{\varrho
(u)\})\bigcup\{(e_3,e_2,\cdots,e_{\rho_G(u)},e_1)\}$;

$\cdots\cdots\cdots\cdots\cdots\cdots\cdots\cdots$;

$\varrho (G_{\rho_G(u)!})=(\varrho (G)\setminus\{\varrho
(u)\})\bigcup\{(e_{\rho_G(u)},\cdots,e_2,e_1,)\}$.

For any integer $i, 1\leq i\leq \rho_G(u)!$, since $|\phi (G)-\phi
(G-e)|\leq 1$ for $\forall e\in E(G)$, we know that $|\phi
(G_{i+1})-\phi (G_i)|\leq 1$. Whence, $|\chi (G_{i+1})-\chi
(G_i)|\leq 1$.

Continuing the above process for every vertex in $G$ we finally
get an embedding of $G$ with the maximum genus $\gamma_M(G)$.
Since in this sequence of embeddings of $G$, the genus of two
successive surfaces differs by at most one, we get that

$$GR^O(G) = [\gamma (G), \gamma_M(G)].\quad\quad \natural$$

The genus problem, i.e., {\it to determine the minimum orientable
or non-orientable genus of a graph} is NP-complete (see [$23$] for
details). Ringel and Youngs got the genus of $K_n$ completely by
{\it current graphs} (a dual form of voltage graphs) as follows.

\vskip 4mm

\no{\bf Theorem $2.3.10$} \ {\it For a complete graph $K_n$ and a
complete bipartite graph $K(m,n)$,  $m,n\geq 3$,}

$$\gamma (K_n) = \lceil\frac{(n-3)(n-4)}{12}\rceil \ and \
\gamma(K(m,n))=\lceil\frac{(m-2)(n-2)}{4}\rceil.$$

\vskip 3mm

Outline proofs for $\gamma(K_n)$ in Theorem $2.3.10$ can be found
in [$42$], [$23$],[$69$] and a complete proof is contained in
$[81]$. For a proof of $\gamma(K(m,n))$ in Theorem $2.3.10$ can be
also found in [$42$], [$23$],[$69$].

For the maximum genus $\gamma_M(G)$ of a graph, the time needed
for computation is bounded by a polynomial function on the number
of $\nu (G)$ ([$23$]). In 1979, Xuong got the following result.

\vskip 4mm

\no{\bf Theorem $2.3.11$}(Xuong 1979) \ {\it Let $G$ be a
connected graph with $n$ vertices and $q$ edges. Then}

$$\gamma_M(G) = \frac{1}{2}(q-n+1)-\frac{1}{2}\min\limits_{T} c_{odd}(G\setminus E(T)),$$

\no{\it where the minimum is taken over all spanning trees $T$ of
$G$ and $c_{odd}(G\setminus E(T))$ denotes the number of
components of $G\setminus E(T)$ with an odd number of edges.}

\vskip 3mm

In 1981, Nebesk\'{y} derived another important formula for the
maximum genus of a graph. For a connected graph $G$ and $A\subset
E(G)$, let $c(A)$ be the number of connected component of
$G\setminus A$ and let $b(A)$ be the number of connected
components $X$ of $G\setminus A$ such that $|E(X)|\equiv
|V(X)|(mod 2)$. With these notations, his formula can be restated
as in the next theorem.

\vskip 4mm

\no{\bf Theorem $2.3.12$}(Nebesk\'{y} 1981) \ {\it Let $G$ be a
connected graph with $n$ vertices and $q$ edges. Then}

$$\gamma_M(G)=\frac{1}{2}(q-n+2)-\max\limits_{A\subseteq E(G)}\{c(A)+b(A)-|A|\}.$$

\vskip 4mm

\no{\bf Corollary $2.3.3$} \ {\it The maximum genus of $K_n$ and
$K(m,n)$ are given by}

$$\gamma_M(K_n)= \lfloor\frac{(n-1)(n-2)}{4}\rfloor \ and \
\gamma_M(K(m,n))=\lfloor\frac{(m-1)(n-1)}{2}\rfloor,$$

\no{\it respectively.}

Now we turn to non-orientable embedding of a graph $G$. For
$\forall e\in E(G)$, we define an {\it edge-twisting surgery}
$\otimes (e)$ to be given the band of $e$ an extra twist such as
that shown in Fig.$26$.

\includegraphics[bb=5 10 350 120]{sgm28.eps}

\vskip 3mm

\c{\bf Fig $2.26$}

\vskip 3mm

Notice that for an embedded graph $G$ on a surface $S$, $e\in
E(G)$, if two sides of $e$ are in two different faces, then
$\otimes (e)$ will make these faces into one and if two sides of
$e$ are in one face, $\otimes (e)$ will divide the one face into
two. This property of $\otimes (e)$ enables us to get the
following result for the crosscap range of a graph.

\vskip 4mm

\no{\bf Theorem $2.3.13$}(Edmonds 1965, Stahl 1978) \ {\it Let $G$
be a connected graph. Then}

$$GR^N(G) = [\widetilde{\gamma} (G),\beta (G)],$$

\no{\it where $\beta (G)=\varepsilon (G)-\nu (G)+1$ is called the
Betti number of $G$.}

\vskip 3mm

{\it Proof} \  It can be checked immediately that
$\widetilde{\gamma}(G)=\widetilde{\gamma}_M(G)=0$ for a tree $G$.
If $G$ is not a tree, we have known there exists a spanning tree
$T$ such that every edge on this tree is type $0$ for any
embedding of $G$.

Let $E(G)\setminus E(T) = \{e_1,e_2,\cdots, e_{\beta (G)}\}$.
Adding the edge $e_1$ to $T$, we get a two faces embedding of
$T+e_1$. Now make edge-twisting surgery on $e_1$. Then we get a
one face embedding of $T+e_1$ on a surface. If we have get a one
face embedding of $T+(e_1+e_2+\cdots +e_i)$, $1\leq i < \beta
(G)$, adding the edge $e_{i+1}$ to $T+(e_1+e_2+\cdots +e_i)$ and
make $\otimes (e_{i+1})$ on the edge $e_{i+1}$. We also get a one
face embedding of $T+(e_1+e_2+\cdots +e_{i+1})$ on a surface
again.

Continuing this process until all edges in $E(G)\setminus E(T)$
have a twist, we finally get a one face embedding of
$T+(E(G)\setminus E(T))=G$ on a surface. Since the number of
twists in each circuit of this embedding of $G$ is $1(mod 2)$,
this embedding is non-orientable with only one face. By the
Euler-Poincar\'{e} formula, we know its genus $\widetilde{g}(G)$

$$\widetilde{g }(G)=2-(\nu (G)-\varepsilon (G)+1)=\beta (G).$$

For a minimum non-orientable embedding ${\mathcal E}_G$ of $G$,
i.e., $\widetilde{\gamma} ({\mathcal E}_G)= \widetilde{\gamma}
(G)$, one can selects an edge $e$ that lies in two faces of the
embedding ${\mathcal E}_G$ and makes $\otimes (e)$. Thus in at
most $\widetilde{\gamma}_M(G)-\widetilde{\gamma}(G)$ steps, one
has obtained all of embeddings of $G$ on every non-orientable
surface $N_s$ with $s\in
[\widetilde{\gamma}(G),\widetilde{\gamma}_M(G)]$. Therefore,

$$GR^N(G) = [\widetilde{\gamma} (G),\beta (G)]\quad\quad \natural$$

\vskip 4mm

\no{\bf Corollary $2.3.4$} \ {\it Let $G$ be a connected graph
with $p$ vertices and $q$ edges. Then}

$$\widetilde{\gamma}_M(G)= q-p+1.$$

\vskip 4mm

\no{\bf Theorem $2.3.14$} \ {For a complete graph $K_n$ and a
complete bipartite graph $K(m,n)$, $m,n\geq 3$,}

$$\widetilde{\gamma}(K_n) = \lceil\frac{(n-3)(n-4)}{6}\rceil$$

\no{\it with an exception value $\widetilde{\gamma}(K_7) = 3$ and}

$$\widetilde{\gamma}(K(m,n))=\lceil\frac{(m-2)(n-2)}{2}\rceil.$$

\vskip 3mm

A complete proof of this theorem is contained in $[81]$, Outline
proofs of Theorem $2.3.14$ can be found in [$42$].

\vskip 4mm

\no{\bf $(2)$ Combinatorial maps}

\vskip 3mm

\no Geometrically, an embedded graph of $G$ on a surface is called
a combinatorial map $M$ and say $G$ underlying $M$.  Tutte found
an algebraic representation for an embedded graph on a locally
orientable surface in 1973 $([98]$, which transfers a geometrical
partition of a surface to a permutation in algebra.

According to the summaries in Liu's $[43]-[44]$, a {\it
combinatorial map} $M = ({\mathcal X} _{\alpha,\beta},\mathcal P)$
is defined to be a permutation $\mathcal{P}$ acting on ${\mathcal
X} _{\alpha,\beta}$ of a disjoint union of quadricells $Kx$ of
$x\in X$, where $X$ is a finite set and
$K=\{1,\alpha,\beta,\alpha\beta \}$ is Klein $4$-group with the
following conditions hold.

($i$) \ $\forall x\in {\mathcal X} _{\alpha,\beta}$, there does
not exist an integer $k$ such that ${\mathcal P }^kx = \alpha x$;

($ii$) \ $\alpha {\mathcal P } = {\mathcal P}^{-1}\alpha$;

($iii$) The group $\Psi_{J}=\left<\alpha,\beta,{\mathcal
P}\right>$ is transitive on ${\mathcal X}_{\alpha,\beta}$.

The {\it vertices} of a combinatorial map are defined to be pairs
of conjugate orbits of ${\mathcal P}$ action on ${\mathcal
X}_{\alpha,\beta}$, {\it edges} to be orbits of $K$ on ${\mathcal
X}_{\alpha,\beta}$ and {\it faces} to be pairs of conjugate orbits
of ${\mathcal P}\alpha\beta$ action on ${\mathcal
X}_{\alpha,\beta}$.

For determining a map $({\mathcal X}_{\alpha,\beta},\mathcal{P})$
is orientable or not, the following condition is needed.\vskip 3mm

($iv$) {\it If the group $\Psi_I = \left<\alpha\beta ,{\mathcal
P}\right>$ is transitive on ${\mathcal X}_{\alpha,\beta}$, then
$M$ is non-orientable. Otherwise, orientable.}\vskip 2mm

For example, the graph $D_{0.4.0}$ (a dipole with $4$ multiple
edges ) on Klein bottle shown in Fig.$27$,

\includegraphics[bb=5 10 350 140]{sgm29.eps}

\vskip 3mm

\c{\bf Fig $2.27$}\vskip 3mm

\no can be algebraic represented by a combinatorial map $M =
({\mathcal X}_{\alpha,\beta},{\mathcal P})$ with

$${\mathcal X}_{\alpha,\beta} = \bigcup\limits_{e\in\{x,y,z,w\}}\{e,\alpha e,\beta
e,\alpha\beta e\},$$

\begin{eqnarray*}
{\mathcal P} &=& (x,y,z,w)(\alpha\beta x,\alpha\beta y,\beta
z,\beta w)\\
&\times& (\alpha x,\alpha w,\alpha z,\alpha y)(\beta x,\alpha\beta
w,\alpha\beta z,\beta y).
\end{eqnarray*}

\no This map has $2$ vertices $v_1= \{(x,y,z,w),(\alpha x,\alpha
w,\alpha z,\alpha y)\}$, $v_2= \{(\alpha\beta x,\alpha\beta
y,\beta z,$ $\beta w),(\beta x,\alpha\beta w,\alpha\beta z,\beta
y)\}$, $4$ edges $e_1= \{x,\alpha x,\beta x,\alpha\beta x\}$,
$e_2= \{y,\alpha y,\beta y,\alpha\beta y\}$, $e_3= \{z,\alpha
z,\beta z,\alpha\beta z\}$, $e_4= \{w,\alpha w,\beta w,\alpha\beta
w\}$ and $2$ faces $f_2= \{(x,\alpha\beta y,z,\beta y,\alpha
x,\alpha\beta w),$ $(\beta x,\alpha w,\alpha\beta x,y,\beta
z,\alpha y)\}$, $f_2= \{(\beta w,\alpha z),(w,\alpha\beta z)\}$.
The Euler characteristic of this map is

$$\chi (M) = 2- 4 + 2=0$$

\no and $\Psi_I = \left<\alpha\beta ,{\mathcal P}\right>$ is
transitive on ${\mathcal X}_{\alpha,\beta}$. Thereby it is a map
of $D_{0.4.0}$ on a Klein bottle with $2$ faces accordance with
its geometrical figure.

The following result was gotten by Tutte in $[98]$, which
establishes a relation for an embedded graph with a combinatorial
map.

\vskip 4mm

\no{\bf Theorem $2.3.15$} \ {\it For an embedded graph $G$ on a
locally orientable surface $S$, there exists one combinatorial map
$M = ({\mathcal X}_{\alpha,\beta},{\mathcal P})$ with an
underlying graph $G$ and for a combinatorial map $M = ({\mathcal
X}_{\alpha,\beta},{\mathcal P})$, there is an embedded graph $G$
underlying $M$ on $S$.}

\vskip 3mm

Similar to the definition of a multi-voltage graph (see $[56]$ for
details), we can define a multi-voltage map and its lifting by
applying a multi-group $\widetilde{\Gamma}=
\bigcup\limits_{i=1}^n\Gamma_i$ with $\Gamma_i=\Gamma_j$ for any
integers $i,j, 1\leq i,j\leq n$.

\vskip 4mm

\no{\bf Definition $2.3.4$} \ {\it Let $\widetilde{\Gamma} =
\bigcup\limits_{i=1}^n\Gamma$ be a finite multi-group with
$\Gamma=\{g_1,g_2,\cdots,g_m\}$ and an operation set
$O(\widetilde{\Gamma})=\{\circ_i | 1\leq i\leq n\}$ and let $M =
({\mathcal X}_{\alpha,\beta},{\mathcal P})$ be a combinatorial
map. If there is a mapping $\psi :{\mathcal
X}_{\alpha,\beta}\rightarrow\widetilde{\Gamma}$ such that}

($i$) \ {\it for $\forall x\in {\mathcal X}_{\alpha,\beta},
\forall \sigma\in K=\{1,\alpha,\beta,\alpha\beta\}$, $\psi(\alpha
x) =\psi(x)$, $\psi(\beta x)=\psi(\alpha\beta x)= \psi(x)^{-1}$;}

($ii$) \ {\it for any face $f=(x,y,\cdots,z)(\beta z,\cdots,\beta
y,\beta x)$, $\psi (f,i)= \psi (x)\circ_i\psi
(y)\circ_i\cdots\circ_i\psi (z)$, where $\circ_i\in
O(\widetilde{\Gamma})$, $1\leq i\leq n$ and $\left<\psi (f,i) |
f\in F(v)\right> = G$ for $\forall v\in V(G)$, where $F(v)$
denotes all faces incident with $v$,}

\no{\it then $(M,\psi)$ is called a multi-voltage map.}

\vskip 3mm

The {\it lifting of a multi-voltage map} is defined in the next
definition.

\vskip 4mm

\no{\bf Definition $2.3.5$} \ {\it For a multi-voltage map
$(M,\psi)$, the lifting map $M^{\psi}=({\mathcal
X}_{\alpha^{\psi},\beta^{\psi}}^{\psi},{\mathcal P}^{\psi})$ is
defined by}

$${\mathcal
X}_{\alpha^{\psi},\beta^{\psi}}^{\psi} = \{x_g | x\in {\mathcal
X}_{\alpha,\beta}, g\in\widetilde{\Gamma}\},$$

$${\mathcal P}^{\psi} = \prod\limits_{g\in\widetilde{\Gamma}}
\prod\limits_{(x,y,\cdots,z)(\alpha z,\cdots,\alpha y,\alpha x)\in
V(M)}(x_g,y_g,\cdots,z_g)(\alpha z_g,\cdots,\alpha y_g,\alpha
x_g),$$

$$\alpha^{\psi}= \prod\limits_{x\in {\mathcal X}_{\alpha,\beta},
g\in\widetilde{\Gamma}}(x_g,\alpha x_g),$$

$$\beta^{\psi} = \prod\limits_{i=1}^m\prod\limits_{x\in {\mathcal X}_{\alpha,\beta}}
(x_{g_i},(\beta x)_{g_i\circ_i\psi(x)})$$

\no{\it with a convention that $(\beta x)_{g_i\circ_i\psi(x)}=
y_{g_i}$ for some quadricells $y\in {\mathcal X}_{\alpha,\beta}$.}

Notice that the lifting $M^{\psi}$ is connected and $\Psi_I^{\psi}
= \left<\alpha^{\psi}\beta^{\psi} ,{\mathcal P}^{\psi}\right>$ is
transitive on ${\mathcal X}_{\alpha^{\psi},\beta^{\psi}}^{\psi}$
if and only if $\Psi_I = \left<\alpha\beta ,{\mathcal P}\right>$
is transitive on ${\mathcal X}_{\alpha,\beta}$. We get a result in
the following.

\vskip 4mm

\no{\bf Theorem $2.3.16$} \ {\it The Euler characteristic $\chi
(M^{\psi})$ of the lifting map $M^{\psi}$ of a multi-voltage map
$(M,\widetilde{\Gamma})$ is}

$$
\chi (M^{\psi})=|\Gamma|(\chi
(M)+\sum\limits_{i=1}^n\sum\limits_{f\in
F(M)}(\frac{1}{o(\psi(f,\circ_i))}-\frac{1}{n})),
$$

\no{\it where $F(M)$ and $o(\psi(f,\circ_i))$ denote the set of
faces in $M$ and the order of $\psi(f,\circ_i)$ in
$(\Gamma;\circ_i)$, respectively.}

\vskip 3mm

{\it Proof} \ By definition the lifting map $M^{\vartheta}$ has
$|\Gamma|\nu (M)$ vertices, $|\Gamma|\varepsilon (M)$ edges.
Notice that each lifting of the boundary walk of a face is a
homogenous lifting by definition of $\beta^{\psi}$. Similar to the
proof of Theorem $2.2.3$, we know that $M^{\vartheta}$ has
$\sum\limits_{i=1}^n\sum\limits_{f\in
F(M)}\frac{|\Gamma|}{o(\psi(f,\circ_i))}$ faces. By the
Eular-Poincar\'{e} formula we get that

\begin{eqnarray*}
\chi (M^{\psi}) &=& \nu (M^{\psi})-\varepsilon
(M^{\psi})+\phi (M^{\psi})\\
&=& |\Gamma|\nu (M)-|\Gamma|\varepsilon (M)+
\sum\limits_{i=1}^n\sum\limits_{f\in
F(M)}\frac{|\Gamma|}{o(\psi(f,\circ_i))}\\
&=& |\Gamma|(\chi (M)-\phi
(M)+\sum\limits_{i=1}^n\sum\limits_{f\in
F(M)}\frac{1}{o(\psi(f,\circ_i))}\\
&=& |G|(\chi (M)+\sum\limits_{i=1}^n\sum\limits_{f\in
F(M)}\frac{1}{o(\psi(f,\circ_i))}-\frac{1}{n}). \quad\quad
\natural
\end{eqnarray*}

Recently, more and more papers concentrated on finding {\it
regular maps} on surface, which are related with {\it discrete
groups, discrete geometry} and {\it crystal physics}. For this
object, an important way is by the voltage assignment on a map. In
this field, general results for automorphisms of the lifting map
are known, see $[45]-[46]$ and $[71]-[72]$ for details. It is also
an interesting problem for applying these multi-voltage maps to
finding non-regular or other maps with some constraint conditions.

Motivated by the Four Color Conjecture, Tait conjectured that {\it
every simple 3-polytope is hamiltonian} in 1880. By Steinitz's a
famous result (see [24]), this conjecture is equivalent to that
{\it every 3-connected cubic planar graph is hamiltonian}. Tutte
disproved this conjecture by giving a 3-connected non-hamiltonian
cubic planar graph with 46 vertices in 1946 and proved that {\it
every $4$-connected planar graph is hamiltonian} in
1956([$97$],[$99$]). In [$56$], Gr\"{u}nbaum conjectured that {\it
each $4$-connected graph embeddable in the torus or in the
projective plane is hamiltonian}. This conjecture had been solved
for the projective plane case by Thomas and Yu in $1994$ ([$93$]).
Notice that the splitting operator $\vartheta$ constructed in the
proof of Theorem $2.1.11$ is a planar operator. Applying Theorem
$2.1.11$ on surfaces we know that {\it for every map $M$ on a
surface, $M^{\vartheta}$ is non-hamiltonian}. In fact, we can
further get an interesting result related with Tait's conjecture.

\vskip 4mm

\no{\bf Theorem $2.3.17$} \ {\it There exist infinite
$3-$connected non-hamiltonian cubic maps on each locally
orientable surface.}

\vskip 3mm

 {\it Proof} \ Notice that there exist $3$-connected
 triangulations on every locally orientable surface $S$. Each dual
 of them is a $3$-connected cubic map on $S$. Now we define a
 splitting operator $\sigma$ as shown in Fig.$2.28$.

\includegraphics[bb=40 10 500 200]{sgm30.eps}

\vskip 3mm

\c{\bf Fig.$2.28$}\vskip 2mm

For a 3-connected cubic map $M$, we prove that $M^{\sigma(v)}$ is
non-hamiltonian for $\forall v\in V(M)$. According to Theorem
$2.1.7$, we only need to prove that there are no $y_1-y_2$, or
$y_1-y_3$, or $y_2-y_3$ hamiltonian path in the nucleus
$N(\sigma(v))$ of operator $\sigma$.

Let $H(z_i)$ be a component of $N(\sigma(v))\backslash
\{z_0z_i,y_{i-1}u_{i+1},y_{i+1}v_{i-1} \}$ which  contains the
vertex $z_i, 1\leq{i}\leq{3}$(all these indices mod 3). If there
exists a  $y_1-y_2$ hamiltonian path $P$ in $N(\sigma(v))$, we
prove that there must be a $u_i-v_i$ hamiltonian path in the
subgraph $H(z_i)$ for an integer $i, 1\leq{i}\leq{3}$.

Since $P$ is a hamiltonian path in $N(\sigma(v))$, there must be
that $v_1y_3u_2 \ {\rm or} \ u_2y_3v_1$ is a subpath of $P$. Now
let $E_1=\{y_1u_3,z_0z_3,y_2v_3 \}$, we know that $|E(P)\bigcap
E_1|=2$. Since $P$ is a $y_1-y_2$ hamiltonian path in the graph
$N(\sigma(v))$, we must have $y_1u_3 \not\in E(P)$ or $y_2v_3
\not\in E(P)$. Otherwise, by $|E(P)\bigcap S_1|=2$ we get that
$z_0z_3 \not\in E(P)$. But in this case, $P$ can not be a $y_1-
y_2$ hamiltonian path in $N(\sigma(v))$, a contradiction.

Assume $y_2v_3 \not\in E(P)$. Then $y_2u_1 \in E(P)$. Let $E_2=\{
u_1y_2,z_1z_0,v_1y_3 \}$. We also know that $|E(P)\bigcap E_2|=2$
by the assumption that $P$ is a hamiltonian path in
$N(\sigma(v))$. Hence $z_0z_1 \not\in E(P)$ and the $v_1-u_1$
subpath in $P$ is a $v_1-u_1$ hamiltonian path in the subgraph
$H(z_1)$.

Similarly, if $y_1u_3 \not\in E(P)$, then $y_1v_2 \in E(P)$. Let
$E_3=\{y_1v_2,z_0z_2,y_3u_2 \}$. We can also get that
$|E(P)\bigcap E_3|=2$ and a $v_2-u_2$ hamiltonian path in the
subgraph $H(z_2)$.

Now if there is a $v_1-u_1$ hamiltonian path in the subgraph
$H(z_1)$, then the graph $H(z_1)+u_1v_1$ must be hamiltonian.
According to the Grinberg's criterion for planar hamiltonian
graphs, we know that

$$ {\phi}'_3 - {\phi}"_3+2({\phi}'_4-{\phi}"_4)
+3({\phi}'_5-{\phi}"_5)+6({\phi}'_8-{\phi}"_8)=0,  \hspace{20mm}
(*) $$

\no{where ${\phi}'_i$ or ${\phi}"_i$ is the number of $i$-gons in
the interior or exterior of a chosen hamiltonian circuit $C$
passing through $u_1v_1$ in the graph $H(z_1)+u_1v_1$. Since it is
obvious that}

$${\phi}'_3={\phi}"_8=1,\quad {\phi}"_3={\phi}'_8=0,  $$

\no{we get that}

$$2({\phi}'_4-{\phi}"_4)+3({\phi}'_5-{\phi}"_5)=5, \hspace{10mm} (**)  $$

\no by (*).

Because ${\phi}'_4+{\phi}"_4=2$, so ${\phi}'_4-{\phi}"_4=0,2 \
{\rm or} \ -2$. Now the valency of $z_1$ in $H(z_1)$ is $2$, so
the $4$-gon containing the vertex $z_1$ must be in the interior of
$C$, that is ${\phi}'_4-{\phi}"_4 \not= -2$. If
${\phi}'_4-{\phi}"_4=0 \ {\rm or} \ {\phi}'_4-{\phi}"_4=2$, we get
$3({\phi}'_5-{\phi}"_5)=5 \ {\rm or} \ 3({\phi}'_5-{\phi}"_5)=1$,
a contradiction.

Notice that $H(z_1)\cong H(z_2)\cong H(z_3)$. If there exists a
$v_2 - u_2$ hamiltonian path in $H(z_2)$, a contradiction can be
also gotten. So there does not exist a $y_1 - y_2$ hamiltonian
path in the graph $N(\sigma(v))$. Similarly , there are no
$y_1-y_3 \ {\rm or} \ y_2-y_3$ hamiltonian paths in the graph
$N(\sigma(v))$. Whence, $M^{\sigma(v)}$ is non-hamiltonian.

Now let $n$ be an integer, $n\geq 1$. We get that

\begin{eqnarray*}
M_1 &=&  (M)^{\sigma(u)},\quad  u\in{V(M)}   ; \\
M_2 &=&  (M_1)^{N(\sigma(v))(v)},\quad  v\in{V(M_1)}   ; \\
\cdots &\cdots& \cdots\cdots\cdots\cdots\cdots\cdots\cdots     ;  \\
M_n &=&  (M_{n-1})^{N(\sigma(v))(w)},\quad w\in{V(M_{n-1})}; \\
\cdots &\cdots& \cdots\cdots\cdots\cdots\cdots\cdots\cdots.
\end{eqnarray*}

\no All of these maps are 3-connected non-hamiltonian cubic maps
on the surface $S$. This completes the proof. \quad\quad
$\natural$

\vskip 4mm

\no{\bf Corollary $2.3.5$} \ {\it There is not a locally
orientable surface on which every 3-connected cubic map is
hamiltonian.}

\vskip 4mm

\no{\bf $2.3.3.$ Multi-Embeddings in an $n$-manifold}

\vskip 3mm

\no We come back to determine multi-embeddings of graphs in this
subsection. Let $S_1,S_2,\cdots,S_k$ be $k$ locally orientable
surfaces and $G$ a connected graph. Define numbers

$$\gamma(G;S_1,S_2,\cdots,S_k)=
\min\{ \sum\limits_{i=1}^k\gamma(G_i)|
G=\biguplus\limits_{i=1}^kG_i, G_i\rightarrow S_i, 1\leq i\leq k
\},$$

$$\gamma_M(G;S_1,S_2,\cdots,S_k)=
\max\{\sum\limits_{i=1}^k\gamma(G_i)|
G=\biguplus\limits_{i=1}^kG_i, G_i\rightarrow S_i, 1\leq i\leq k
\}.$$

\no and the {\it multi-genus range $GR(G;S_1,S_2,\cdots,S_k)$} by

$$GR(G;S_1,S_2,\cdots,S_k)=\{ \sum\limits_{i=1}^kg(G_i)|G=\biguplus\limits_{i=1}^kG_i,
G_i\rightarrow S_i, 1\leq i\leq k \},$$

\no where $G_i$ is embeddable on a surface of genus $g(G_i)$. Then
we get the following result.

\vskip 4mm

\no{\bf Theorem $2.3.18$} \ {\it Let $G$ be a connected graph and
let $S_1,S_2,\cdots,S_k$ be locally orientable surfaces with empty
overlapping. Then }

$$GR(G;S_1,S_2,\cdots,S_k)=
[\gamma(G;S_1,S_2,\cdots,S_k),\gamma_M(G;S_1,S_2,\cdots,S_k)].$$

\vskip 3mm

{\it Proof} \ Let $G=\biguplus\limits_{i=1}^kG_i, G_i\rightarrow
S_i, 1\leq i\leq k$. We prove that there are no gap in the
multi-genus range from
$\gamma(G_1)+\gamma(G_2)+\cdots+\gamma(G_k)$ to
$\gamma_M(G_1)+\gamma_M(G_2)+\cdots+\gamma_M(G_k)$. According to
Theorems $2.3.8$ and $2.3.12$, we know that the genus range
$GR^O(G_i)$ or $GR^N(G)$ is $[\gamma(G_i),\gamma_M(G_i)]$ or
$[\widetilde{\gamma}(G_i),\widetilde{\gamma}_M(G_i)]$ for any
integer $i, 1\leq i\leq k$. Whence, there exists a multi-embedding
of $G$ on $k$ locally orientable surfaces $P_1,P_2,\cdots,P_k$
with $g(P_1)=\gamma (G_1)$, $g(P_2)=\gamma (G_2)$,$\cdots$,
$g(P_k)=\gamma (G_k)$. Consider the graph $G_1$, then $G_2$, and
then $G_3$, $\cdots$ to get multi-embedding of $G$ on $k$ locally
orientable surfaces step by step. We get a multi-embedding of $G$
on $k$ surfaces with genus sum at least being an unbroken interval
$[\gamma(G_1)+\gamma(G_2)+\cdots+\gamma(G_k),
\gamma_M(G_1)+\gamma_M(G_2)+\cdots+\gamma_M(G_k)]$ of integers.

By definitions of $\gamma(G;S_1,S_2,\cdots,S_k)$ and
$\gamma_M(G;S_1,S_2,\cdots,S_k)$, we assume that
$G=\biguplus\limits_{i=1}^kG'_i, G'_i\rightarrow S_i, 1\leq i\leq
k$ and $G=\biguplus\limits_{i=1}^kG''_i, G''_i\rightarrow S_i,
1\leq i\leq k$ attain the extremal values
$\gamma(G;S_1,S_2,\cdots,S_k)$ and
$\gamma_M(G;S_1,S_2,\cdots,S_k)$, respectively. Then we know that
the multi-embedding of $G$ on $k$ surfaces with genus sum is at
least an unbroken intervals $[\sum\limits_{i=1}^k\gamma(G'_i),
\sum\limits_{i=1}^k\gamma_M(G'_i)]$ and
$[\sum\limits_{i=1}^k\gamma(G''_i),
\sum\limits_{i=1}^k\gamma_M(G''_i)]$ of integers.

Since

$$\sum\limits_{i=1}^kg(S_i)\in[\sum\limits_{i=1}^k\gamma(G'_i),
\sum\limits_{i=1}^k\gamma_M(G'_i)]\bigcap
[\sum\limits_{i=1}^k\gamma(G''_i),
\sum\limits_{i=1}^k\gamma_M(G''_i)],$$

\no we get that

$$GR(G;S_1,S_2,\cdots,S_k)=
[\gamma(G;S_1,S_2,\cdots,S_k),\gamma_M(G;S_1,S_2,\cdots,S_k)].$$

\no This completes the proof.\quad\quad $\natural$

For multi-embeddings of a complete graph, we get the following
result.

\vskip 4mm

\no{\bf Theorem $2.3.19$} \ {\it Let $P_1,P_2,\cdots,P_k$ and
$Q_1,Q_2,\cdots,Q_k$ be respective $k$ orientable and
non-orientable surfaces of genus$\geq 1$. A complete graph $K_n$
is multi-embeddable in $P_1,P_2,\cdots,P_k$ with empty overlapping
if and only if}

$$\sum\limits_{i=1}^k\lceil\frac{3+\sqrt{16g(P_i)+1}}{2}\rceil\leq
n\leq\sum\limits_{i=1}^k\lfloor\frac{7+\sqrt{48g(P_i)+1}}{2}\rfloor$$

\no{\it and is multi-embeddable in $Q_1,Q_2,\cdots,Q_k$ with empty
overlapping if and only if}

$$\sum\limits_{i=1}^k\lceil 1+\sqrt{2g(Q_i)}\rceil\leq
n\leq\sum\limits_{i=1}^k\lfloor\frac{7+\sqrt{24g(Q_i)+1}}{2}\rfloor.$$

\vskip 3mm

{\it Proof} \ According to Theorem $2.3.9$ and Corollary $2.3.2$,
we know that the genus $g(P)$ of an orientable surface $P$ on
which a complete graph $K_n$ is embeddable satisfies

$$\lceil\frac{(n-3)(n-4)}{12}\rceil\leq g(P)\leq\lfloor\frac{(n-1)(n-2)}{4}\rfloor,$$

\no i.e.,

$$\frac{(n-3)(n-4)}{12}\leq g(P)\leq\frac{(n-1)(n-2)}{4}.$$

If $g(P)\geq 1$, we get that

$$\lceil\frac{3+\sqrt{16g(P)+1}}{2}\rceil\leq n
\leq\lfloor\frac{7+\sqrt{48g(P)+1}}{2}\rfloor.$$

Similarly, if $K_n$ is embeddable on a non-orientable surface $Q$,
then

$$\lceil\frac{(n-3)(n-4)}{6}\rceil\leq g(Q)\leq\lfloor\frac{(n-1)^2}{2}\rfloor,$$

\no i.e.,

$$\lceil 1+\sqrt{2g(Q)}\rceil\leq n\leq\lfloor\frac{7+\sqrt{24g(Q)+1}}{2}\rfloor.$$

Now if $K_n$ is multi-embeddable in $P_1,P_2,\cdots,P_k$ with
empty overlapping, then there must exists a partition
$n=n_1+n_2+\cdots+n_k$, $n_i\geq 1, 1\leq i\leq k$. Since each
vertex-induced subgraph of a complete graph is still a complete
graph, we know that for any integer $i, 1\leq i\leq k$,

$$\lceil\frac{3+\sqrt{16g(P_i)+1}}{2}\rceil\leq n_i
\leq\lfloor\frac{7+\sqrt{48g(P_i)+1}}{2}\rfloor.$$

\no Whence, we know that

$$\sum\limits_{i=1}^k\lceil\frac{3+\sqrt{16g(P_i)+1}}{2}\rceil\leq
n\leq\sum\limits_{i=1}^k\lfloor\frac{7+\sqrt{48g(P_i)+1}}{2}\rfloor.\quad\quad
(*)$$

On the other hand, if the inequality (*) holds, we can find
positive integers $n_1,n_2,\cdots,n_k$ with $n=n_1+n_2+\cdots+n_k$
and

$$\lceil\frac{3+\sqrt{16g(P_i)+1}}{2}\rceil\leq n_i
\leq\lfloor\frac{7+\sqrt{48g(P_i)+1}}{2}\rfloor.$$

\no for any integer $i, 1\leq i\leq k$. This enables us to
establish a partition $K_n=\biguplus\limits_{i=1}^kK_{n_i}$ for
$K_n$ and embed each $K_{n_i}$ on $P_i$ for $1\leq i\leq k$.
Therefore, we get a multi-embedding of $K_n$ in
$P_1,P_2,\cdots,P_k$ with empty overlapping.

Similarly, if $K_n$ is multi-embeddable in $Q_1,Q_2,\cdots Q_k$
with empty overlapping, there must exists a partition
$n=m_1+m_2+\cdots+m_k$, $m_i\geq 1, 1\leq i\leq k$ and

$$\lceil 1+\sqrt{2g(Q_i)}\rceil\leq m_i\leq\lfloor\frac{7+\sqrt{24g(Q_i)+1}}{2}\rfloor.$$

\no for any integer $i, 1\leq i\leq k$. Whence, we get that

$$\sum\limits_{i=1}^k\lceil 1+\sqrt{2g(Q_i)}\rceil\leq
n\leq\sum\limits_{i=1}^k\lfloor\frac{7+\sqrt{24g(Q_i)+1}}{2}\rfloor.
\quad\quad (**)$$

Now if the inequality (**) holds, we can also find positive
integers $m_1,m_2,\cdots,m_k$ with $n=m_1+m_2+\cdots+m_k$ and

$$\lceil 1+\sqrt{2g(Q_i)}\rceil\leq m_i\leq\lfloor\frac{7+\sqrt{24g(Q_i)+1}}{2}\rfloor.$$

\no for any integer $i, 1\leq i\leq k$. Similar to those of
orientable cases, we get a multi-embedding of $K_n$ in
$Q_1,Q_2,\cdots,Q_k$ with empty overlapping.\quad\quad $\natural$

\vskip 4mm

\no{\bf Corollary $2.3.6$} \ {\it A complete graph $K_n$ is
multi-embeddable in $k, k\geq 1$ orientable surfaces of genus $p,
p\geq 1$ with empty overlapping if and only if}

$$\lceil\frac{3+\sqrt{16p+1}}{2}\leq\frac{n}{k}
\leq\lfloor\frac{7+\sqrt{48p+1}}{2}\rfloor$$

\no{\it and is multi-embeddable in $l, l\geq 1$ non-orientable
surfaces of genus $q, q\geq 1$ with empty overlapping if and only
if}

$$\lceil 1+\sqrt{2q}\rceil\leq
\frac{n}{k}\leq\lfloor\frac{7+\sqrt{24q+1}}{2}\rfloor.$$

\vskip 4mm

\no{\bf Corollary $2.3.7$} \ {\it A complete graph $K_n$ is
multi-embeddable in $s, s\geq 1$ tori with empty overlapping if
and only if}

$$4s\leq n\leq 7s$$

\no{\it and is multi-embeddable in $t, t\geq 1$ projective planes
with empty overlapping if and only if}

$$3t\leq n\leq 6t.$$

\vskip 3mm

Similarly, the following result holds for a complete bipartite
graph $K(n,n)$.

\vskip 4mm

\no{\bf Theorem $2.3.20$} \ {\it \it Let $P_1,P_2,\cdots,P_k$ and
$Q_1,Q_2,\cdots,Q_k$ be respective $k$ orientable and $k$
non-orientable surfaces of genus$\geq 1$. A complete bipartite
graph $K(n,n)$ is multi-embeddable in $P_1,P_2,\cdots,P_k$ with
empty overlapping if and only if}

$$\sum\limits_{i=1}^k\lceil1+\sqrt{2g(P_i)}\rceil\leq n
\leq\sum\limits_{i=1}^k\lfloor 2+2\sqrt{g(P_i)}\rfloor$$

\no{\it and is multi-embeddable in $Q_1,Q_2,\cdots,Q_k$ with empty
overlapping if and only if}

$$\sum\limits_{i=1}^k\lceil 1+\sqrt{g(Q_i)}\rceil\leq n
\leq\sum\limits_{i=1}^k\lfloor 2+\sqrt{2g(Q_i)}\rfloor.$$

\vskip 4mm

{\it Proof} \ Similar to the proof of Theorem $2.3.18$, we get
this result. \quad\quad $\natural$

\vskip 4mm

\no{\bf $2.3.4.$ Classification of graphs in an $n$-manifold}

\vskip 3mm

\no By Theorem $2.3.1$ we can give a combinatorial definition for
a graph embedded in an $n$-manifold, i.e., a {\it manifold graph}
similar to the Tutte's definition for a map.

\vskip 4mm

\no{\bf Definition $2.3.6$} \ {\it For any integer $n, n\geq 2$,
an $n$-dimensional manifold graph $^n{\mathcal G}$ is a pair
$^n{\mathcal G} = ({\mathcal E}_{\Gamma},\mathcal L)$ in where a
permutation ${\mathcal L}$ acting on ${\mathcal E}_{\Gamma}$ of a
disjoint union $\Gamma x =\{\sigma x| \sigma\in\Gamma\}$ for
$\forall x\in E$, where $E$ is a finite set and $\Gamma =\{\mu,o |
\mu^2=o^n=1, \mu o=o\mu\}$ is a commutative group of order $2n$
with the following conditions hold.}

($i$) \ {\it $\forall x\in {\mathcal E}_K$, there does not exist
an integer $k$ such that ${\mathcal L }^kx = o^i x$ for $\forall
i, 1\leq i\leq n-1$;}

($ii$) \ $\mu {\mathcal L } = {\mathcal L}^{-1}\mu$;

($iii$) \ {\it The group $\Psi_{J}=\left<\mu,o,{\mathcal
L}\right>$ is transitive on ${\mathcal E}_{\Gamma}$.}

\vskip 3mm

According to ($i$) and ($ii$), a {\it vertex $v$ of an
$n$-dimensional manifold graph} is defined to be an $n$-tuple
$\{(o^ix_1,o^ix_2,\cdots,o^ix_{s_l(v)})(o^iy_1,o^iy_2,\cdots,o^iy_{s_2(v)})
\cdots(o^iz_1,o^iz_2,$ $\cdots,o^iz_{s_{l(v)}(v)}); 1\leq i\leq
n\}$ of permutations of ${\mathcal L}$ action on ${\mathcal
E}_{\Gamma}$, edges to be these orbits of $\Gamma$ action on
${\mathcal E}_{\Gamma}$. The number
$s_1(v)+s_2(v)+\cdots+s_{l(v)}(v)$ is called the {\it valency of
$v$}, denoted by $\rho_G^{s_1,s_2,\cdots,s_{l(v)}}(v)$. The
condition ($iii$) is used to ensure that an $n$-dimensional
manifold graph is connected. Comparing definitions of a map with
an $n$-dimensional manifold graph, the following result holds.

\vskip 4mm

\no{\bf Theorem $2.3.21$} \ {\it For any integer $n, n\geq 2$,
every $n$-dimensional manifold graph $^n{\mathcal G} = ({\mathcal
E}_{\Gamma},\mathcal L)$ is correspondent to a unique map
$M=({\mathcal E}_{\alpha,\beta},{\mathcal P})$ in which each
vertex $v$ in $^n{\mathcal G}$ is converted to $l(v)$ vertices
$v_1,v_2,\cdots,v_{l(v)}$ of $M$. Conversely, a map $M=({\mathcal
E}_{\alpha,\beta},{\mathcal P})$ is also correspondent to an
$n$-dimensional manifold graph $^n{\mathcal G} = ({\mathcal
E}_{\Gamma},\mathcal L)$ in which $l(v)$ vertices
$u_1,u_2,\cdots,u_{l(v)}$ of $M$ are converted to one vertex $u$
of $^n{\mathcal G}$.}

\vskip 3mm

Two $n$-dimensional manifold graphs $^n{\mathcal G}_1=({\mathcal
E}_{\Gamma_1}^1,{\mathcal L}_1)$ and $^n{\mathcal G}_2=({\mathcal
E}_{\Gamma_2}^2,{\mathcal L}_2)$ are said to be {\it isomorphic}
if there exists a one-to-one mapping $\kappa:{\mathcal
E}_{\Gamma_1}^1\rightarrow{\mathcal E}_{\Gamma_2}^2$ such that
$\kappa\mu=\mu\kappa, \kappa o=o\kappa$ and $\kappa {\mathcal
L}_1={\mathcal L}_2\kappa$. If ${\mathcal
E}_{\Gamma_1}^1={\mathcal E}_{\Gamma_2}^2={\mathcal E}_{\Gamma}$
and ${\mathcal L}_1={\mathcal L}_2={\mathcal L}$, an isomorphism
between $^n{\mathcal G}_1$ and $^n{\mathcal G}_2$ is called an
automorphism of $^n{\mathcal G} = ({\mathcal E}_{\Gamma},\mathcal
L)$. It is immediately that all automorphisms of $^n{\mathcal G}$
form a group under the composition operation. We denote this group
by ${\rm Aut}^n{\mathcal G}$.

It is obvious that for two isomorphic $n$-dimensional manifold
graphs $^n{\mathcal G}_1$ and $^n{\mathcal G}_2$, their underlying
graphs $G_1$ and $G_2$ are isomorphic. For an embedding
$^n{\mathcal G}=({\mathcal E}_{\Gamma},{\mathcal L})$ in an
$n$-dimensional manifold and $\forall \zeta\in{\rm
Aut}_{\frac{1}{2}}G$, an induced action of $\zeta$ on ${\mathcal
E}_{\Gamma}$ is defined by

$$ \zeta(gx)=g\zeta(x)$$

\no for $\forall x\in{\mathcal E}_{\Gamma}$ and $\forall
g\in\Gamma$. Then the following result holds.

\vskip 4mm

\no{\bf Theorem $2.3.22$} \  \ \ { ${\rm Aut}^n{\mathcal G}\preceq
{\rm Aut}_{\frac{1}{2}}G\times\left<\mu\right>$.}

\vskip 3mm

{\it Proof} \ First we prove that two $n$-dimensional manifold
graphs $^n{\mathcal G}_1=({\mathcal E}_{\Gamma_1}^1,{\mathcal
L}_1)$ and$^n{\mathcal G}_2=({\mathcal E}_{\Gamma_2}^2,{\mathcal
L_2})$ are isomorphic if and only if there is an element $\zeta\in
{\rm Aut}_{\frac{1}{2}}\Gamma$ such that ${\mathcal
L}_1^{\zeta}={\mathcal L}_2$ or ${\mathcal L}_2^{-1}$.

If there is an element $\zeta\in {\rm Aut}_{\frac{1}{2}}\Gamma$
such that ${\mathcal L}_1^{\zeta}={\mathcal L}_2$, then the
$n$-dimensional manifold graph $^n{\mathcal G}_1$ is isomorphic to
$^n{\mathcal G}_2$ by definition. If ${\mathcal
L}_1^{\zeta}={\mathcal L}_2^{-1}$, then ${\mathcal
L}_1^{\zeta\mu}={\mathcal L}_2$. The $n$-dimensional manifold
graph $^n{\mathcal G}_1$ is also isomorphic to $^n{\mathcal G}_2$.

By the definition of an isomorphism $\xi$ between $n$-dimensional
manifold graphs $^n{\mathcal G}_1$ and $^n{\mathcal G}_2$, we know
that

$$\mu\xi(x)=\xi\mu (x), \ o\xi (x)=\xi o(x) \ {\rm
and} \ {\mathcal L}_1^{\xi }(x)= {\mathcal L}_2(x).$$

\no $\forall x\in {\mathcal E}_{\Gamma}$. By definition these
conditions

$$o\xi (x)=\xi o (x) \ {\rm and} \ {\mathcal L}_1^{\xi }(x)= {\mathcal L}_2(x).$$

\no are just the condition of an automorphism $\xi$ or $\alpha\xi$
on $X_{\frac{1}{2}}(\Gamma)$. Whence, the assertion is true.

Now let ${\mathcal E}_{\Gamma_1}^1={\mathcal
E}_{\Gamma_2}^2={\mathcal E}_{\Gamma}$ and ${\mathcal
L}_1={\mathcal L}_2={\mathcal L}$. We know that

$${\rm Aut}^n{\mathcal G}\preceq {\rm
Aut}_{\frac{1}{2}}G\times\left<\mu\right>. \quad\quad \natural$$

Similar to combinatorial maps, the action of an automorphism of a
manifold graph on ${\mathcal E}_{\Gamma}$ is fixed-free.

\vskip 4mm

\no{\bf Theorem $2.3.23$} \ {\it Let $^n{\mathcal G}=({\mathcal
E}_{\Gamma},{\mathcal L})$ be an $n$-dimensional manifold graph.
Then $({\rm Aut}^n{\mathcal G})_x$ is trivial for $\forall x\in
{\mathcal E}_{\Gamma}$.}

\vskip 3mm

{\it Proof} \ For $\forall g\in ({\rm Aut}^n{\mathcal G})_x$, we
prove that $g(y) = y$ for $\forall y\in{\mathcal E}_{\Gamma}$. In
fact, since the group $\Psi_{J}=\left<\mu,o,{\mathcal L}\right>$
is transitive on ${\mathcal E}_{\Gamma}$, there exists an element
$\tau\in\Psi_{J}$ such that $y = \tau(x)$. By definition we know
that every element in $\Psi_{J}$ is commutative with automorphisms
of $^n{\mathcal G}$. Whence, we get that

$$g(y)= g(\tau(x))=\tau(g(x))=\tau(x)=y.$$

\no i.e., $({\rm Aut}^n{\mathcal G})_x$ is trivial. \quad\quad
$\natural$

\vskip 4mm

\no{\bf Corollary $2.3.8$} \ {\it Let $M=({\mathcal
X}_{\alpha,\beta},{\mathcal P})$ be a map. Then for $\forall
x\in{\mathcal X}_{\alpha,\beta}$, $({\rm Aut}M)_x$ is trivial.}

\vskip 3mm

For an $n$-dimensional manifold graph $^n{\mathcal G}=({\mathcal
E}_{\Gamma},{\mathcal L})$, an $x\in{\mathcal E}_{\Gamma}$ is said
a {\it root} of $^n{\mathcal G}$. If we have chosen a root $r$ on
an $n$-dimensional manifold graph $^n{\mathcal G}$, then
$^n{\mathcal G}$ is called a {\it rooted $n$-dimensional manifold
graph}, denoted by $^n{\mathcal G}^r$. Two rooted $n$-dimensional
manifold graphs $^n{\mathcal G}^{r_1}$ and $^n{\mathcal G}^{r_2}$
are said to be {\it isomorphic} if there is an isomorphism
$\varsigma$ between them such that $\varsigma(r_1)=r_2$. Applying
Theorem $2.3.23$ and Corollary $2.3.1$, we get an enumeration
result for $n$-dimensional manifold graphs underlying a graph $G$
in the following.

\vskip 4mm

\no{\bf Theorem $2.3.24$} \ {\it For any integer $n, n\geq 3$, the
number $r_n^S(G)$ of rooted $n$-dimensional manifold graphs
underlying a graph $G$ is}

$$r_n^S(G)=\frac {n\varepsilon(G)
\prod\limits_{v\in{V(G)}}\rho_G(v)!}{|{Aut_{\frac{1}{2}}G }|}.
$$

\vskip 3mm

{\it Proof} \ Denote the set of all non-isomorphic $n$-dimensional
manifold graphs underlying a graph $G$ by ${\mathcal G}^S(G)$. For
an $n$-dimensional graph $^n{\mathcal G}=({\mathcal
E}_{\Gamma},{\mathcal L})\in{\mathcal G}^S(G)$, denote the number
of non-isomorphic rooted $n$-dimensional manifold graphs
underlying $^n{\mathcal G}$ by $r(^n{\mathcal G})$. By a result in
permutation groups theory, for $\forall x\in{\mathcal E}_{\Gamma}
$ we know that

$$|{\rm Aut}^n{\mathcal G}|= |({\rm Aut}^n{\mathcal G})_x||x^{{\rm Aut}^n{\mathcal G}}|.$$

According to Theorem $2.3.23$, $|({\rm Aut}^n{\mathcal G})_x|=1$.
Whence, $|x^{{\rm Aut}^n{\mathcal G}}|=|{\rm Aut}^n{\mathcal G}|$.
However there are $|{\mathcal E}_{\Gamma}|= 2n\varepsilon(G)$
roots in $^n{\mathcal G}$ by definition. Therefore, the number of
non-isomorphic rooted $n$-dimensional manifold graphs underlying
an $n$-dimensional graph $^n{\mathcal G}$ is

$$r(^n{\mathcal G})= \frac{|{\mathcal E}_{\Gamma}|}{|{\rm Aut}^n{\mathcal G}|}
= \frac{2n\varepsilon(G)}{|{\rm Aut}^n{\mathcal G}|}.$$

\no Whence, the number of non-isomorphic rooted $n$-dimensional
manifold graphs underlying a graph $G$ is

$$
r_n^S(G)=\sum\limits_{^n{\mathcal G}\in{\mathcal
G}^S(G)}\frac{2n\varepsilon (G)}{|{\rm Aut}^n{\mathcal G}|}.
$$

\no According to Theorem $2.3.22$, ${\rm Aut}^n{\mathcal G}\preceq
{\rm Aut}_{\frac{1}{2}}G\times\left<\mu\right>$. Whence $\tau\in
{\rm Aut}^n{\mathcal G}$ for $^n{\mathcal G}\in{\mathcal G}^S(G)$
if and only if $\tau\in({\rm
Aut}_{\frac{1}{2}}G\times\left<\mu\right>)_{^n{\mathcal G}}$.
Therefore, we know that ${\rm Aut}^n{\mathcal G}=({\rm Aut}
_{\frac{1}{2}}G\times\left<\mu\right>)_{^n{\mathcal G}}$. Because
of $|{\rm Aut}_{\frac{1}{2}}G\times\left<\mu\right>| = |({\rm Aut}
_{\frac{1}{2}}G\times\left<\mu\right>)_{^n{\mathcal G}}|
|^n{\mathcal G}^{{\rm Aut}_{\frac{1}{2}}
G\times\left<\mu\right>}|$, we get that

$$
|^n{\mathcal G}^{{\rm Aut}
_{\frac{1}{2}}G\times\left<\mu\right>}|= \frac{2|{\rm
Aut}_{\frac{1}{2}}G|}{|{\rm Aut}^n{\mathcal G}|}.
$$

\no Therefore,

\begin{eqnarray*}
r_n^S(G) &=& \sum\limits_{^n{\mathcal G}\in{\mathcal
G}^S(G)}\frac{2n\varepsilon (G)}{|{\rm Aut}^n{\mathcal G}|}\\
&=& \frac{2n\varepsilon (G)}{|{\rm Aut}
_{\frac{1}{2}}G\times\left<\mu\right>|} \sum\limits_{^n{\mathcal
G}\in{\mathcal G}^S(G)}\frac{|{\rm Aut}
_{\frac{1}{2}}G\times\left<\mu\right>|}
{|{\rm Aut}^n{\mathcal G}|}\\
&=&  \frac{2n\varepsilon (G)}{|{\rm Aut}
_{\frac{1}{2}}G\times\left<\mu\right>|} \sum\limits_{^n{\mathcal
G}\in{\mathcal G}^S(G)}{|^n{\mathcal G}^{{\rm Aut}
_{\frac{1}{2}}G\times\left<\mu\right>|}}\\
&=&  \frac {n\varepsilon (G) \prod\limits_{v \in V(G)}
\rho_G(v)!}{|{\rm Aut_{\frac{1}{2}}}G|}
\end{eqnarray*}

\no by applying Corollary $2.3.1$. \quad\quad $\natural$

Notice the fact that an embedded graph in a $2$-dimensional
manifolds is just a map. Then Definition $3.6$ is converted to
Tutte's definition for combinatorial maps in this case. We can
also get an enumeration result for rooted maps on surfaces
underlying a graph $G$ by applying Theorems $2.3.7$ and $2.3.23$
in the following.

\vskip 4mm

\no{\bf Theorem $2.3.25$}([66],[67]) \ {\it The number
$r^L(\Gamma)$ of rooted maps on locally orientable surfaces
underlying a connected graph $G$ is}

$$r^L(G)=\frac {2^{\beta (G) +1}\varepsilon (G)
\prod\limits_{v\in{V(G)}}(\rho (v)-1)!}{|{Aut_{\frac{1}{2}}G }|},
$$

\no{\it where $\beta (G)=\varepsilon(G)-\nu(G)+1$ is the Betti
number of $G$.}

\vskip 3mm

Similarly, for a graph $G=\bigoplus\limits_{i=1}^lG_i$ and a
multi-manifold $\widetilde{M}=\bigcup\limits_{i=1}^l {\bf
M}^{l_i}$, choose $l$ commutative groups
$\Gamma_1,\Gamma_2,\cdots,\Gamma_l$, where
$\Gamma_i=\left<\mu_i,o_i|\mu_i^2=o^{h_i}=1\right>$ for any
integer $i, 1\leq i\leq l$. Consider permutations acting on
$\bigcup\limits_{i=1}^l{\mathcal E}_{\Gamma_i}$, where for any
integer $i,1\leq i\leq l$, ${\mathcal E}_{\Gamma_i}$ is a disjoint
union $\Gamma_i x =\{\sigma_i x| \sigma_i\in\Gamma\}$ for $\forall
x\in E(G_i)$. Similar to Definition $2.3.6$, we can also get a
multi-embedding of $G$ in
$\widetilde{M}=\bigcup\limits_{i=1}^l{\bf M}^{h_i}$.

\vskip 5mm

\no{\bf \S $2.4$ \ Multi-Spaces on Graphs}

\vskip 4mm

\no A Smarandache multi-space is a union of $k$ spaces
$A_1,A_2,\cdots,A_k$ for an integer $k, k\geq 2$ with some
additional constraint conditions. For describing a finite
algebraic multi-space, graphs are a useful way. All graphs
considered in this section are directed graphs.

\vskip 4mm

\no{\bf $2.4.1.$ A graph model for an operation system}

\vskip 3mm

\no A graph is called a {\it directed graph} if there is an
orientation on its every edge. A directed graph
$\overrightarrow{G}$ is called an {\it Euler graph} if we can
travel all edges of $\overrightarrow{G}$ alone orientations on its
edges with no repeat starting at any vertex $u\in
V($\overrightarrow{G}$)$ and come back to $u$. For a directed
graph $\overrightarrow{G}$, we use the convention that the
orientation on the edge $e$ is $u\rightarrow v$ for $\forall
e=(u,v)\in E(\overrightarrow{G})$ and say that $e$ is {\it
incident from $u$} and {\it incident to $v$}. For $u\in
V(\overrightarrow{G})$, the {\it outdegree
$\rho_{\overrightarrow{G}}^+(u)$} of $u$ is the number of edges in
$\overrightarrow{G}$ incident from $u$ and the {\it indegree
$\rho_{\overrightarrow{G}}^-(u)$} of $u$ is the number of edges in
$\overrightarrow{G}$ incident to $u$. Whence, we know that

$$\rho_{\overrightarrow{G}}^+(u)+\rho_{\overrightarrow{G}}^-(u)=
\rho_{\overrightarrow{G}}(u).$$

\no It is well-known that a graph $\overrightarrow{G}$ is Eulerian
if and only if
$\rho_{\overrightarrow{G}}^+(u)=\rho_{\overrightarrow{G}}^-(u)$
for $\forall u\in V(\overrightarrow{G})$, seeing examples in
$[11]$ for details. For a multiple $2$-edge $(a,b)$, if two
orientations on edges are both to $a$ or both to $b$, then we say
it to be a {\it parallel multiple $2$-edge}. If one orientation is
to $a$ and another is to $b$, then we say it to be an {\it
opposite multiple $2$-edge}.

Now let $(A; \circ)$ be an algebraic system with operation
¡°$\circ$¡±. We associate a {\it weighted graph $G[A]$} for $(A;
\circ)$ defined as in the next definition.

\vskip 4mm

\no{\bf Definition $2.4.1$} \ {\it Let $(A; \circ)$ be an
algebraic system. Define a weighted graph $G[A]$ associated with
$(A; \circ)$ by}

$$V(G[A]) \ = \ A $$

\no{\it and}

$$E(G[A]) = \{(a,c) \ with \ weight \ \circ b \ | \ if \ a\circ b=c \ for \ \forall a,b,c\in A\}$$

\no{\it as shown in {\rm Fig.$2.29$}.}

\vskip 3mm

\includegraphics[bb=40 10 500 60]{sgm31.eps}

\vskip 3mm

\c{\bf Fig.$2.29$}\vskip 2mm

\vskip 3mm

For example, the associated graph $G[Z_4]$ for the commutative
group $Z_4$ is shown in Fig.$2.30$.

\vskip 3mm

\includegraphics[bb=40 10 500 140]{sgm32.eps}

\vskip 3mm

\c{\bf Fig.$2.30$}\vskip 2mm

The advantage of Definition $2.4.1$ is that for any edge in
$G[A]$, if its vertices are a,c with a weight $\circ b$, then
a$\circ b=c$ and vice versa, if a$\circ b=c$, then there is one
and only one edge in $G[A]$ with vertices $a,c$ and weight $\circ
b$. This property enables us to find some structure properties of
$G[A]$ for an algebraic system $(A; \circ)$.

\vskip 3mm

\no{\bf $P1.$} \ {\it $G[A]$ is connected if and only if there are
no partition $A=A_1\bigcup A_2$ such that for $\forall a_1\in
A_1$, $\forall a_2\in A_2$, there are no definition for $a_1\circ
a_2$ in $(A; \circ)$.}

\vskip 2mm

If $G[A]$ is disconnected, we choose one component $C$ and let
$A_1= V(C)$. Define $A_2= V(G[A])\setminus V(C)$. Then we get a
partition $A=A_1\bigcup A_2$ and for $\forall a_1\in A_1$,
$\forall a_2\in A_2$, there are no definition for $a_1\circ a_2$
in $(A; \circ)$, a contradiction and vice versa.

\vskip 3mm

\no{\bf $P2.$} \ {\it If there is a unit ${\bf 1}_A$ in $(A;
\circ)$, then there exists a vertex ${\bf 1}_A$ in $G[A]$ such
that the weight on the edge $({\bf 1}_A,x)$ is $\circ x$ if ${\bf
1}_A\circ x$ is defined in $(A; \circ)$ and vice versa.}

\vskip 3mm

\no{\bf $P3.$} \ {\it For $\forall a\in A$, if $a^{-1}$ exists,
then there is an opposite multiple $2$-edge $({\bf 1}_A, a)$ in
$G[A]$ with weights $\circ a$ and $\circ a^{-1}$, respectively and
vice versa.}

\vskip 3mm

\no{\bf $P4.$} \ {\it For $\forall a,b\in A$ if $a\circ b=b\circ
a$, then there are edges $(a,x)$ and $(b,x)$, $x\in A$ in $(A;
\circ)$ with weights $w(a,x)=\circ b$ and $w(b,x)=\circ a$,
respectively and vice versa.}

\vskip 3mm

\no{\bf $P5.$} \ {\it If the cancellation law holds in $(A;
\circ)$, i.e., for $\forall a,b,c\in A$, if $a\circ b=a\circ c$
then $b=c$, then there are no parallel multiple $2$-edges in
$G[A]$ and vice versa.}

\vskip 2mm

The property $P2,P3,P4$ and $P5$ are gotten by definition. Each of
these cases is shown in Fig.$2.31(1),(2),(3)$ and $(4)$,
respectively.

\includegraphics[bb=40 10 500 140]{sgm33.eps}

\vskip 3mm

\c{\bf Fig.$2.31$}

\vskip 4mm

\no{\bf Definition $2.4.2$} \ {\it An algebraic system $(A;
\circ)$ is called to be a one-way system if there exists a mapping
$\varpi: A\rightarrow A$ such that if $a\circ b\in A$, then there
exists a unique $c\in A$, $c\circ \varpi(b)\in A$. $\varpi$ is
called a one-way function on $(A; \circ)$.}

\vskip 3mm

We have the following results for an algebraic system $(A; \circ)$
with its associated weighted graph $G[A]$.

\vskip 4mm

\no{\bf Theorem $2.4.1$} \ {\it Let $(A; \circ)$ be an algebraic
system with a associated weighted graph $G[A]$. Then}

($i$) \ {\it if there is a one-way function $\varpi$ on $(A;
\circ)$, then $G[A]$ is an Euler graph, and vice versa, if $G[A]$
is an Euler graph, then there exist a one-way function $\varpi$ on
$(A; \circ)$.}

($ii$) \ {\it if $(A; \circ)$ is a complete algebraic system, then
the outdegree of every vertex in $G[A]$ is $|A|$; in addition, if
the cancellation law holds in $(A; \circ)$, then $G[A]$ is a
complete multiple $2$-graph with a loop attaching at each of its
vertices such that each edge between two vertices in $G[A]$ is an
opposite multiple $2$-edge, and vice versa.}

\vskip 3mm

{\it Proof} \ ($i$) \ Assume $\varpi$ is a one-way function
$\varpi$ on $(A; \circ)$. By definition there exists $c\in A$,
$c\circ \varpi(b)\in A$ for $\forall a\in A$, $a\circ b\in A$.
Thereby there is a one-to-one correspondence between edges from
$a$ with edges to $a$. That is,
$\rho_{G[A]}^+(a)=\rho_{G[A]}^-(a)$ for $\forall a\in V(G[A])$.
Therefore, $G[A]$ is an Euler graph.

Now if $G[A]$ is an Euler graph, then there is a one-to-one
correspondence between edges in $E^-=\{e_i^-; 1\leq i\leq k\}$
from a vertex $a$ with edges $E^+=\{e_i^+; 1\leq i\leq k\}$ to the
vertex $a$. For any integer $i, 1\leq i\leq k$, define $\varpi:
w(e_i^-)\rightarrow w(e_i^+)$. Therefore, $\varpi$ is a
well-defined one-way function on $(A; \circ)$.

($ii$) \ If $(A; \circ)$ is complete, then for $\forall a\in A$
and $\forall b\in A$, $a\circ b\in A$. Therefore,
$\rho_{\overrightarrow{G}}^+(a)= |A|$ for any vertex $a\in
V(G[A])$.

If the cancellation law holds in $(A; \circ)$, by $P5$ there are
no parallel multiple $2$-edges in $G[A]$. Whence, each edge
between two vertices is an opposite $2$-edge and weights on loops
are $\circ {\bf 1}_A$.

By definition, if $G[A]$ is a complete multiple $2$-graph with a
loop attaching at each of its vertices such that each edge between
two vertices in $G[A]$ is an opposite multiple $2$-edge, we know
that $(A; \circ)$ is a complete algebraic system with the
cancellation law holding by the definition of $G[A]$.\quad\quad
$\natural$

\vskip 4mm

\no{\bf Corollary $2.4.1$} \ {\it Let $\Gamma$ be a semigroup.
Then $G[\Gamma]$ is a complete multiple $2$-graph with a loop
attaching at each of its vertices such that each edge between two
vertices in $G[A]$ is an opposite multiple $2$-edge.}

\vskip 3mm

Notice that in a group $\Gamma$, $\forall g\in\Gamma$, if
$g^2\not= {\bf 1}_{\Gamma}$, then $g^{-1}\not=g$. Whence, all
elements of order$> 2$ in $\Gamma$ can be classified into pairs.
This fact enables us to know the following result.

\vskip 4mm

\no{\bf Corollary $2.4.2$} \ {\it Let $\Gamma$ be a group of even
order. Then there are opposite multiple $2$-edges in $G[\Gamma]$
such that weights on its $2$ directed edges are the same.}

\vskip 4mm

\no{\bf $2.4.2.$ Multi-Spaces on graphs}

\vskip 3mm

\no Let $\widetilde{\Gamma}$ be a Smarandache multi-space. Its
associated weighted graph is defined in the following.

\vskip 4mm

\no{\bf Definition $2.4.3$} \ {\it Let
$\widetilde{\Gamma}=\bigcup\limits_{i=1}^n\Gamma_i$ be an
algebraic multi-space with $(\Gamma_i; \circ_i)$ being an
algebraic system for any integer $i, 1\leq i\leq n$. Define a
weighted graph $G(\widetilde{\Gamma})$ associated with
$\widetilde{\Gamma}$ by}

$$G(\widetilde{\Gamma})=\bigcup\limits_{i=1}^nG[\Gamma_i],$$

\no{\it where $G[\Gamma_i]$ is the associated weighted graph of
$(\Gamma_i; \circ_i)$ for $1\leq i\leq n$.}

\vskip 3mm

For example, the weighted graph shown in Fig.$2.32$ is
correspondent with a multi-space
$\widetilde{\Gamma}=\Gamma_1\bigcup\Gamma_2\bigcup\Gamma_3$, where
$(\Gamma_1; +)= (Z_3,+)$, $\Gamma_2=\{e,a,b\}$,
$\Gamma_3=\{1,2,a,b\}$ and these operations ¡°$\cdot$¡±on
$\Gamma_2$ and ¡°$\circ$¡± on $\Gamma_3$ are shown in tables
$2.4.1$ and $2.4.2$.

\includegraphics[bb=40 10 500 140]{sgm34.eps}

\vskip 3mm

\c{\bf Fig.$2.32$}

\vskip 3mm

\begin{center}
\begin{tabular}{c|ccc}
$\cdot$ & \ $e$ \ & \ $a$ \ & \ $b$ \ \\ \hline
$e$ & \ $e$ \ & \ $a$ \ & \ $b$ \ \\
$a$ & \ $a$ \ & \ $b$ \ & \ $e$ \ \\
$b$ & \ $b$ \ & \ $e$ \ & \ $a$ \ \\
\end{tabular}
\end{center}

\vskip 2mm

\c{\bf table $2.4.1$}

\vskip 3mm

\begin{center}
\begin{tabular}{c|cccc}
$\circ$ & \ $1$ \ & \ $2$ \ & \ $a$ \ & \ $b$ \ \\ \hline
$1$ & \ * \ & \ $a$ \ & \ $b$ \ & \ * \ \\
$2$ & \ $b$ \ & \ * \ & \ * \ & \ $a$ \ \\
$a$ & \ * \ & \ * \ & \ * \ & \ $1$ \ \\
$b$ & \ * \ & \ * \ & \ $2$ \ & \ * \ \\
\end{tabular}
\end{center}

\vskip 2mm

\c{\bf table $2.4.2$}

\vskip 3mm

Notice that the correspondence between the multi-space
$\widetilde{\Gamma}$ and the weighted graph
$G[\widetilde{\Gamma}]$ is one-to-one. We immediately get the
following result.

\vskip 4mm

\no{\bf Theorem $2.4.2$} \ {\it The mappings $\pi:
\widetilde{\Gamma}\rightarrow G[\widetilde{\Gamma}]$ and
$\pi^{-1}:G[\widetilde{\Gamma}]\rightarrow \widetilde{\Gamma}$ are
all one-to-one.}

\vskip 3mm

According to Theorems $2.4.1$ and $2.4.2$, we get some
consequences in the following.

\vskip 4mm

\no{\bf Corollary $2.4.3$} \ {\it Let
$\widetilde{\Gamma}=\bigcup\limits_{i=1}^n\Gamma_i$ be a
multi-space with an algebraic system $(\Gamma_i; \circ_i)$ for any
integer $i, 1\leq i\leq n$. If for any integer $i, 1\leq i\leq n$,
$G[\Gamma_i]$ is a complete multiple $2$-graph with a loop
attaching at each of its vertices such that each edge between two
vertices in $G[\Gamma_i]$ is an opposite multiple $2$-edge, then
$\widetilde{\Gamma}$ is a complete multi-space.}

\vskip 4mm

\no{\bf Corollary $2.4.4$} \ {\it Let
$\widetilde{\Gamma}=\bigcup\limits_{i=1}^n\Gamma_i$ be a
multi-group with an operation set
$O(\widetilde{\Gamma})=\{\circ_i; 1\leq i\leq n \}$. Then there is
a partition $G[\widetilde{\Gamma}]=\bigcup\limits_{i=1}^nG_i$ such
that each $G_i$ being a complete multiple $2$-graph attaching with
a loop at each of its vertices such that each edge between two
vertices in $V(G_i)$ is an opposite multiple $2$-edge for any
integer $i, 1\leq i\leq n$.}

\vskip 4mm

\no{\bf Corollary $2.4.5$} \ {\it Let $F$ be a body. Then $G[F]$
is a union of two graphs $K^2(F)$ and $K^2(F^*)$, where $K^2(F)$
or $K^2(F^*)$ is a complete multiple $2$-graph with vertex set $F$
or $F^*=F\setminus\{0\}$ and with a loop attaching at each of its
vertices such that each edge between two different vertices is an
opposite multiple $2$-edge.}

\vskip 4mm

\no{\bf $2.4.3.$ Cayley graphs of a multi-group}

\vskip 3mm

\no Similar to the definition of Cayley graphs of a finite
generated group, we can also define {\it Cayley graphs of a finite
generated multi-group}, where a multi-group
$\widetilde{\Gamma}=\bigcup\limits_{i=1}^n\Gamma_i$ is said to be
{\it finite generated} if the group $\Gamma_i$ is finite generated
for any integer $i, 1\leq i\leq n$, i.e.,
$\Gamma_i=\left<x_i,y_i,\cdots,z_{s_i}\right>$. We denote by
$\widetilde{\Gamma}=\left<x_i,y_i,\cdots,z_{s_i}; 1\leq i\leq
n\right>$ if $\widetilde{\Gamma}$ is finite generated by
$\{x_i,y_i,\cdots,z_{s_i}; 1\leq i\leq n\}$.

\vskip 4mm

\no{\bf Definition $2.4.4$} \ {\it Let
$\widetilde{\Gamma}=\left<x_i,y_i,\cdots,z_{s_i}; 1\leq i\leq
n\right>$ be a finite generated multi-group,
$\widetilde{S}=\bigcup\limits_{i=1}^nS_i$, where
$1_{\Gamma_i}\not\in S_i$, $\widetilde{S}^{-1}=\{a^{-1} | a\in
\widetilde{S}\}=\widetilde{S}$ and $\left<S_i\right>=\Gamma_i$ for
any integer $i, 1\leq i\leq n$. A Cayley graph
$Cay(\widetilde{\Gamma}:\widetilde{S})$ is defined by}

$$V(Cay(\widetilde{\Gamma}:\widetilde{S})) = \widetilde{\Gamma}$$

\no{\it and}

$$E(Cay(\widetilde{\Gamma}:\widetilde{S})) =
\{(g,h)| \ if \ there \ exists \ an \ integer \ i,
g^{-1}\circ_ih\in S_i,1\leq i\leq n\}.$$

\vskip 3mm

By Definition $2.4.4$, we immediately get the following result for
Cayley graphs of a finite generated multi-group.

\vskip 4mm

\no{\bf Theorem $2.4.3$} \ {\it For a Cayley graph
$Cay(\widetilde{\Gamma}:\widetilde{S})$ with
$\widetilde{\Gamma}=\bigcup\limits_{i=1}^n\Gamma_i$ and
$\widetilde{S}=\bigcup\limits_{i=1}^nS_i$,}

$$Cay(\widetilde{\Gamma}:\widetilde{S}) = \bigcup\limits_{i=1}^nCay(\Gamma_i:S_i).$$

\vskip 3mm

It is well-known that {\it every Cayley graph of order$\geq 3$ is
$2$-connected}. But in general, a Cayley graph of a multi-group is
not connected. For the connectedness of Cayley graphs of
multi-groups, we get the following result.

\vskip 4mm

\no{\bf Theorem $2.4.4$} \ {\it A Cayley graph
$Cay(\widetilde{\Gamma}:\widetilde{S})$ with
$\widetilde{\Gamma}=\bigcup\limits_{i=1}^n\Gamma_i$ and
$\widetilde{S}=\bigcup\limits_{i=1}^nS_i$ is connected if and only
if for any integer $i, 1\leq i\leq n$, there exists an integer $j,
1\leq j\leq n$ and $j\not=i$ such that
$\Gamma_i\bigcap\Gamma_j\not=\emptyset$.}

\vskip 3mm

{\it Proof} \ According to Theorem $2.4.3$, if there is an integer
$i, 1\leq i\leq n$ such that $\Gamma_i\bigcap\Gamma_j=\emptyset$
for any integer $j, 1\leq j\leq n$, $j\not=i$, then there are no
edges with the form $(g_i,h)$, $g_i\in\Gamma_i$,
$h\in\widetilde{\Gamma}\setminus\Gamma_i$. Thereby
$Cay(\widetilde{\Gamma}:\widetilde{S})$ is not connected.

Notice that $Cay(\widetilde{\Gamma}:\widetilde{S}) =
\bigcup\limits_{i=1}^nCay(\Gamma_i:S_i)$. Not loss of generality,
we assume that $g\in\Gamma_k$ and $h\in\Gamma_l$, where $1\leq
k,l\leq n$ for any two elements $g,h\in\widetilde{\Gamma}$. If
$k=l$, then there must exists a path connecting $g$ and $h$ in
$Cay(\widetilde{\Gamma}:\widetilde{S})$.

Now if $k\not=l$ and for any integer $i, 1\leq i\leq n$, there is
an integer $j, 1\leq j\leq n$ and $j\not=i$ such that
$\Gamma_i\bigcap\Gamma_j\not=\emptyset$, then we can find integers
$i_1,i_2,\cdots,i_s$, $1\leq i_1,i_2,\cdots,i_s\leq n$ such that

$$\Gamma_k\bigcap\Gamma_{i_1}\not=\emptyset,$$

$$\Gamma_{i_1}\bigcap\Gamma_{i_2}\not=\emptyset,$$

$$\cdots\cdots\cdots\cdots\cdots\cdots,$$

$$\Gamma_{i_s}\bigcap\Gamma_{l}\not=\emptyset.$$

\no Thereby we can find a path connecting $g$ and $h$ in
$Cay(\widetilde{\Gamma}:\widetilde{S})$ passing through these
vertices in $Cay(\Gamma_{i_1}:S_{i_1})$,
$Cay(\Gamma_{i_2}:S_{i_2})$, $\cdots$, and
$Cay(\Gamma_{i_s}:S_{i_s})$. Therefore,
$Cay(\widetilde{\Gamma}:\widetilde{S})$ is connected. \quad\quad
$\natural$

The following theorem is gotten by the definition of a Cayley
graph and Theorem $2.4.4$.

\vskip 4mm

\no{\bf Theorem $2.4.5$} \ {\it If $\widetilde{\Gamma} =
\bigcup\limits_{i=1}^n\Gamma$ with $|\Gamma|\geq 3$, then a Cayley
graph $Cay(\widetilde{\Gamma}:\widetilde{S})$}

($i$)\ {\it is an $|\widetilde{S}|$-regular graph;}

($ii$) \ {\it the edge connectivity
$\kappa(Cay(\widetilde{\Gamma}:\widetilde{S}))\geq 2n$.}

\vskip 3mm

{\it Proof} \ The assertion ($i$) is gotten by the definition of
$Cay(\widetilde{\Gamma}:\widetilde{S})$. For ($ii$) since every
Cayley graph of order$\geq 3$ is $2$-connected, for any two
vertices $g,h$ in  $Cay(\widetilde{\Gamma}:\widetilde{S})$, there
are at least $2n$ edge disjoint paths connecting $g$ and $h$.
Whence, the edge connectivity
$\kappa(Cay(\widetilde{\Gamma}:\widetilde{S}))\geq 2n$. \quad\quad
$\natural$

Applying multi-voltage graphs, we get a structure result for
Cayley graphs of a finite multi-group similar to that of Cayley
graphs of a finite group.

\vskip 4mm

\no{\bf Theorem $2.4.6$} \ {\it For a Cayley graph
$Cay(\widetilde{\Gamma}:\widetilde{S})$ of a finite multi-group
$\widetilde{\Gamma}=\bigcup\limits_{i=1}^n\Gamma_i$ with
$\widetilde{S}=\bigcup\limits_{i=1}^nS_i$, there is a
multi-voltage bouquet
$\varsigma:B_{|\widetilde{S}|}\rightarrow\widetilde{S}$ such that
$Cay(\widetilde{\Gamma}:\widetilde{S})\cong
(B_{|\widetilde{S}|})^{\varsigma}$.}

\vskip 3mm

{\it Proof} \ Let $\widetilde{S}=\{s_i ; 1\leq i\leq
|\widetilde{S}|\}$ and $E(B_{|\widetilde{S}|})=\{L_i; 1\leq i\leq
|\widetilde{S}|\}$. Define a multi-voltage graph on a bouquet
$B_{|\widetilde{S}|}$ by

$$\varsigma: L_i\rightarrow s_i, \ \ 1\leq i\leq
|\widetilde{S}|.$$

\no Then we know that there is an isomorphism $\tau$ between
$(B_{|\widetilde{S}|})^{\varsigma}$ and
$Cay(\widetilde{\Gamma}:\widetilde{S})$ by defining $\tau(O_g)=g$
for $\forall g\in\widetilde{\Gamma}$, where
$V(B_{|\widetilde{S}|})=\{O\}$. \quad\quad $\natural$

\vskip 4mm

\no{\bf Corollary $2.4.6$} \ {\it For a Cayley graph
$Cay(\Gamma:S)$ of a finite group $\Gamma$, there exists a voltage
bouquet $\alpha: B_{|S|}\rightarrow S$ such that
$Cay(\Gamma:S)\cong (B_{|S|})^{\alpha}$.}

\vskip 4mm

\no{\bf \S $2.5$ \ Graph Phase Spaces}

\vskip 5mm

\no The behavior of a graph in an $m$-manifold is related with
theoretical physics since it can be viewed as a model of
$p$-branes in M-theory both for a microcosmic and macrocosmic
world. For more details one can see in Chapter $6$. This section
concentrates on surveying some useful fundamental elements for
graphs in $n$-manifolds.

\vskip 4mm

\no{\bf $2.5.1.$ Graph phase in a multi-space}

\vskip 3mm

\no For convenience, we introduce some notations used in this
section in the following. \vskip 3mm

$\widetilde{\bf M}$ -- a multi-manifold $\widetilde{\bf
M}=\bigcup\limits_{i=1}^n{\bf M}^{n_i}$, where ${\bf M}^{n_i}$ is
an $n_i$-manifold, $n_i\geq 2$. For multi-manifolds, see also
those materials in Subsection $1.5.4$.

$\overline{u}\in\widetilde{\bf M}$ -- a point $\overline{u}$ of
$\widetilde{\bf M}$.

${\mathcal G}$ -- a graph $G$ embedded in $\widetilde{\bf M}$.

$C(\widetilde{\bf M})$ -- the set of smooth mappings
$\omega:\widetilde{\bf M}\rightarrow \widetilde{\bf M}$,
differentiable at each point $\overline{u}$ in $\widetilde{\bf
M}$.

\vskip 2mm

Now we define the phase of a graph in a multi-space.

\vskip 4mm

\no{\bf Definition $2.5.1$} \ {\it Let ${\mathcal G}$ be a graph
embedded in a multi-manifold $\widetilde{\bf M}$. A phase of
${\mathcal G}$ in $\widetilde{\bf M}$ is a triple $({\mathcal
G};\omega,\Lambda)$ with an operation $\circ$ on
$C(\widetilde{M})$, where $\omega: V(G)\rightarrow
C(\widetilde{\bf M})$ and $\Lambda: E({\mathcal G})\rightarrow
C(\widetilde{\bf M})$ such that
$\Lambda(\overline{u},\overline{v})=\frac{\omega(\overline{u})\circ\omega(\overline{v})}
{\parallel\overline{u}-\overline{v}\parallel}$ for $\forall
(\overline{u},\overline{v})\in E({\mathcal G})$, where
$\parallel\overline{u}\parallel$ denotes the norm of
$\overline{u}$.}

\vskip 3mm

For examples, the complete graph $K_4$ embedded in ${\bf R}^3$ has
a phase as shown in Fig.$2.33$, where $g\in C({\bf R}^3)$ and
$h\in C({\bf R}^3)$.

\includegraphics[bb=10 10 500 120]{sgm35.eps}

\vskip 3mm

\c{\bf Fig.$2.33$}

\vskip 2mm

Similar to the definition of a adjacent matrix on a graph, we can
also define matrixes on graph phases .

\vskip 4mm

\no{\bf Definition $2.5.2$} \ {\it Let $({\mathcal
G};\omega,\Lambda)$ be a phase and $A[G]=[a_{ij}]_{p\times p}$ the
adjacent matrix of a graph $G$ with $V(G)=\{v_1,v_2,\cdots,v_p\}$.
Define matrixes $V[{\mathcal G}]=[V_{ij}]_{p\times p}$ and
$\Lambda[{\mathcal G}]=[\Lambda_{ij}]_{p\times p}$ by}

$$V_{ij}=\frac{\omega(\overline{v}_{i})}{\parallel\overline{v}_i-\overline{v}_j\parallel}
\ if \ a_{ij}\not=0; \ otherwise, V_{ij}=0$$

\no{\it and}

$$\Lambda_{ij}=\frac{\omega(\overline{v}_i)\circ\omega(\overline{v}_j)}
{\parallel\overline{v}_i-\overline{v}_j\parallel^2} \ if \
a_{ij}\not=0; \ otherwise, \Lambda_{ij}=0,$$

\no{\it where ¡°$\circ$¡± is an operation on $C(\widetilde{M})$.}

\vskip 3mm

For example, for the phase of $K_4$ in Fig.$2.33$, if choice
$g(u)=(x_1,x_2,x_3)$, $g(v)=(y_1,y_2,y_3)$, $g(w)=(z_1,z_2,z_3)$,
$g(o)=(t_1,t_2,t_3)$ and $\circ=\times$, the multiplication of
vectors in ${\bf R}^3$, then we get that

\[
V({\mathcal G})=\left[\begin{array}{cccc} 0 &
\frac{g(u)}{\rho(u,v)} & \frac{g(u)}{\rho(u,w)} &
\frac{g(u)}{\rho(u,o)}\\
\frac{g(v)}{\rho(v,u)} & 0 & \frac{g(v)}{\rho(v,w)} &
\frac{g(v)}{\rho(v,t)}\\
\frac{g(w)}{\rho(w,u)} & \frac{g(w)}{\rho(w,v)} & 0 &
\frac{g(w)}{\rho(w,o)}\\
\frac{g(o)}{\rho(o,u)} & \frac{g(o)}{\rho(o,v)} &
\frac{g(o)}{\rho(o,w)} & 0
\end{array}
\right]
\]

\no where

$$\rho(u,v)=\rho(v,u)=\sqrt{(x_1-y_1)^2+(x_2-y_2)^2+(x_3-y_3)^2},$$

$$\rho(u,w)=\rho(w,u)=\sqrt{(x_1-z_1)^2+(x_2-z_2)^2+(x_3-z_3)^2},$$

$$\rho(u,o)=\rho(o,u)=\sqrt{(x_1-t_1)^2+(x_2-t_2)^2+(x_3-t_3)^2},$$

$$\rho(v,w)=\rho(w,v)=\sqrt{(y_1-z_1)^2+(y_2-z_2)^2+(y_3-z_3)^2},$$

$$\rho(v,o)=\rho(o,v)=\sqrt{(y_1-t_1)^2+(y_2-t_2)^2+(y_3-t_3)^2},$$

$$\rho(w,o)=\rho(o,w)=\sqrt{(z_1-t_1)^2+(z_2-t_2)^2+(z_3-t_3)^2}$$

\no and

\[
\Lambda({\mathcal G})=\left[\begin{array}{cccc} 0 &
\frac{g(u)\times g(v)}{\rho^2(u,v)} & \frac{g(u)\times
g(w)}{\rho^2(u,w)} &
\frac{g(u)\times g(o)}{\rho^2(u,o)}\\
\frac{g(v)\times g(u)}{\rho^2(v,u)} & 0 & \frac{g(v)\times
g(w)}{\rho^2(v,w)} &
\frac{g(v\times g(o)}{\rho^2(v,o)}\\
\frac{g(w)\times g(u)}{\rho^2(w,u)} & \frac{g(w)\times
g(v)}{\rho^2(w,v)} & 0 &
\frac{g(w)\times g(o)}{\rho^2(w,o)}\\
\frac{g(o)\times g(u)}{\rho^2(o,u)} & \frac{g(o)\times
g(v)}{\rho^2(o,v)} & \frac{g(o)\times g(w)}{\rho^2(o,w)} & 0
\end{array}
\right].
\]

\vskip 4mm

\no where

$$g(u)\times g(v)=(x_2y_3-x_3y_2,x_3y_1-x_1y_3,x_1y_2-x_2y_1),$$

$$g(u)\times g(w)=(x_2z_3-x_3z_2,x_3z_1-x_1z_3,x_1z_2-x_2z_1),$$

$$g(u)\times g(o)=(x_2t_3-x_3t_2,x_3t_1-x_1t_3,x_1t_2-x_2t_1),$$

$$g(v)\times g(u)=(y_2x_3-y_3x_2,y_3x_1-y_1x_3,y_1x_2-y_2x_1),$$

$$g(v)\times g(w)=(y_2z_3-y_3z_2,y_3z_1-y_1z_3,y_1z_2-y_2z_1),$$

$$g(v)\times g(o)=(y_2t_3-y_3t_2,y_3t_1-y_1t_3,y_1t_2-y_2t_1),$$

$$g(w)\times g(u)=(z_2x_3-z_3x_2,z_3x_1-z_1x_3,z_1x_2-z_2x_1),$$

$$g(w)\times g(v)=(z_2y_3-z_3y_2,z_3y_1-z_1y_3,z_1y_2-z_2y_1),$$

$$g(w)\times g(o)=(z_2t_3-z_3t_2,z_3t_1-z_1t_3,z_1t_2-z_2t_1),$$

$$g(o)\times g(u)=(t_2x_3-t_3x_2,t_3x_1-t_1x_3,t_1x_2-t_2x_1),$$

$$g(o)\times g(v)=(t_2y_3-t_3y_2,t_3y_1-t_1y_3,t_1y_2-t_2y_1),$$

$$g(o)\times g(w)=(t_2z_3-t_3z_2,t_3z_1-t_1z_3,t_1z_2-t_2z_1).$$

For two given matrixes $A=[a_{ij}]_{p\times p}$ and
$B=[b_{ij}]_{p\times p}$, the {\it star product} ¡°$\ast$¡± on an
operation ¡°$\circ$¡± is defined by $A\ast B=[a_{ij}\circ
b_{ij}]_{p\times p}$. We get the following result for matrixes
$V[{\mathcal G}]$ and $\Lambda[{\mathcal G}]$.

\vskip 4mm

\no{\bf Theorem $2.5.1$} \ \ \ $V[{\mathcal G}]\ast V^t[{\mathcal
G}] = \Lambda[{\mathcal G}]$.

\vskip 3mm

{\it Proof} \ Calculation shows that each $(i,j)$ entry in
$V[{\mathcal G}]\ast V^t[{\mathcal G}]$ is

$$\frac{\omega(\overline{v}_{i})}{\parallel\overline{v}_i-\overline{v}_j\parallel}
\circ\frac{\omega(\overline{v}_{j})}{\parallel\overline{v}_j-\overline{v}_i\parallel}
=\frac{\omega(\overline{v}_{i})\circ\omega(\overline{v}_{j})}
{\parallel\overline{v}_i-\overline{v}_j\parallel^2}=\Lambda_{ij},$$

\no where $1\leq i,j\leq p$. Therefore, we get that

$$V[{\mathcal G}]\ast V^t[{\mathcal G}] = \Lambda[{\mathcal G}]. \quad\quad \natural$$

An operation called {\it addition on graph phases} is defined in
the next.

\vskip 4mm

\no{\bf Definition $2.5.3$} \ {\it For two phase spaces
$({\mathcal G}_1;\omega_1,\Lambda_1)$, $({\mathcal
G}_2;\omega_2,\Lambda_2)$ of graphs $G_1, G_2$ in $\widetilde{M}$
and two operations ¡°$\bullet$¡± and ¡°$\circ$¡± on
$C(\widetilde{M})$, their addition is defined by}

$$({\mathcal
G}_1;\omega_1,\Lambda_1)\bigoplus({\mathcal
G}_2;\omega_2,\Lambda_2)= ({\mathcal G}_1\bigoplus{\mathcal
G}_2;\omega_1\bullet\omega_2,\Lambda_1\bullet\Lambda_2),$$

\no{\it where $\omega_1\bullet\omega_2: V({\mathcal
G}_1\bigcup{\mathcal G}_2)\rightarrow C(\widetilde{M})$
satisfying}

\[ \omega_1\bullet\omega_2(\overline{u}) = \left\{
\begin{array}{ll}
\omega_1(\overline{u})\bullet\omega_2(\overline{u}), & \ if \
\overline{u}\in V({\mathcal G}_1)\bigcap V({\mathcal G}_2),\\
\omega_1(\overline{u}), & \ if \ \overline{u}\in V({\mathcal
G}_1)\setminus V({\mathcal G}_2),\\
\omega_2(\overline{u}), & \ if \ \overline{u}\in V({\mathcal
G}_2)\setminus V({\mathcal G}_1).
\end{array}
\right.
\]

\no{\it and}

$$\Lambda_1\bullet\Lambda_2(\overline{u},\overline{v})
=\frac{\omega_1\bullet\omega_2(\overline{u})\circ\omega_1\bullet\omega_2(\overline{v})}
{\parallel\overline{u}-\overline{v}\parallel^2}$$

\no{\it for $(\overline{u},\overline{v})\in E({\mathcal
G}_1)\bigcup E({\mathcal G}_2)$}

\vskip 3mm

The following result is immediately gotten by Definition $2.5.3$.

\vskip 4mm

\no{\bf Theorem $2.5.2$} \ {\it For two given operations
¡°$\bullet$¡± and ¡°$\circ$¡± on $C(\widetilde{M})$, all graph
phases in $\widetilde{M}$ form a linear space on the field $Z_2$
with a phase $\bigoplus$ for any graph phases $({\mathcal
G}_1;\omega_1,\Lambda_1)$ and $({\mathcal
G}_2;\omega_2,\Lambda_2)$ in $\widetilde{M}$.}

\vskip 4mm

\no{\bf $2.5.2.$ Transformation of a graph phase}

\vskip 4mm

\no{\bf Definition $2.5.4$} \ {\it Let $({\mathcal
G}_1;\omega_1,\Lambda_1)$ and $({\mathcal
G}_2;\omega_2,\Lambda_2)$ be graph phases of graphs $G_1$ and
$G_2$ in a multi-space $\widetilde{M}$ with operations
¡°$\circ_1,\circ_2$¡±, respectively. If there exists a smooth
mapping $\tau\in C(\widetilde{M})$ such that }

$$\tau:({\mathcal
G}_1;\omega_1,\Lambda_1)\rightarrow ({\mathcal
G}_2;\omega_2,\Lambda_2),$$

\no{\it i.e., for $\forall\overline{u}\in V({\mathcal G}_1)$,
$\forall (\overline{u},\overline{v})\in E({\mathcal G}_1)$,
$\tau({\mathcal G}_1)= {\mathcal G}_2$,
$\tau(\omega_1(\overline{u}))=\omega_2(\tau(\overline{u}))$ and
$\tau(\Lambda_1(\overline{u},\overline{v}))=\Lambda_2(\tau(\overline{u},\overline{v}))$,
then we say $({\mathcal G}_1;\omega_1,\Lambda_1)$ and $({\mathcal
G}_2;\omega_2,\Lambda_2)$ are transformable and $\tau$ a transform
mapping.}

\vskip 3mm

For examples, a transform mapping $t$ for embeddings of $K_4$ in
${\bf R}^3$ and on the plane is shown in Fig.$2.34$

\includegraphics[bb=10 10 500 130]{sgm36.eps}

\vskip 3mm

\c{\bf Fig.$2.34$}

\vskip 2mm

\vskip 4mm

\no{\bf Theorem $2.5.3$} \ {\it Let $({\mathcal
G}_1;\omega_1,\Lambda_1)$ and $({\mathcal
G}_2;\omega_2,\Lambda_2)$ be transformable graph phases with
transform mapping $\tau$. If $\tau$ is one-to-one on ${\mathcal
G}_1$ and ${\mathcal G}_2$, then ${\mathcal G}_1$ is isomorphic to
${\mathcal G}_2$.}

\vskip 3mm

{\it Proof} \ By definitions, if $\tau$ is one-to-one on
${\mathcal G}_1$ and ${\mathcal G}_2$, then $\tau$ is an
isomorphism between ${\mathcal G}_1$ and ${\mathcal
G}_2$.\quad\quad $\natural$

A very useful case among transformable graph phases is that one
can find parameters $t_1,t_2,\cdots,t_q, q\geq 1$, such that each
vertex of a graph phase is a smooth mapping of
$t_1,t_2,\cdots,t_q$, i.e., for $\forall
\overline{u}\in\widetilde{M}$, we consider it as
$\overline{u}(t_1,t_2,\cdots,t_q)$. In this case, we introduce two
conceptions on graph phases.

\vskip 4mm

\no{\bf Definition $2.5.5$} \ {\it For a graph phase $({\mathcal
G};\omega,\Lambda)$, define its capacity $Ca({\mathcal
G};\omega,\Lambda)$ and entropy $En({\mathcal G};\omega,\Lambda)$
by }

$$Ca({\mathcal
G};\omega,\Lambda)= \sum\limits_{\overline{u}\in V({\mathcal
G})}\omega(\overline{u})$$

\no{\it and}

$$En({\mathcal G};\omega,\Lambda)=\log(\prod\limits_{\overline{u}\in V({\mathcal
G})}\parallel\omega(\overline{u})\parallel).$$

\vskip 3mm

Then we know the following result.

\vskip 4mm

\no{\bf Theorem $2.5.4$} \ {\it For a graph phase $({\mathcal
G};\omega,\Lambda)$, its capacity $Ca({\mathcal
G};\omega,\Lambda)$ and entropy $En({\mathcal G};\omega,\Lambda)$
satisfy the following differential equations}

$${\rm d}Ca({\mathcal
G};\omega,\Lambda)= \frac{\partial Ca({\mathcal
G};\omega,\Lambda)}{\partial u_i}{\rm d}u_i \ \ and \ \ {\rm
d}En({\mathcal G};\omega,\Lambda)= \frac{\partial En({\mathcal
G};\omega,\Lambda)}{\partial u_i}{\rm d}u_i,$$

\no{\it where we use the Einstein summation convention, i.e., a
sum is over $i$ if it is appearing both in upper and lower
indices.}

\vskip 3mm

{\it Proof} \ Not loss of generality, we assume
$\overline{u}=(u_1,u_2,\cdots,u_p)$ for
$\forall\overline{u}\in\widetilde{M}$. According to the invariance
of differential form, we know that

$${\rm d}\omega = \frac{\partial\omega}{\partial u_i} {\rm
d}u_i.$$

By the definition of the capacity $Ca({\mathcal
G};\omega,\Lambda)$ and entropy $En({\mathcal G};\omega,\Lambda)$
of a graph phase, we get that

\begin{eqnarray*}
{\rm d}Ca({\mathcal G};\omega,\Lambda) &=&
\sum\limits_{\overline{u}\in V({\mathcal G})}{\rm
d}(\omega(\overline{u}))\\
&=& \sum\limits_{\overline{u}\in V({\mathcal
G})}\frac{\partial\omega(\overline{u})}{\partial u_i} {\rm d}u_i=
\frac{\partial(\sum\limits_{\overline{u}\in V({\mathcal
G})}\omega(\overline{u}))}{\partial u_i} {\rm d}u_i\\
&=& \frac{\partial Ca({\mathcal G};\omega,\Lambda)}{\partial
u_i}{\rm d}u_i.
\end{eqnarray*}

Similarly, we also obtain that

\begin{eqnarray*}
{\rm d}En({\mathcal G};\omega,\Lambda) &=&
\sum\limits_{\overline{u}\in V({\mathcal G})}{\rm
d}(\log\parallel\omega(\overline{u})\parallel)\\
&=& \sum\limits_{\overline{u}\in V({\mathcal
G})}\frac{\partial\log|\omega(\overline{u})|}{\partial u_i} {\rm
d}u_i= \frac{\partial(\sum\limits_{\overline{u}\in V({\mathcal
G})}\log\parallel\omega(\overline{u})\parallel)}{\partial u_i} {\rm d}u_i\\
&=& \frac{\partial En({\mathcal G};\omega,\Lambda)}{\partial
u_i}{\rm d}u_i.
\end{eqnarray*}

\no This completes the proof. \quad\quad $\natural$\vskip 3mm

In a $3$-dimensional Euclid space we can get more concrete results
for graph phases $({\mathcal G};\omega,\Lambda)$. In this case, we
get some formulae in the following by choice
$\overline{u}=(x_1,x_2,x_3)$ and $\overline{v}=(y_1,y_2,y_3)$.

$$\omega(\overline{u})= (x_1,x_2,x_3) \ {\rm for} \
\forall\overline{u}\in V({\mathcal G}),$$

$$\Lambda(\overline{u},\overline{v})=
\frac{x_2y_3-x_3y_2,x_3y_1-x_1y_3,x_1y_2-x_2y_1}{(x_1-y_1)^2+(x_2-y_2)^2+(x_3-y_3)^2}
\ {\rm for} \ \forall(\overline{u},\overline{v})\in E({\mathcal
G}),$$

$$Ca({\mathcal
G};\omega,\Lambda)=(\sum\limits_{\overline{u}\in V({\mathcal
G})}x_1(\overline{u}),\sum\limits_{\overline{u}\in V({\mathcal
G})}x_2(\overline{u}),\sum\limits_{\overline{u}\in V({\mathcal
G})}x_3(\overline{u}))$$

\no{and}

$$En({\mathcal
G};\omega,\Lambda)=\sum\limits_{\overline{u}\in V({\mathcal
G})}\log(x_1^2(\overline{u})+x_2^2(\overline{u})+x_3^2(\overline{u}).$$

\vskip 8mm

\no{\bf \S $2.6$ \ Remarks and Open Problems}

\vskip 3mm

\no{\bf $2.6.1$} \ A graphical property $P(G)$ is called to be
{\it subgraph hereditary} if for any subgraph $H\subseteq G$, $H$
posses $P(G)$ whenever $G$ posses the property $P(G)$. For
example, the properties: {\it $G$ is complete} and {\it the vertex
coloring number $\chi(G)\leq k$} both are subgraph hereditary. The
hereditary property of a graph can be generalized by the following
way.

Let $G$ and $H$ be two graphs in a space $\widetilde{M}$. If there
is a smooth mapping $\varsigma$ in $C(\widetilde{M})$ such that
$\varsigma(G)=H$, then we say $G$ and $H$ are {\it equivalent in
$\widetilde{M}$}. Many conceptions in graph theory can be included
in this definition, such as {\it graph homomorphism, graph
equivalent}, $\cdots$, etc.

\vskip 4mm

\no{\bf Problem $2.6.1$} \ {\it Applying different smooth mappings
in a space such as smooth mappings in ${\bf R}^3$ or ${\bf R}^4$
to classify graphs and to find their invariants.}

\vskip 4mm

\no{\bf Problem $2.6.2$} \ {\it Find which parameters already
known in graph theory for a graph is invariant or to find the
smooth mapping in a space on which this parameter is invariant.}

\vskip 3mm

\no{\bf $2.6.2$} \ As an efficient way for finding regular
covering spaces of a graph, voltage graphs have been gotten more
attentions in the past half-century by mathematicians. Works for
regular covering spaces of a graph can seen in $[23]$, $[45]-[46]$
and $[71]-[72]$. But few works are found in publication for
irregular covering spaces of a graph. The multi-voltage graph of
type $1$ or type $2$ with multi-groups defined in Section $2.2$
are candidate for further research on irregular covering spaces of
graphs.

\vskip 4mm

\no{\bf Problem $2.6.3$} \ {\it Applying multi-voltage graphs to
get the genus of a graph with less symmetries.}

\vskip 4mm

\no{\bf Problem $2.6.4$} \ {\it Find new actions of a multi-group
on a graph, such as the left subaction and its contribution to
topological graph theory. What can we say for automorphisms of the
lifting of a multi-voltage graph?}

\vskip 3mm

There is a famous conjecture for Cayley graphs of a finite group
in algebraic graph theory, i.e., {\it every connected Cayley graph
of order$\geq 3$ is hamiltonian}. Similarly, we can also present a
conjecture for Cayley graphs of a multi-group.

\vskip 4mm

\no{\bf Conjecture $2.6.1$} \ {\it Every Cayley graph of a finite
multi-group $\widetilde{\Gamma}=\bigcup\limits_{i=1}^n\Gamma_i$
with order$\geq 3$ and $|\bigcap\limits_{i=1}^n\Gamma_i|\geq 2$ is
hamiltonian.}

\vskip 3mm

\no{\bf $2.6.3$} As pointed out in $[56]$, for applying
combinatorics to other sciences, a good idea is pullback measures
on combinatorial objects, initially ignored by the classical
combinatorics and reconstructed or make a combinatorial
generalization for the classical mathematics, such as, the
algebra, the differential geometry, the Riemann geometry, $\cdots$
and the mechanics, the theoretical physics, $\cdots$. For this
object, a more natural way is to put a graph in a metric space and
find its good behaviors. The problem discussed in Sections $2.3$
is just an elementary step for this target. More works should be
done and more techniques should be designed. The following open
problems are valuable to research for a researcher on
combinatorics.

\vskip 4mm

\no{\bf Problem $2.6.5$} \ {\it Find which parameters for a graph
can be used to a graph in a space. Determine combinatorial
properties of a graph in a space.}

\vskip 3mm

Consider a graph in an Euclid space of dimension $3$. All of its
edges are seen as a structural member, such as steel bars or rods
and its vertices are hinged points. Then we raise the following
problem.

\vskip 4mm

\no{\bf Problem $2.6.6$} \ {\it Applying structural mechanics to
classify what kind of graph structures are stable or unstable.
Whether can we discover structural mechanics of dimension$\geq 4$
by this idea?}

\vskip 3mm

We have known the orbit of a point under an action of a group, for
example, a torus is an orbit of $Z\times Z$ action on a point in
${\bf R}^3$. Similarly, we can also define an {\it orbit of a
graph in a space} under an action on this space.

{\it Let ${\mathcal G}$ be a graph in a multi-space
$\widetilde{M}$ and $\Pi$ a family of actions on $\widetilde{M}$.
Define an orbit $Or({\mathcal G})$ by}

$$Or({\mathcal
G})=\{\pi({\mathcal G})| \ \forall\pi\in\Pi\}.$$

\vskip 4mm

\no{\bf Problem $2.6.7$} \ {\it Given an action $\pi$, continuous
or discontinuous on a space $\widetilde{M}$, for example ${\bf
R}^3$ and a graph ${\mathcal G}$ in $\widetilde{M}$, find the
orbit of ${\mathcal G}$ under the action of $\pi$. When can we get
a closed geometrical object by this action?}

\vskip 4mm

\no{\bf Problem $2.6.8$} \ {\it Given a family ${\mathcal A}$ of
actions, continuous or discontinuous on a space $\widetilde{M}$
and a graph ${\mathcal G}$ in $\widetilde{M}$, find the orbit of
${\mathcal G}$ under these actions in ${\mathcal A}$. Find the
orbit of a vertex or an edge of ${\mathcal G}$ under the action of
${\mathcal G}$, and when are they closed?}

\vskip 3mm

\no{\bf $2.6.4$} \ The central idea in Section $2.4$ is that a
graph is equivalent to Smarandache multi-spaces. This fact enables
us to investigate Smarandache multi-spaces possible by a
combinatorial approach. Applying infinite graph theory (see $[94]$
for details), we can also define an infinite graph for an infinite
Smarandache multi-space similar to Definition $2.4.3$.

\vskip 4mm

\no{\bf Problem $2.6.9$} \ {\it Find its structural properties of
an infinite graph of an infinite Smarandache multi-space.}

\vskip 3mm

\no{\bf $2.6.5$} \ There is an alternative way for defining
transformable graph phases, i.e., by homotopy groups in a
topological space, which is stated as follows.

Let $({\mathcal G}_1;\omega_1,\Lambda_1)$ and $({\mathcal
G}_2;\omega_2,\Lambda_2)$ be two graph phases. If there is a
continuous mapping $H:C(\widetilde{M})\times I\rightarrow
C(\widetilde{M})\times I$, $I=[0,1]$ such that
$H(C(\widetilde{M}),0)=({\mathcal G}_1;\omega_1,\Lambda_1)$ and
$H(C(\widetilde{M}),1)=({\mathcal G}_2;\omega_2,\Lambda_2)$, then
$({\mathcal G}_1;\omega_1,\Lambda_1)$ and $({\mathcal
G}_2;\omega_2,\Lambda_2)$ are said two {\it transformable graph
phases}.

Similar to topology, we can also introduce product on homotopy
equivalence classes and prove that all homotopy equivalence
classes form a group. This group is called a {\it fundamental
group} and denote it by $\pi({\mathcal G};\omega,\Lambda)$. In
topology there is a famous theorem, called the {\it Seifert and
Van Kampen theorem} for characterizing fundamental groups
$\pi_1({\mathcal A})$ of topological spaces ${\mathcal A}$
restated as follows (see $[92]$ for details).

\vskip 3mm

{\it Suppose ${\mathcal E}$ is a space which can be expressed as
the union of path-connected open sets ${\mathcal A}$, ${\mathcal
B}$ such that ${\mathcal A}\bigcap{\mathcal B}$ is path-connected
and $\pi_1({\mathcal A})$ and $\pi_1({\mathcal B})$ have
respective presentations}

$$\left<a_1,\cdots,a_m; r_1,\cdots,r_n\right>,$$

$$\left<b_1,\cdots,b_m; s_1,\cdots,s_n\right>$$

\no{\it while $\pi_1({\mathcal A}\bigcap{\mathcal B}$) is finitely
generated. Then $\pi_1({\mathcal E}$) has a presentation}

$$\left<a_1,\cdots,a_m,b_1,\cdots,b_m; r_1,\cdots,r_n,s_1,\cdots,s_n,
u_1=v_1,\cdots, u_t=v_t\right>,$$

\no{\it where $u_i,v_i, i=1,\cdots, t$ are expressions for the
generators of $\pi_1({\mathcal A}\bigcap{\mathcal B}$) in terms of
the generators of $\pi_1({\mathcal A})$ and $\pi_1({\mathcal B})$
respectively.}

\vskip 2mm

Then there is a problem for the fundamental group $\pi({\mathcal
G};\omega,\Lambda)$ of a graph phase $({\mathcal
G};\omega,\Lambda)$.

\vskip 4mm

\no{\bf Problem $2.6.10$}\ {\it Find a result similar to the
Seifert and Van Kampen theorem for the fundamental group of a
graph phase.}

\end{document}